\documentclass[final]{siamltex}

\usepackage{latexsym, amssymb, enumerate, amsmath, hyperref, pdftricks, framed,color,enumitem,tikz}

\renewcommand{\theequation}{\arabic{section}.\arabic{equation}}

\pdfoutput=1

\def\Y_#1{y_{\!#1}}

\def\E{\mathbb{E}}
\def\P{\mathbb{P}}

\def\TV{\mathrm{TV}}

\def\tr{\operatorname{trace}}
\def\divergence{\operatorname{div}}
\newcommand{\bfi}[1]{\textbf{\textit{#1}}}
\newcommand{\argmin}{\arg\!\min}

\newtheorem{thm}{Theorem}[section]
\newtheorem{prop}[thm]{Proposition}

\newtheorem{cor}[thm]{Corollary}

\newtheorem{assumption}{Assumption}[section]

\newtheorem{algorithm}[theorem]{Algorithm}

\definecolor{shadecolor}{gray}{0.96}  

\title{Metropolis integration schemes\\ for self-adjoint diffusions}
\author{Nawaf Bou-Rabee\thanks {Department of Mathematical Sciences, 
Rutgers University Camden, 311 N 5th Street, Camden, NJ 08102, USA, (nawaf.bourabee@rutgers.edu).
The work of N.~B-R. was supported by NSF grant DMS-1212058.}
\and Aleksandar Donev\thanks{Courant Institute of Mathematical Sciences,
New York University, 251 Mercer Street, New York, NY 10012-1185, (donev@courant.nyu.edu).
A. Donev was supported in part by the Air Force Office of Scientific 
Research under grant number FA9550-12-1-0356.
}
\and Eric Vanden-Eijnden\thanks {Courant Institute of Mathematical Sciences,
New York University, 251 Mercer Street, New York, NY 10012-1185, (eve2@cims.nyu.edu).}
}

\begin{document}
\maketitle

\begin{abstract} 
We present explicit methods for simulating diffusions whose generator is self-adjoint with respect to a known (but possibly not normalizable) density. These methods exploit this property and combine an optimized Runge-Kutta algorithm with a Metropolis-Hastings Monte-Carlo scheme. The resulting numerical integration scheme is shown to be weakly accurate at finite noise and to gain higher order accuracy in the small noise limit. It also permits to avoid computing explicitly certain terms in the equation, such as the divergence of the mobility tensor, which can be tedious to calculate.  Finally, the scheme is shown to be ergodic with respect to the exact equilibrium probability distribution of the diffusion when it exists. These results are illustrated on several examples including a Brownian dynamics simulation of DNA in a solvent.  In this example, the proposed scheme  is able to accurately compute dynamics at time step sizes that are an order of magnitude (or more) larger than those permitted with commonly used explicit predictor-corrector schemes.  \end{abstract}



%
%

\pagestyle{myheadings}
\thispagestyle{plain}
\markboth{BOU-RABEE, DONEV, and VANDEN-EIJNDEN}{Metropolis integrators for self-adjoint diffusions}

\section{Introduction}\label{intro}

This paper is about simulation of diffusions of the type: 
\begin{equation} 
  \label{sde}
 \boxed{
  d Y = \underset{\text{deterministic drift}}{\underbrace{- M( Y )DU( Y ) dt}} + \underset{\text{heat bath}}{\underbrace{\beta^{-1} \operatorname{div} M( Y ) dt + \sqrt{2 \beta^{-1}} B( Y ) d W}}
  }
\end{equation}
where we have introduced the following.

\bigskip

\begin{center}
\begin{tabular}{cc}
\hline
$Y(t) \in \mathbb{R}^{n}$ &  state of the system \\
$U: \mathbb{R}^n \to \mathbb{R}$  &  potential energy \\
$DU: \mathbb{R}^n \to \mathbb{R}^n$&   force  \\
$M : \mathbb{R}^n \to \mathbb{R}^{n \times n}$ & mobility  matrix \\
$(\operatorname{div} M)_i = \sum_{j=1}^n \partial M_{i j}/\partial x_j$ & divergence of mobility matrix \\
$B : \mathbb{R}^n \to \mathbb{R}^{n \times n}$ &  noise coefficient matrix \\
$W(t) \in \mathbb{R}^n$ &   Brownian motion  \\
$\beta$ &  inverse temperature factor \\
\hline
\end{tabular}
\end{center}

\bigskip

\noindent
Let $B(x)^T$ denote the transpose of the real matrix $B(x)$.  Assuming that \begin{equation}
M(x) = B(x) B(x)^T  \qquad \text{for all $x \in \mathbb{R}^n$} 
\end{equation}
this dynamics defines a Markov process whose generator $L$ is {\em self-adjoint} with respect to the following density
\begin{equation}  
  \label{eq:nu}
  \nu(x) = \exp\left(- \beta U(x) \right) \;.
\end{equation}
Indeed since the action of $L$ on any suitable test function $f(x)$ can be written as
\begin{equation}
  \label{eq:generator}
  (Lf)(x) = \beta^{-1}\, \nu^{-1}(x) \operatorname{div} \left( \nu(x) M( x ) D f(x) \right) 
\end{equation}
it follows from integration by parts that
\begin{equation}
\label{eq:idb}
\langle Lf, g \rangle_{\nu} = \langle f, L g \rangle_{\nu} \qquad \text {for all suitable test functions $f(x)$ and $g(x)$}  
\end{equation}
where $\langle \cdot , \cdot  \rangle_{\nu}$ denotes an $L^2$-inner product weighted by the density $\nu(x)$: \[
\langle f , g  \rangle_{\nu} = \int_{\mathbb{R}^n} f(x) g(x) \nu(x) dx \;.
\]   This property implies that the diffusion is  {\em $\nu$-symmetric}~\cite{Ke1978}, in the sense that
\begin{equation}
  \label{eq:db}
  \nu(x) p_t(x,y) = \nu(y) p_t(y,x) \qquad \text {for all $t>0$} 
\end{equation}
where $p_t(x,y)$ denotes the transition density of $Y(t)$ given that $Y(0) = x$, i.e., \[
\int_A p_t(x,y) dy = \P ( Y(t) \in A \mid Y(0) = x ) 
\]  for any measurable set $A$. In fact, \eqref{eq:idb} corresponds to an 
infinitesimal version of  \eqref{eq:db}.  We stress that \eqref{eq:db} can hold even if $\nu(x)$ is not normalizable: in this 
case, the solution to \eqref{sde} reaches no statistical steady state, i.e., it is not ergodic.  If, however, the density~\eqref{eq:nu} 
is normalizable, 
\begin{equation}
  \label{eq:normalize}
  Z = \int_{\mathbb{R}^n} \exp\left(- \beta U(x) \right) dx < \infty 
\end{equation}
%
then the diffusion \eqref{sde} is ergodic with respect to $\nu(x)/Z$.   
This normalized density is called the equilibrium probability density of the diffusion and is also referred to as the \textit{Boltzmann-Gibbs density}.
When \eqref{eq:normalize} holds, the dynamics observed at equilibrium is time-reversible, and the property~\eqref{eq:db} is referred to as the {\em detailed balance condition} in the physics literature.  Our main purpose here is to exploit the properties above, in particular~\eqref{eq:db}, to design stable and accurate numerical integrators for~\eqref{sde} that work even in situations when the density $\nu(x)$ in~\eqref{eq:nu} is not normalizable. 

The importance of this objective stems from the wide range of applications where diffusions of the type~\eqref{sde} serve as dynamical models. For example, \eqref{sde} is used to model the evolution of coarse-grained molecular systems in regimes where the details of the microscopic interactions can be represented by a heat bath.  Multiscale models of polymers in a Stokesian solvent fit this category, a context in which \eqref{sde} is referred to as Brownian dynamics (BD) with hydrodynamic interactions \cite{ErMc1978, Fi1978, Ot1996, weinan2011principles}.  BD has been used to quantitatively simulate the non-equilibrium dynamics of biopolymers in dilute solutions \cite{La2005, Sh2005}, and to validate and clarify experimental findings for bacteriophage DNA dynamics in extensional, shear, and mixed planar flows \cite{LaHuSmCh1999, HuShLa2000, HuShBaCh2002, JedeGr2002, ScBaShCh2003, ScShCh2004}.  The dynamics of DNA molecules in confined solutions with complex geometries has also been studied using BD to guide the design of microfluidic devices to manipulate these molecules \cite{ChLa2002, JeDiScGrde2003a,JeDiScGrde2003b,JeScdeGr2004, WoShKh2004a, WoShKh2004b,HeMadeGr2006a, HeMadeGr2006b,HedeGr2007,HeMadeGr2008, Gr2011, de2011, ZhdeGr2012}.   To be clear, the proposed integrator is not  able to simulate diffusions that  do not satisfy the $\nu$-symmetry condition, e.g., polymers in shear flows.

Because of the relevance of~\eqref{sde} to BD, this is also the context in which most work on the development of numerical schemes to simulate this equation has been devoted~\cite{Ot1996,HuOt1998}. Perhaps the most widely used among these is the so-called Fixman scheme, which is an explicit predictor-corrector scheme that also has the advantage that it avoids the explicit computation of the divergence of the mobility matrix~\cite{Fi1978}. Unfortunately, explicit BD time integrators like the Fixman scheme are often insufficient in applications because the interplay between the drift and the noise terms in~\eqref{sde} can induce numerical instabilities or artifacts such as systematic drifts \cite{CaGa1991, ReFrGa1992}.  This is a well known problem with explicit discretizations of nonlinear diffusions \cite{Ta2002, HiMaSt2002, MiTr2005,Hi2011,HuJeKl2012}, and it is especially acute if the potential force in~\eqref{sde} is stiff, which is typically the case in applications.  The standard solution to this problem is to introduce implicitness.  In a sequence of works, this strategy was investigated to deal with the steep potentials that enforce finite extensibility of bond lengths in bead-spring models \cite{Fi1986, HeOt1997, JedeGr2000, SoKhWoHuSh2002, HsLiLa2003}.  The leading semi-implicit predictor-corrector scheme emerging from this effort is numerically stable at much larger time step sizes than explicit predictor-corrector schemes but it requires a nonlinear solve at every step, which is time-consuming and raises convergence questions about the methods used to solve this nonlinear system that are not easy to answer, in general.  

In this paper, following the strategy introduced in~\cite{BoVa2010}, we adopt a probabilistic approach based on the $\nu$-symmetry property in~\eqref{eq:db} to build stable explicit integrators for~\eqref{sde} that avoid implicitness and do not assume a specific form of the potential forces. Specifically, we propose an integrator for~\eqref{sde} with the following properties.

\medskip

\begin{description}
\item[(P1)] it samples the exact equilibrium probability density of the SDE when this density exists (i.e., when $\nu(x)$ in~\eqref{eq:nu} is normalizable);
\item[(P2)] it generates a weakly accurate approximation to the solution of~\eqref{sde} at constant temperature; 
\item[(P3)] it acquires higher order accuracy in the small noise limit, $\beta \to \infty$; and,
\item[(P4)] it avoids computing the divergence of the mobility matrix. 
\end{description}

\medskip

\noindent
We stress that Properties (P2)-(P4) hold even when the density $\nu(x)$ in~\eqref{eq:nu} is not normalizable, i.e., \eqref{sde} is out of equilibrium and not ergodic with respect to any invariant distribution. Thus, even though the scheme we propose involves a Monte-Carlo step, its aim and philosophy are very different from Monte-Carlo methods where a discretized version of~\eqref{sde} is used as a proposal step but whose only goal is to sample a target distribution with no concern for the dynamics of~\eqref{sde} \cite{MeRoRoTeTe1953,Ha1970, RoDoFr1978, DuKePeRo1987, Ho1991, KePe2001, Li2008, AkRe2008, AkBoRe2009, LeRoSt2010A}. Compared to the method proposed in~\cite{BoVa2010}, the main novelty of the scheme introduced here stems from Properties (P3) and (P4).  The point of  (P3) is that the finite time dynamics of $Y(t)$ in the small noise limit, e.g., how the solution moves to the critical points of $U(x)$, is a limiting situation that the approximation should be able to handle. To achieve (P3) we need to control the rate of rejections in the Monte Carlo step in the small noise limit which is nontrivial because the density~$\nu(x)$ in~\eqref{eq:nu} that enters the key relation~\eqref{eq:db} becomes singular in this limit.    Similarly, (P4) constrains the type of proposal moves we can use in the scheme, and this property is important in situations where the mobility matrix does not have an explicit expression, e.g., when the solvent is confined \cite{Gr2011}.  In the sequel, we use a BD model of DNA in a solvent to show that the new scheme achieves the same accuracy as explicit predictor-corrector schemes with time step sizes that are an order of magnitude (or more) larger.  Beside BD, this scheme should also be useful 
in the simulation of polymer conformational transitions in electrophoretic flows \cite{kim2006brownian,randall2006methods,kim2007design,trahan2009simulation,hsieh2011simulation,hsieh2012simulation}, in other contexts where the effective dynamics of a set of coarse-grained variables satisfies~\eqref{sde} \cite{EVa2004, MaVa2006, MaFiVaCi2006, EVa2010, LeLe2010}, in applications to Bayesian inference \cite{GiCa2011}, etc.

\subsection*{Organization of the paper}  The new scheme, with a structure that immediately shows that it satisfies Property (P4) is introduced in \S\ref{mainresults}.  In \S\ref{sec:numerics}, we present numerical examples with comparisons to the Fixman algorithm.  A proof of the ergodicity Property (P1) of the scheme is provided in \S\ref{sec:ergodicity}.  The weak accuracy Property (P2) is proven in \S\ref{sec:accuracy} using~\eqref{eq:db}.  The behavior of the scheme in the small noise limit, specifically Property (P3), is investigated in \S\ref{sec:deterministicaccuracy}.  A conclusion is given in \S\ref{sec:conclusion}.

\section{The integrator and its main properties}\label{mainresults}

Following~\cite{BoVa2010}, the scheme introduced in this paper combines a proposal step obtained via time-discretization of~\eqref{sde} with an accept/reject Monte-Carlo step. 
The detailed algorithm is given below in terms of vector and matrix-valued variables $G_h : \mathbb{R}^{n} \to \mathbb{R}^{n}$ and $B_h: \mathbb{R}^{n} \to \mathbb{R}^{n \times n}$, respectively, whose explicit forms will be specified shortly in such a way that Properties~(P1)--(P4) are met.   

\medskip


\begin{shaded}
\begin{algorithm}[Metropolis Integrator] \label{MetropolisIntegrator}
Given the current position $X_0$ at time $t$ the algorithm proposes 
an updated position $X_1^{\star}$ at time $t+h$ for some time step size 
$h>0$ via
\begin{equation}
\label{proposal} 
\begin{cases}
X_1^{\star} = \tilde X_{1} + h\, G_h ( \tilde X_1)  + ( \tilde X_1 - X_0 )  \\
\tilde X_1 =  X_0 + \sqrt{\frac{h}{2}} B_h(X_0) \xi 
\end{cases}
\end{equation}
Here $\xi\in\mathbb{R}^n$ denotes a Gaussian random vector with 
mean zero and covariance $\E( \xi_i \xi_j ) = \beta^{-1} \delta_{ij}$.  The 
\bfi{proposal move} $X_1^{\star}$ is then accepted or rejected by 
taking as actual update for the position at time $t+h$ the value
\begin{equation}  \label{actualupdate}
  X_1 = \gamma X_1^{\star} + (1-\gamma) X_0 
\end{equation}
where $\gamma$ is a Bernoulli random variable which 
assumes the value 1 with probability $\alpha_h(X_0, \xi)$ 
and value 0 with probability $1-\alpha_h(X_0, \xi)$.
The function $\alpha_h(X_0, \xi)$ is known as the 
\bfi{acceptance probability} and is given by
\begin{equation}
  \label{alphah}
  \begin{aligned}
    & \alpha_h(X_0, \xi)  \\
    & \quad = \min\left(1, \frac{\det(B_h(X_0))}{\det(B_h(X_1^{\star}))}
      \exp\left(- \beta \left[ \frac{| \eta |^2}{2} 
          - \frac{ | \xi |^2}{2} + U(X_1^{\star}) - U(X_0) \right] \right) \right)  
  \end{aligned} 
\end{equation} 
where $\eta$ satisfies: 
\begin{equation}
 \label{eq:linsolv}
 B_h(X_1^{\star}) \eta = B_h(X_0) \xi + \sqrt{2 h}\, G_h(\tilde X_1)  
\end{equation}
\end{algorithm}
\end{shaded}


\medskip

\begin{figure}
\centering
\begin{tikzpicture}[scale=1.5]
\draw[dotted]    (0,1.5) -- (2,0) -- (3.33,1.33);
\filldraw[color=black,fill=white] (0,1.5) circle (0.05);	
\filldraw[color=black,fill=white] (2,0) circle (0.05);	
\filldraw[color=black,fill=white] (3.33,1.33) circle (0.05);	
\draw[->, ultra thick](0,1.5) -- (0.25,1.3125);
\draw[->, ultra thick](2,0) -- (2.25,0.25);
\node[black,scale=1.5] at (0.1,1.80) {$X_0$};
\node[black,scale=1.5] at (3.3,-0.3) {$\tilde X_1  {\scriptstyle  = X_0 + \sqrt{\frac{h}{2}} B_h(X_0) \xi }$};
\node[black,scale=1.5] at (4.9,1.55) {$X_1^{\star}  {\scriptstyle = \tilde X_1+ \sqrt{\frac{h}{2}} B_h(X_1^{\star}) \eta} $};
\node[black, scale=1.25,fill=white] at (0.54,1.05) {$B_h(X_0) \xi$};
\node[black, scale=1.25,fill=white] at (2.5,0.5) {$B_h({X_1}^{\star}) \eta$};
\end{tikzpicture}
\caption{A diagram depicting the points $X_0$, $\tilde X_1$ and $X_1^{\star}$ appearing in Algorithm~\ref{MetropolisIntegrator},
and how they are related by the vectors $B_h(X_0) \xi$ and $B_h(X_1^{\star}) \eta$ introduced in \eqref{proposal} and \eqref{eq:linsolv}, respectively.}
\label{proposal_move_diagram}
\end{figure}
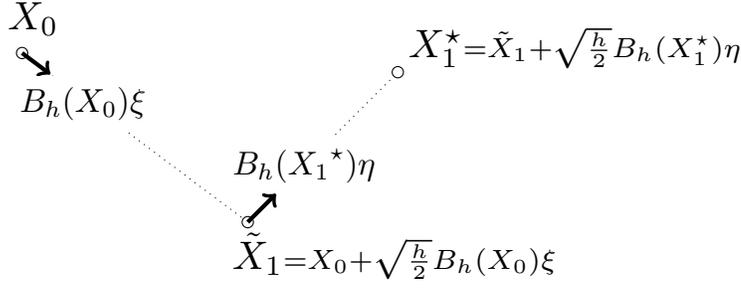

We refer to Algorithm~\ref{MetropolisIntegrator} as Metropolis integrator in the rest of the paper.
Figure~\ref{proposal_move_diagram} illustrates how the proposal move $X_1^{\star}$ in \eqref{proposal} is related to the initial state $X_0$ and the Gaussian random vector $\tilde X_1$.  We remark that the calculation of the acceptance probability in \eqref{alphah} usually requires a Cholesky factorization to determine $B_h(X_1^{\star})$ and its determinant; and a linear solve to determine~$\eta$ via~\eqref{eq:linsolv}.  If $B_h(x) B_h(x)^T$ is positive definite -- which is the case in all of the examples considered in this paper -- this linear system has a unique solution.   Let us now discuss how to choose $G_h(x)$ and $B_h(x)$ in~\eqref{proposal} by requiring that Properties~(P1)--(P4) are satisfied. It is useful to look at these properties sequentially and see the constraints on $G_h(x)$ and $B_h(x)$ they entail.

\subsection{Property (P1): Ergodicity} 
Under quite general assumptions on $\nu(x)$, $G_h(x)$ and $B_h(x)$ (see Assumptions~\ref{P1_nu_assumption},~\ref{P1_Gh_assumption}, and \ref{P1_Bh_assumption}),
it will be shown in \S\ref{sec:ergodicity} that Algorithm~\ref{MetropolisIntegrator} is ergodic with respect to the equilibrium probability distribution with density $\nu(x)/Z$.  Thus, one can stably generate an infinitely long trajectory from the method with the right stationary distribution.  Moreover, the ergodic theorem implies that for every initial state and sufficiently small time step size:
\begin{equation}
  \label{eq:ergo1}
 \frac{1}{T} \int_0^T f(X_{ \lfloor t/h \rfloor } ) dt \to \frac1Z\int_{\mathbb{R}^n} f(x) \nu(x) dx \;, \quad \text{as $T \to \infty$}  \;, \quad \text{a.s.}
\end{equation}
where $f(x)$ is any suitable test function \cite{mattingly2010convergence}. Note that the acceptance probability function in \eqref{alphah} does not involve derivatives of either $G_h(x)$ or $B_h(x)$ because of a symmetry in the law of the proposal move.  Specifically, given the Gaussian random vector $\tilde X_1$, the proposal move satisfies the equation: \[
X_1^{\star} + X_0 = 2 \tilde X_1 + h G_h( \tilde X_1 ) 
\] which is invariant under the step reversal transformation $(X_0, X_1^{\star}) \mapsto (X_1^{\star}, X_0)$.   This invariance is a key step in the proof of Theorem~\ref{thm:ergodicityofscheme}, where the $\nu$-symmetry and ergodicity of the Metropolis Integrator are proved.  It is not surprising that Property (P1) holds for quite general $G_h(x)$ and $B_h(x)$ (for example $G_h(x) = 0$) since the $\nu$-symmetry property~\eqref{eq:db} is enforced by the accept/reject step in~Algorithm~\ref{MetropolisIntegrator} as long as the correct acceptance probability $\alpha_h(X_0,\xi)$ is used.   (Note that the transition probability distribution of  Algorithm~\ref{MetropolisIntegrator} has no density with respect to Lebesgue measure, so~\eqref{eq:db} must be reinterpreted in terms of distributions: this is a minor technical difficulty on which we will dwell upon in \S\ref{sec:ergodicity}.)  This idea is the essence of the Metropolis-Hastings method. More remarkable is the next observation, relevant to Property (P2).

\subsection{Property (P2): Weak Accuracy} 
For every sufficiently regular $G_h(x)$ and $B_h(x)$ (see Assumptions~\ref{P2_Gh_assumption} and \ref{P2_Bh_assumption}) satisfying 
\begin{equation} \label{P2condition}
B_h(x)B_h(x)^T = M(x) + \mathcal{O}(h) \;, \qquad \text{for all $x \in \mathbb{R}^n$} 
\end{equation}
where $M(x)$ is the mobility matrix entering~\eqref{sde}, Algorithm~\ref{MetropolisIntegrator}  is weakly accurate on finite time intervals:
\begin{equation}
  \label{weakacc}
| \E_x ( f(Y( \lfloor t/h \rfloor h))) - \E_x (f (X_{\lfloor t/h \rfloor }) ) | \le C(T, G_h) h^{1/2} 
\end{equation}
for every $t \in [0, T]$, $x \in \mathbb{R}^n$, and sufficiently regular $U(x)$ (see Assumption~\ref{P2_U_Assumption}).
Here $\E_x$ denotes the expectation conditional on the initial state being $x$. The observation above is remarkable in that it holds for any sufficiently regular $G_h(x)$, including the trivial $G_h(x)=0$. As we will see in \S\ref{sec:accuracy}, this really is a consequence of an (infinitesimal) fluctuation-dissipation theorem: as long as the noise term is handled accurately in the proposal move~\eqref{proposal} (which it is if \eqref{P2condition} holds), weak accuracy follows from the fact that the only diffusion satisfying the $\nu$-symmetry property~\eqref{eq:db} (which it does through Property (P1)) is the one with the correct drift term. Note that this immediately opens the door to schemes where the divergence of $M(x)$ does not need to be computed.  In fact, with $G_h(x)=0$, the calculation of no part of the drift is necessary. Of course, this trivial choice is not the best (nor even a good) one in terms of accuracy, and to enforce Property~(P3), we will have to use more specific $G_h(x)$ and $B_h(x)$.

\subsection{Property (P3): Second-order Accuracy in Small Noise Limit}   

Let $G_h(x)$ and $B_h(x)$ be the following two-stage Runge-Kutta combinations:
\begin{equation} \label{rk2_Gh}
\begin{cases}
G_h(x) = - b_1 M(x) DU(x)  - b_2 M(x) DU(x_1) \\
\qquad \qquad - b_3 M(x_1) DU(x) - b_4 M(x_1) DU(x_1)  \\
 x_1 = x - a_{12} h M(x) DU(x)   \\
 \end{cases}
 \end{equation}
\begin{equation} \label{rk2_Bh}
\begin{cases}
B_h(x) B_h(x)^T = d_1 M(x) + d_2 M(\bar x_1)  \\
\bar x_1 =  x + c_{12} h M(x) DU(x) 
\end{cases}
\end{equation}
with parameter values: 
\begin{equation} \label{rk2_optimal_parameters}
d_1 = 1/4 \;, ~~  d_2 = 3/4 \;,~~  b_1=5/8 \;,~~  b_2= b_3=-3/8 \;,~~  b_4=9/8\;, ~~  \&~~ c_{12}=a_{12}=2/3 \;.
\end{equation}
Algorithm~\ref{MetropolisIntegrator} operated with this choice of $G_h(x)$ and $B_h(x)$ satisfies:
\begin{equation}
  \label{eq:strong_acc}
\lim_{\beta \to \infty}  ( \E_x | Y( \lfloor t/h \rfloor h) - X_{\lfloor t/h \rfloor } |^2 )^{1/2} \le C(T) h^2
\end{equation}
for every $t \in [0, T]$ and $x \in \mathbb{R}^n$ and sufficiently regular $U(x)$ and $M(x)$ 
(see Assumptions~\ref{P3_U_Assumption} and \ref{P3_M_Assumption}).

\medskip

When the mobility matrix is constant, then $B_h(x) = B(x)$ and $G_h(x)$ reduces to the so-called Ralston Runge-Kutta method \cite{ralston1962runge}.  The above statement, which is one of the main results of this paper, is established in \S\ref{sec:deterministicaccuracy} using an asymptotic analysis of the acceptance probability function in the deterministic limit as $\beta \to \infty$.   This analysis reveals that the parameter values in \eqref{rk2_optimal_parameters} are optimal, that is, they are the only choice that yield Property (P3).  The validity of this statement requires assumptions whose sufficiency is discussed in the numerical examples in \S\ref{sec:numerics}.  When these conditions are violated, \S\ref{second_order} proposes a modification to the Metropolis integrator which involves adapting the parameter $a_{12}$ appearing in \eqref{rk2_optimal_parameters}.   The small noise limit is relevant when the deterministic drift dominates the dynamics of \eqref{sde}, which is the case in BD applications when the P\'{e}clet number is moderate to high or when the system is driven by an external flow.  (The P\'{e}clet number compares the work done by the potential force to the thermal energy $\beta^{-1}$.)   As a by-product of (P3), in \S\ref{sec:accuracy} we show that the Metropolis integrator is $3/2$-weakly accurate when the mobility matrix is constant.

\subsection{Additional remarks on integrator}  
Note that Property~(P4) automatically follows when $G_h(x)$ and $B_h(x)$ are given by \eqref{rk2_Gh} and~\eqref{rk2_Bh}, respectively, since at no point do we need to calculate the divergence of $M(x)$.  Note also, that with this choice of $G_h(x)$ and $B_h(x)$, the proposal move in Algorithm~\ref{MetropolisIntegrator} involves internal stages.  If an internal stage variable is not in the domain of definition of the mobility matrix or force, e.g., assumes a non-physical value, it is straightforward to show that any extension of these functions results in an algorithm that satisfies the $\nu$-symmetry condition~\eqref{eq:db}.   Hence, we suggest using the trivial extension where the mobility matrix is set equal to the identity matrix and the  force is set equal to zero.   Independent of the extension chosen, the energy is taken to be infinite at non-physical states, which from \eqref{alphah} implies that non-physical proposal moves are always  rejected.

Let us end this section by stressing that the idea of using Monte-Carlo to perform BD simulation is not new and goes back at least to \cite{KiYoMaWa1991, KiYoMaWa1992}.  In place of the proposal move \eqref{proposal}, these papers use the Ermak-McCammon scheme, which corresponds to a forward Euler discretization of \eqref{sde}. This scheme  reduces to the Metropolis-adjusted Langevin algorithm (MALA) when the mobility matrix is constant \cite{RoTw1996A, RoTw1996B, BoHa2013}.  MALA is a special case of the smart and hybrid Monte-Carlo algorithms, which are older and more general sampling methods \cite{RoDoFr1978, DuKePeRo1987}.  However, the Metropolized Ermak-McCammon scheme has two drawbacks: it involves the divergence of the mobility tensor, and worse, as illustrated in the next section, there are important situations where the acceptance probability in the scheme breaks down in the small noise limit.   The proposed integrator gets around these problems.

Next, we will use several examples to illustrate these properties of the Metropolis integrator then prove all the statements made in this section in \S\ref{sec:ergodicity},~\S\ref{sec:accuracy} and~\S\ref{sec:deterministicaccuracy}.

\section{Numerical Examples} \label{sec:numerics}

Unless  otherwise indicated, in this section we work with the Metropolis integrator, Algorithm \ref{MetropolisIntegrator}, operated with $G_h(x)$ and $B_h(x)$ given by  \eqref{rk2_Gh} and~\eqref{rk2_Bh} with parameter values \eqref{rk2_optimal_parameters}.  Test problems consist of the following self-adjoint diffusions: 

\smallskip

\begin{description}
\item[(E0)] Brownian particle with a heavy-tailed stationary density;
\item[(E1)] Brownian particle with a tilted square well potential energy;
\item[(E2)] 1D bead-spring chain in a confined solvent;
\item[(E3)] 3D bead-spring chain in an unbounded solvent; and,
\item[(E4)] Brownian particle with a two-dimensional double-well potential energy.
\end{description}

\smallskip

\noindent
Examples (E0) and (E1) involve a non-normalizable density. Example (E1) shows that standard explicit integrators may not detect properly features of a potential with jumps, which leads to potentially unnoticed but large errors in dynamic quantities such as mean first passage times.  As shown in (E2) and (E3), standard integrators can also fail when the potential contains a hard-core or Lennard-Jones-type component.  Example (E2) is not physically realistic, however, it is simple enough to be used for intensive numerical tests and yet complex enough to capture some of the essential features which make \eqref{sde} challenging to simulate including multiplicative noise, nontrivial $\divergence{M}(x)$, and steep potentials.   On the other hand, example (E3) is physically relevant to polymeric fluid simulations but, since the solvent is unbounded, $\divergence{M}(x) = 0$ in this case.   In sum, (E1) -- (E3) demonstrate two possible modes of failure of standard explicit integrators.  (E4) illustrates the properties of the proposed scheme in the small noise limit.   
In the numerical tests that follow, we assess accuracy and stability as a function of the time-step size parameter.  In practice, however, the computational cost/time of the algorithms should also be assessed, but since this cost is problem-specific, we do not carry out such assessments in this paper.

In the numerical tests that follow, we will compare the Metropolis integrator to the Fixman scheme.  There are variants of this predictor-corrector scheme,
and among these we implemented the following trapezoidal discretization: \begin{equation}
\begin{cases}
\tilde{X}_{1} = X_0 - h M(X_0)  D U(X_0)  + \sqrt{2 h} B(X_0) \xi \\
X_{1} = X_0 - \frac{h}{2} \left( M(X_0) D U(X_0)+  M(\tilde X_1) D U(\tilde X_1) \right) \\
\qquad \qquad +  \frac{\sqrt{2 h}}{2} \left( M(X_0) + M(\tilde X_1) \right) B(X_0)^{-T} \xi
\end{cases}
\tag{Fixman}
\end{equation}
where $\xi\in\mathbb{R}^n$ is a Gaussian random vector with 
mean zero and covariance $\E( \xi_i \xi_j ) = \beta^{-1} \delta_{ij}$.  This method is weakly first order accurate, second-order deterministically accurate, and second-order weakly accurate for the special case of additive noise. Notice that the scheme has the nice feature that it avoids computing the divergence of the mobility matrix.  The approximation to the noise in the Fixman integrator can be derived from a two-step discretization of the `kinetic stochastic integral', which is introduced in \cite{HuOt1998}.

\subsection{Example (E0): Brownian particle with a heavy-tailed stationary density} 
\begin{figure}[Ht!]
\centering
\includegraphics[width=1.0\textwidth]{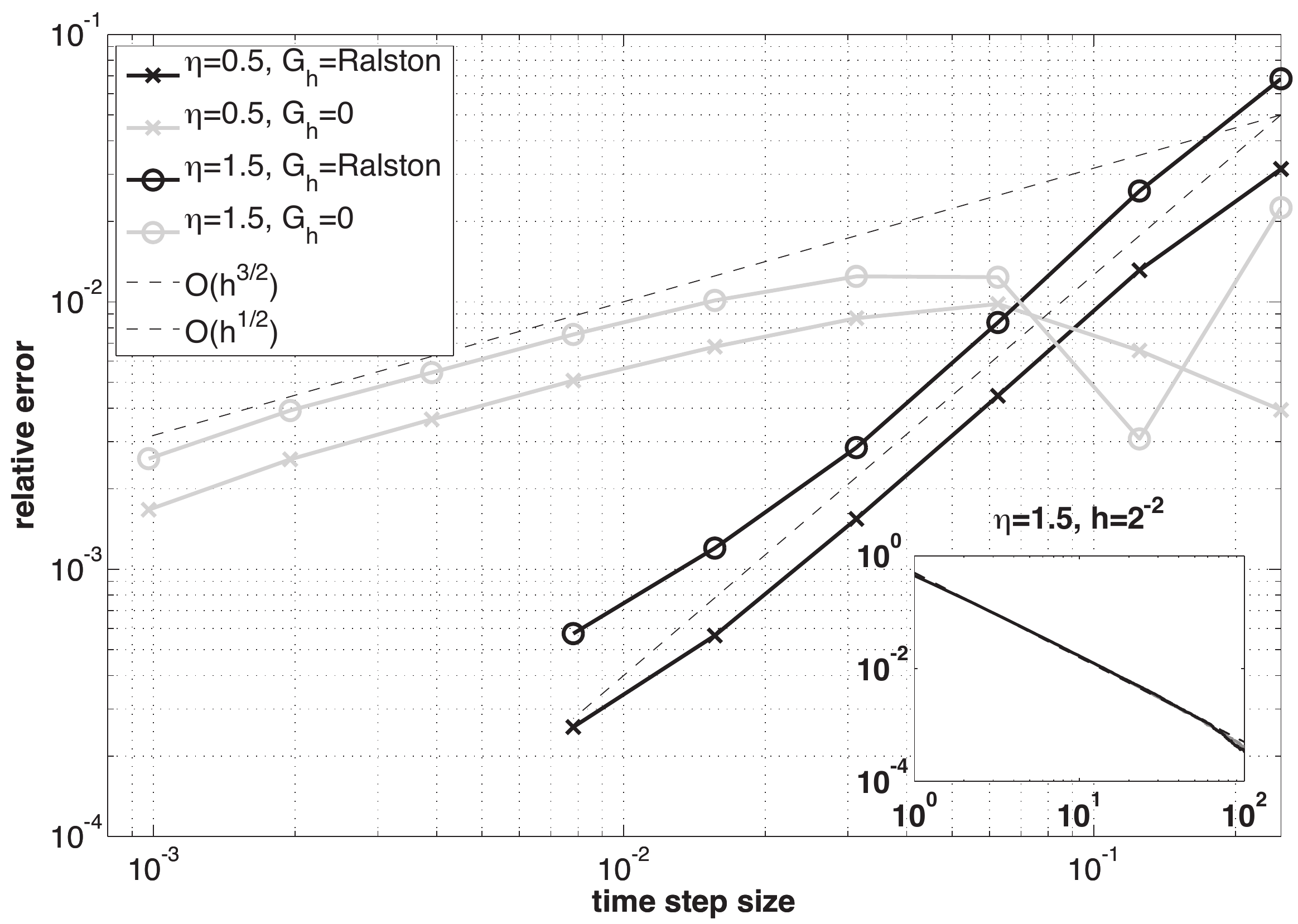}  
\caption{ \small {\bf Brownian Particle with a Heavy-tailed Stationary Density.} 
The main figure graphs the relative error in the Metropolis integrators' estimate of the mean-squared displacement of the
particle at time $T=1$ with initial condition $x=2$ as a function of the time step size $h$ for two different values of $\eta$ and using a
sample average with $10^9$ samples.  The black and gray colors indicate the Metropolis integrator with $G_h(x)$ chosen to be the Ralston Runge-Kutta combination \eqref{rk2_Gh}
and zero, respectively.  To derive a benchmark solution, the associated Fokker-Planck equation was numerically solved 
yielding the following approximations: $\E_x Y(1)^2 \approx 6.0487504$ ($\approx  4.7229797$) when $\eta=0.5$ ( $=1.5$), respectively.
The dashed lines in the main figure are $\mathcal{O}(h^{3/2})$ and  $\mathcal{O}(h^{1/2})$ reference slopes.  
When $\eta=1.5$ the stationary density of this process is normalizable, and the Metropolis integrators exactly reproduce this probability
density function, as confirmed in the inset.  Specifically, the  log-log plot  in the inset shows that the approximate densities of the 
Metropolis integrators (solid black and gray) overlap the true stationary probability density of the process (dashed).
 }
\label{fig:E0_heavy_tailed}
\end{figure}

To confirm that the proposed integrator applies to diffusions with non-normalizable densities $\nu(x)$, consider
a Brownian particle on the interval $\{ x \ge 1 \}$ governed by \eqref{sde} with $M=1$, $\beta=1$,  and the following  
potential energy: \begin{equation} \label{heavy_tailed}
U(x) = \eta \log(x) \;,  \quad (x \ge 1) \;, \quad (\eta > 0) 
\end{equation}
where $\eta$ is a positive parameter.  If $\eta>1$ the stationary density $\nu(x) = \exp(-U(x))$ is normalizable with $Z = \int_1^{\infty} \nu(x) dx = 1/(1-\eta)$ and the rate of convergence to equilibrium is sub-exponential.  In contrast, if $\eta \le 1$ the stationary density is no longer integrable; see Theorem 2.1 in \cite{hairer2009hot} for more details. 
Note that if a proposal move $X_1^{\star}$ in Algorithm~\ref{MetropolisIntegrator} is not in the interval $\{ x \ge 1 \}$, then this move is rejected by the algorithm.  
Figure~\ref{fig:E0_heavy_tailed} confirms that the Metropolis integrator is finite-time weakly accurate for a value of $\eta$ above and below one, and for two different choices of $G_h(x)$, namely the trivial choice $G_h(x)=0$ and the Ralston Runge-Kutta combination given in \eqref{rk2_Gh} with parameter values \eqref{rk2_optimal_parameters}. 
The $\mathcal{O}(h^{3/2})$ rate of convergence and improved accuracy of the latter, which uses a higher-order approximation to the drift, emphasizes that $G_h(x)=0$ is not a good choice.   In fact, Theorem~\ref{thm:finite_time_accuracy_additive_noise} proves that when the mobility matrix is constant, the Metropolis integrator with any second-order Runge-Kutta combination in \eqref{rk2_Gh} is weakly $3/2$-order accurate at constant temperature.  The inset in the figure confirms that the Metropolis integrators reproduce the stationary probability density of the diffusion when it exists, that is, when $\eta$ is above one.  Since all of the derivatives of the potential in \eqref{heavy_tailed} are continuous and bounded on $\{ x \ge 1 \}$, a standard integrator that can deal with moves that leave the interval $\{ x \ge 1 \}$ is sufficient in this example.  


\subsection{Example (E1): Brownian particle with a tilted square well potential energy}  The following example emphasizes that the Metropolis integrator, Algorithm~\ref{MetropolisIntegrator},  applies to self-adjoint diffusions even when the density $\nu(x)$ is not normalizable.  To introduce this example, it helps to consider simulating a Brownian particle moving in a regularized, periodic square well potential.    In order to adequately resolve the jump in the potential, a standard scheme requires a sufficiently small time step size.  Without this resolution the scheme's equilibrium distribution will be essentially uniform, and its estimate of dynamic quantities associated to crossing the square potential barrier will be inaccurate.    These predictions are confirmed in the following numerical experiment.    Inspired by \cite{ReVaLiHaRuPe2001}, we introduce a constant tilting force $F > 0$ so that the particle drifts to the right intermittently stopping at the jumps in the potential.  With this tilting force, Stratonovich was able to derive a formula for the mean first passage time, which we use below to benchmark and test the Fixman and Metropolis integrator \cite{St1958}.  At this point it is worth mentioning that when the mobility matrix in \eqref{sde} is constant, the Fixman scheme reduces to a second-order trapezoidal discretization of the drift and a first-order approximation to the noise.   

\begin{figure}[Ht!]
\centering
\includegraphics[width=1.0\textwidth]{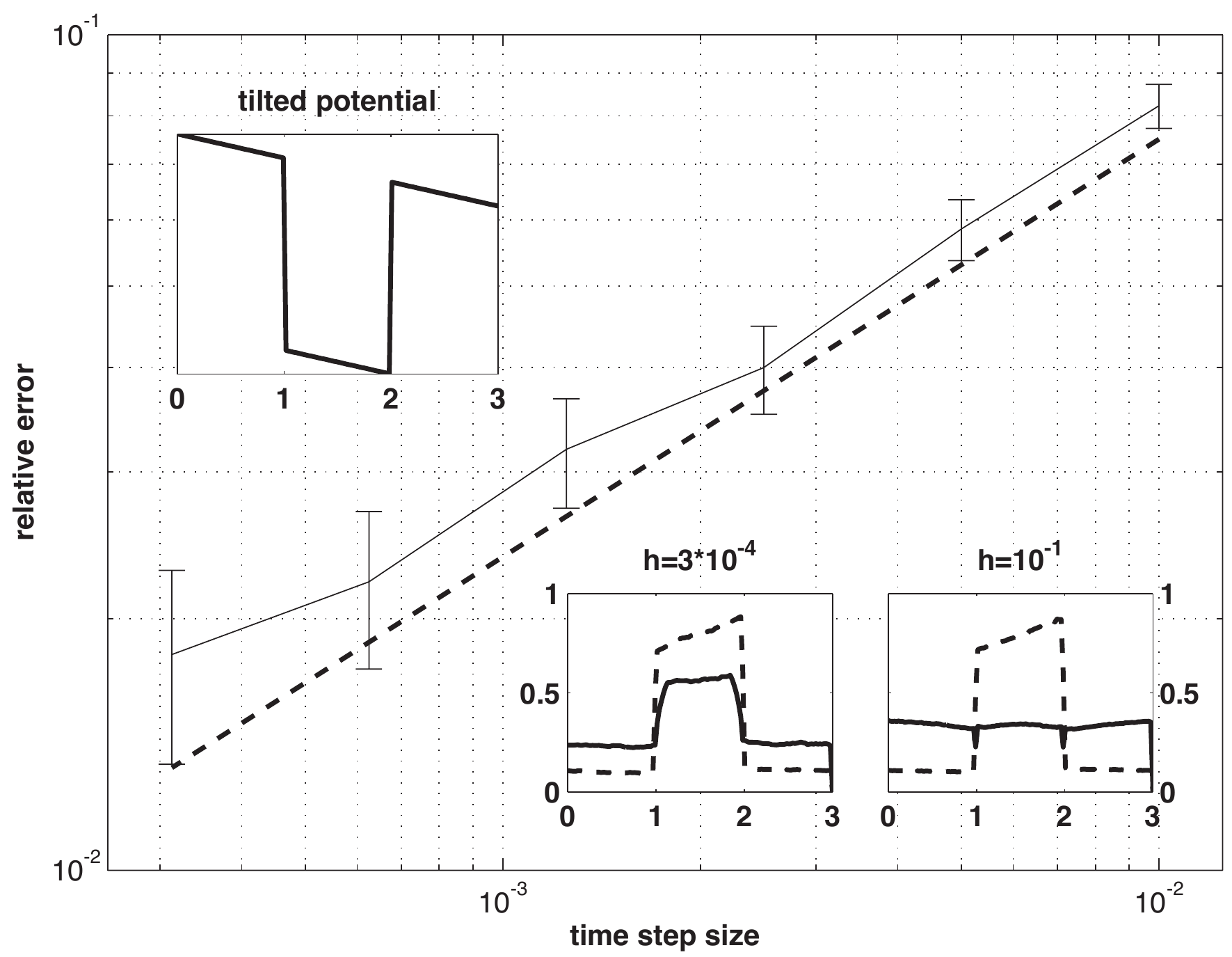}  
\caption{ \small {\bf Brownian Particle with a Tilted Square Well Potential.} 
The main figure graphs the relative error in the mean first passage time as a function of the time step size 
for the Metropolis integrator (solid).   The error bars represent a $99\%$ confidence interval.
The dashed line is an $\mathcal{O}(h^{1/2})$ reference slope.
The inset on the upper left shows one cycle of the tilted square well potential.  
The insets on the lower right plot the approximate density of the Fixman (solid) and  
Metropolis integrator (dashed) obtained by wrapping the points along a trajectory 
around the boundaries of a period of the square well potential.   By comparing the densities at the 
coarse and fine time step size, notice that the density of the Metropolis integrator is converged at the coarse
time step size, while the Fixman scheme at both the coarse and fine time step size is unable to detect properly 
the jump in the potential.  This difference is reflected in the main figure which shows that the Metropolis 
integrator is  about $2 \%$ accurate with a time step size $h = 10^{-3}$, while the Fixman scheme is  inaccurate ($>10 \%$ 
relative error) at the time step sizes tested.
 }
\label{fig:tiltedsquarewell}
\end{figure}

To be more precise, let $U(x)$ be a periodic, square well potential given by:
\begin{align*}
U(x) = \tanh\left( \frac{(x\bmod{3}) -2}{\epsilon} \right) - \tanh\left( \frac{(x\bmod{3}) -1}{\epsilon} \right)  
\end{align*}
where $\epsilon$ is a smoothness parameter, and $U(x) = U(x+3)$ for all $x \in \mathbb{R}$.  The period in this function is selected so that the jumps in $U(x)$ over one cycle $[0, 3]$  occur at $x=1$ and $x=2$.  A Brownian particle moving in a tilted square well potential satisfies an equation of the form \eqref{sde} with mobility equal to unity and, \[
dY = - \tilde U^{\prime}(Y) dt + \sqrt{2 \beta^{-1}} dW \;,  \quad Y(0) \in \mathbb{R} 
\]
where we have introduced the following tilted potential energy function: \[
\tilde U(x) = U(x) - F x \;.
\]  Observe that when $F=0$ the potential $\tilde U(x)$ reduces to $U(x)$.   For every $F \in \mathbb{R}$, it is straightforward to use \eqref{eq:generator} to show that the generator of $Y(t)$ is self-adjoint with respect to the density $\nu(x) = \exp(- \beta \tilde U(x) )$.    When $F=0$ the stationary density is normalizable over one period of $U(x)$.  When $F>0$ the constant tilting force gives rise to a downward tilt in the potential $\tilde U(x)$ as shown in the northwest inset in  Figure~\ref{fig:tiltedsquarewell}, and so, $\nu(x)$ is not normalizable since $\int_{\mathbb{R}} \nu(x) dx = + \infty$.     The first passage time of $Y(t)$ from any $x_0 \in \mathbb{R}$ to $x_0 + 3$ is defined as: \[
\tau = \inf \{ \; t>0 ~:~ Y(0) = x_0\;, ~ Y(t) \ge x_0 + 3 \} \;.
\] The numerical experiments that follow consist of launching an integrator with initial condition $X_0=x_0$ and time step size $h$, and terminating the simulation at the first time step $n$ where $X_n \ge x_0 + 3$.  To avoid systematically overestimating the first passage time, we use the approximation, $\tau \approx n h - h/2$.  As a side note, we mention that the accuracy of this approximation to the mean first passage time (which goes like $\mathcal{O}(\sqrt{h})$ \cite{gobet2004exact,gobet2007discrete,gobet2010stopped}) can be improved upon by accounting for the probability that the particle reaches $x_0 + 3$ in between each discrete step~\cite{Ma1999}.

The numerical and physical parameters used in the numerical experiments are given in Table~\ref{tab:tiltedsquarewell}. To visualize the long-term behavior of the schemes, it helps to plot the probability density of points produced by each scheme modulo a period of the square well potential.  An approximation to this density is shown in the southeast insets in Figure~\ref{fig:tiltedsquarewell}.    At both the coarse and fine time step size tested, the Fixman scheme underestimates the barrier height, and consequently, numerical tests show that it grossly underestimates, e.g., the mean first passage time between wells.   This underestimate persists unless its time step size is small enough to resolve the barrier ($h < 10^{-4}$).  The Metropolis integrator is able to capture the features of the potential even at the coarse time step size $h=10^{-1}$, and its approximation to the mean first passage time is about $2 \%$ accurate with a time step size $h = 10^{-3}$.   


 \setlength\tabcolsep{1pt}

 \begin{table}[Ht!]
\centering
\begin{tabular}{|c|c|c|}
\hline
\multicolumn{1}{|c|}{\bf Parameter} & \multicolumn{1}{|c|}{\bf Description} &  \multicolumn{1}{|c|}{\bf Value}   \\
\hline
\hline
\multicolumn{3}{|c|}{ {\em Physical Parameters} } \\
\hline
\hline
$F$ & constant tilting force & $0.25$ \\
  \hline
$\beta$ & inverse temperature factor & $1$  \\
  \hline
  $\epsilon$ & \begin{tabular}{c} square well \\ smoothness parameter \end{tabular} & $0.001$  \\
\hline
\hline
\multicolumn{3}{|c|}{{\em  Numerical Parameters }} \\
\hline
\hline
$h$ & time step size & \begin{tabular}{c} $\{ 0.01, 0.05, 0.025,$ \\ 
        $ 0.0125, 0.00625,   0.003125 \}$ \end{tabular} \\
\hline
$N_s$ & \begin{tabular}{c} \# of first passage \\ time samples \end{tabular} & $10^5$ \\
\hline
\end{tabular} 
\caption{  {\bf Simulation Parameters for a Brownian Particle with a Tilted Square Well Potential.}   }
\label{tab:tiltedsquarewell}
\end{table}


\subsection{Example (E2): 1D bead-spring chain with hydrodynamic interactions}  

\begin{figure}[Ht!]
\centering
\includegraphics[width=1.0\textwidth]{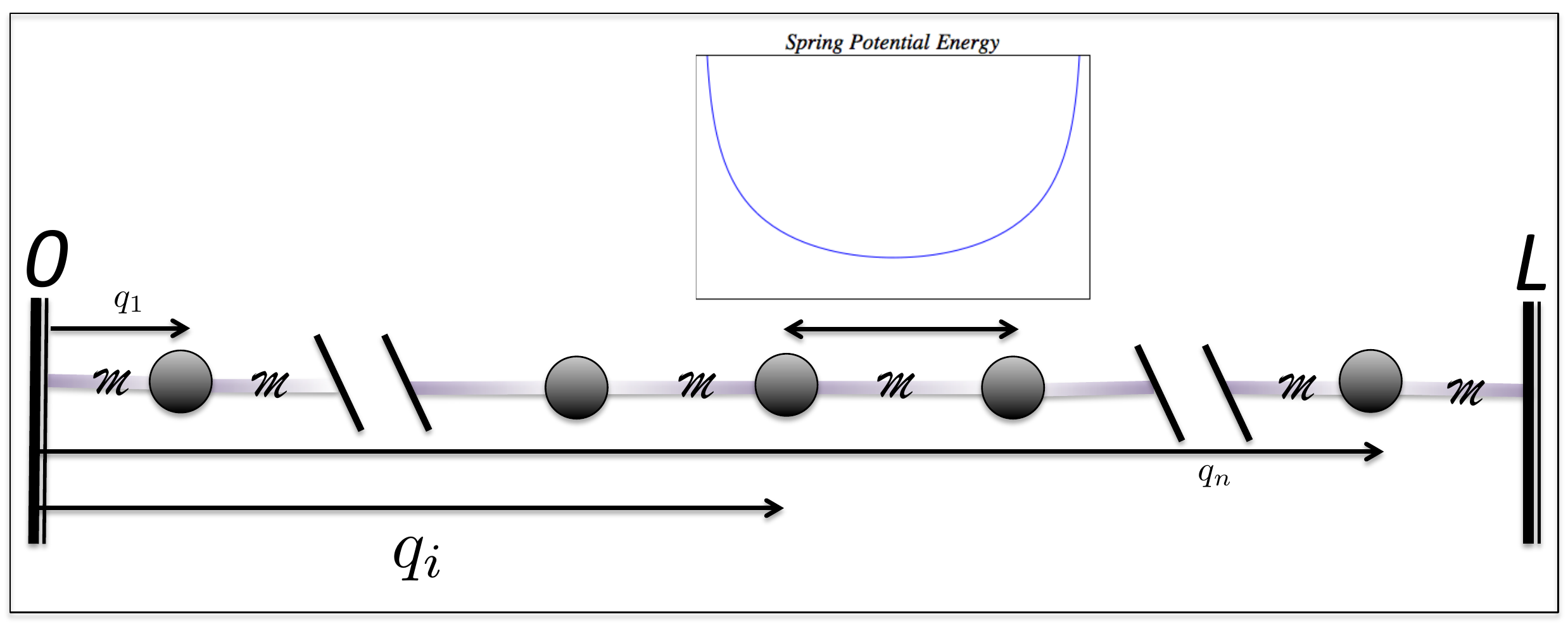}  
\caption{ \small  This figure shows a one-dimensional bead spring chain in a solvent.  The chain consists of $n$ beads connected by $n+1$ springs confined to an interval $[0, L]$.  The spring potential energies are modeled by a two-sided FENE potential \eqref{fene} that is plotted in the inset.  The singularities in this potential enforce excluded volume between beads and finite  extension of the spring length.   The solvent velocity satisfies a steady-state Stokes equation with point sources of drag located at the bead positions.  }
\label{fig:1Dbeadspringchain}
\end{figure}

This example confirms Properties (P2) and (P3) of the Metropolis integrator.
Consider a 1D bead-spring chain consisting of $n$ beads and $n+1$ springs confined to an interval $[0, L]$ and immersed in a `fictitious solvent' with viscosity $\mu$ as illustrated in Figure~\ref{fig:1Dbeadspringchain}.  Incompressibility implies that the velocity of a 1D solvent is constant, and so to obtain nontrivial hydrodynamic interactions, we remove this physical constraint.  Even though the resulting solvent does not satisfy incompressibility, this example captures some of the key elements of hydrodynamic interactions -- including a non-trivial $\divergence M (x)$ -- while being simple enough to permit detailed numerical studies.

We begin by writing this system as a self-adjoint diffusion of the form \eqref{sde} with a normalizable density $\nu(x)$.  Assume that the bead and solvent inertia are negligible.   Let $q_i$ and $F_i$  denote the $ith$ bead position and force, respectively; 
let $u(q) \in \mathbb{R}$ denote the solvent velocity for $q \in [0, L]$; and, let  $u_i = u(q_i)$ for $i=1, \cdots, n$. Order the particle positions so that: \[
 0 < q_1 < \cdots < q_n < L  \;.
\]
This ordering is maintained by the following soft-core spring potential:
\begin{equation}
U_{\text{FENE}}(x) = -    \frac{\epsilon}{2}  \log\left( 1 - \left( \frac{x}{2 \ell} \right)^2 - \left(1 - \frac{x}{2 \ell} \right)^2 \right)  
\label{fene}
\end{equation}   
where the acronym stands for finitely extensible nonlinear elastic (FENE) spring force law.
There are two parameters in this potential: an energy constant $\epsilon$ and bond length at rest $\ell$.  The linear behavior of $U_{\text{FENE}}(x)$ about its resting position $x=\ell$ is given by a Hookean spring with stiffness $\epsilon/\ell^2$.  Moreover, this potential possesses a singularity at $x=0$ that prevents interbead collisions, and a second singularity at $x=2 \ell$ to enforce finite extension of the spring length.  Potentials of this type play an important role in capturing the right non-Newtonian behavior of a dilute polymer solution using bead-spring models \cite{BiArHaCu1987, Ot1996}.  

In addition to neighboring spring interactions, the particles are coupled by non-bonded interactions mediated by a fictitious solvent.   Since the bead and solvent inertia are negligible, a balance between the solvent viscous force and the bead spring forces leads to the following equations for the solvent velocity: \begin{equation} \label{1Dstokes}
\mu \frac{\partial^2 u}{\partial x^2} = - \sum_{i=1}^n F_i \; \delta (q - q_i)   \;,~~u(0)=u(L)=0 
\end{equation}
where $\delta( q- q_i )$ is a point mass located at the position of the $ith$ bead $q_i$, and $F_i$ is the force on the $ith$ bead.   In this equation the location of the beads and the forces on the beads are known, and the solvent velocity is unknown.  Let $x = (q_1, \cdots, q_n)^T$, and in terms of which, the mobility matrix $M(x)$ in \eqref{sde} comes from solving \eqref{1Dstokes} and it describes the linear relationship between the solvent velocities at the bead positions and the bead forces:   \begin{equation} \label{1Dmobility}
\sum_{j=1}^n M_{ij} (x) F_j = u_i \;.
\end{equation}
The deterministic drift in \eqref{sde} arises from assuming the solvent velocity matches the bead velocities at the bead positions.  A heat bath is then added to this drift to ensure that the density of the stationary distribution of the resulting equations is proportional to $\nu(x)$ in \eqref{eq:nu}, where the total potential energy is just the sum of the $n+1$ spring potential energies: \[
U(x) = \sum_{i=1}^{n+1} U_{\text{FENE}}( q_i - q_{i-1} ) \;, ~~  \text{where $q_0 = 0$ and $q_{n+1} = L$} \;.
\]
Now we show how to solve \eqref{1Dstokes} for the solvent velocity, and in the process, derive a procedure for exactly evaluating the mobility matrix in this specific case.

In terms of the friction matrix $\Gamma(x) = M^{-1}(x)$, the linear relationship~\eqref{1Dmobility} can be written as:
 \begin{equation} \label{1Dfriction}
F_i = \sum_{j=1}^n \Gamma_{ij} u_j \;.
\end{equation}
Given the instantaneous forces and positions of the particles, the solution to \eqref{1Dstokes} can be derived as follows.  In between the particles, the fluid velocity is linear since the solvent experiences no force there.  At the $ith$ particle position, there is a discontinuity in the derivative of the fluid velocity: \begin{equation} \label{1Dkink}
\frac{\partial u}{\partial q}(q = q_i^+)  - \frac{\partial u}{\partial q} (q= q_i^-) =  - \frac{F_i}{\mu} 
\end{equation}
for $i=1, \cdots, n$.  With $u_0=u(q_0)=0$ and $u_{n+1} = u(q_{n+1})=0$, notice that \eqref{1Dkink} implies: \begin{equation}
\frac{u_{i+1} - u_i}{ q_{i+1} -q_i} - \frac{u_i - u_{i-1}}{q_i - q_{i-1}} =  - \frac{F_i}{\mu} \;.
\end{equation}  
Comparing this equation with \eqref{1Dfriction}, it is evident that $\Gamma$ is a tridiagonal matrix with entries given by: \begin{equation}
\Gamma_{ij} = \begin{cases}
 \frac{\mu}{q_{i+1} - q_i} + \frac{\mu}{q_i - q_{i-1}} \;,  & \text{if $i=j$} \;, \\
- \frac{\mu}{q_{i+1} - q_i} \;, & \text{if $i-j=1$} \;, \\
- \frac{\mu}{q_i - q_{i-1}} \;,  & \text{if $i-j=-1$} \;.
\end{cases} 
\end{equation}
These expressions are explicit and straightforward to implement in the Fixman scheme or the Metropolis integrator, Algorithm~\ref{MetropolisIntegrator}.  Note that the friction and mobility matrices are symmetric positive-definite, and the friction matrix is tridiagonal while the mobility matrix is not tridiagonal, in general.   

\begin{figure}[Ht!]
\centering
\includegraphics[width=1.0\textwidth]{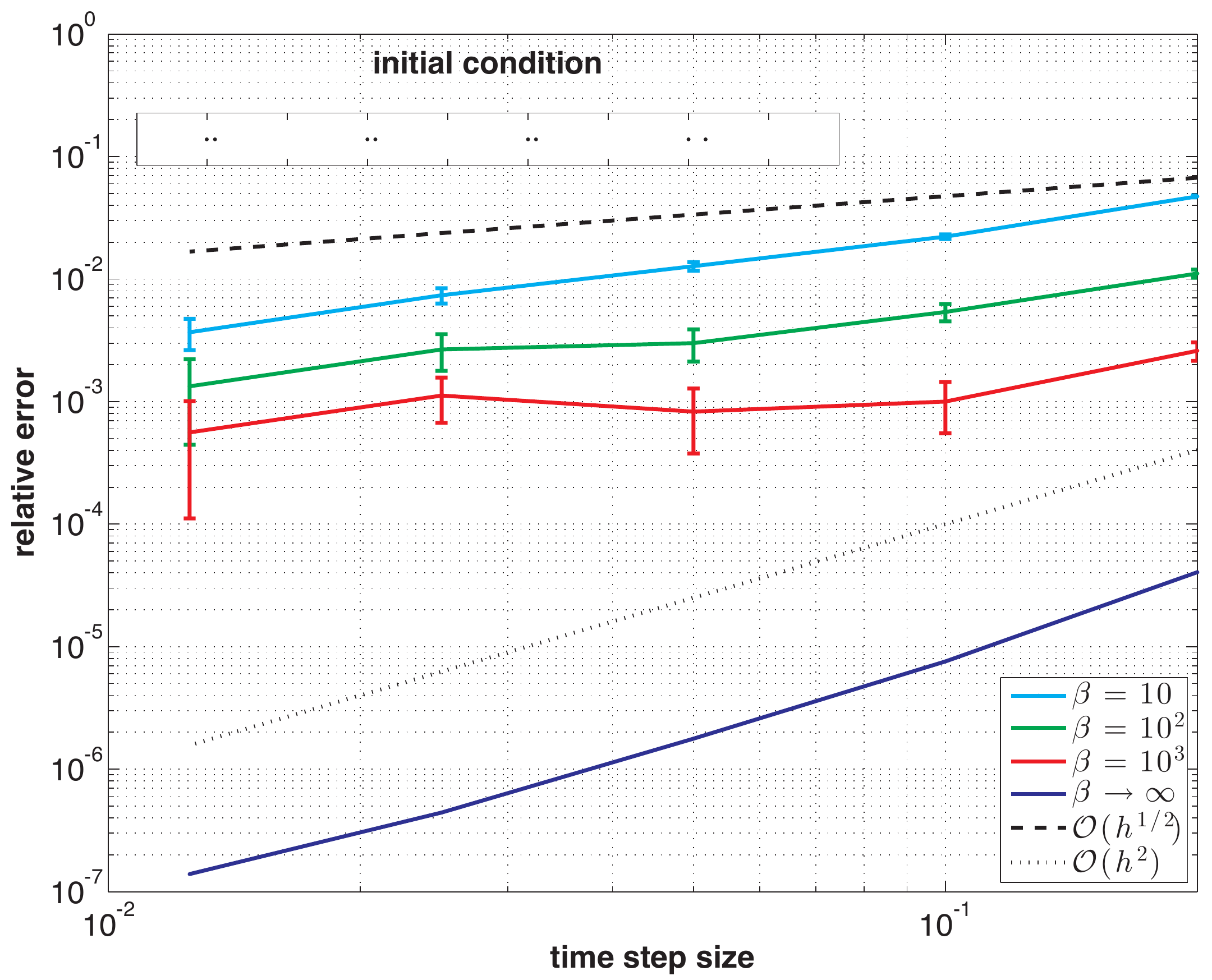}  
\caption{ \small {\bf One-dimensional bead-spring chain in a solvent.}
The main figure shows  the relative error in the mean position of the particles \eqref{observable1} computed using the Metropolis integrator at the $\beta$ values shown in the legend. 
The error bars represent a $95\%$ confidence interval.
The northwest inset shows the initial condition (dots) and equilibrium positions (tick marks) of the beads.  Observe that the Metropolis integrator remains within $1\%$ accurate at the large time step size $h=0.1$ for a wide range of temperatures.  This figure confirms that the Metropolis integrator is weakly accurate at constant temperature Property (P2), and second-order accurate in the small noise limit Property (P3).  The curve $\beta \to \infty$ corresponds to a $\beta$ that is large enough ($\beta = 10^{12}$) so that the $\mathcal{O}(h^2)$ rate of convergence becomes visible.
}
\label{fig:brownianlattice}
\end{figure}

For the numerical experiment, consider an eight bead chain that is initially compressed  as shown in the northwest inset of Figure~\ref{fig:brownianlattice}.    
Setting  $\mu$, $L$, and $\epsilon$ equal to unity is equivalent to rescaling the system by a characteristic length, time and energy scale.  This leaves the inverse temperature factor $\beta$ as a free physical parameter that we vary in the numerical experiments.   These simulation parameter values are provided in Table~\ref{tab:brownianlattice}.  
The equilibrium positions of the beads are indicated by the tick marks in this inset.  We use the Fixman and Metropolis integrator, Algorithm~\ref{MetropolisIntegrator}, to simulate the relaxation of the system from this state.  Specifically, the schemes compute the expected value of the average position of the particles as the lattice relaxes towards equilibrium:  \begin{equation} \label{observable1}
\E_x \left( \frac{1}{n} \sum_{i=1}^n q_i(t) \right) 
 \end{equation} for $t \in [0, T]$, where $T=50$ time units was found to be sufficient to capture the relaxation dynamics for all parameter values tested.  

The relative error of the Metropolis integrator in computing \eqref{observable1} is shown in Figure~\ref{fig:brownianlattice} for a range of temperatures.  Near the deterministic limit, we used a forward Euler scheme operated at a very small time step size to obtain a benchmark solution ($h<10^{-5}$), and observe from Figure~\ref{fig:brownianlattice},  that the Metropolis integrator is accurate for time step sizes that exceed $h =0.05$.    As the temperature is increased, and in the absence of an analytical solution, we used Richardson extrapolation to obtain a benchmark solution.  Using this benchmark we estimated the error of the Metropolis integrator at various temperatures.   These graphs show that although the relative error increases with increasing temperature, the Metropolis integrator remains within $1\%$ error for a time step size $h=0.1$ from the deterministic limit, $\beta=10^{12}$, to  the moderate temperature, $\beta=10$. 

To deal with stochastic instabilities in the Fixman scheme,  we employed `the method of rejecting exploding trajectories' \cite{MiTr2005}.   In this stabilization technique, if a non-physical move is generated by the Fixman scheme the entire sample path is rejected.  At the very low temperature of $\beta=1000$, and if $h\le 0.0125$, the Fixman scheme rejects a negligible amount of trajectories ($<1\%$) and is as accurate as the Metropolis integrator operated at $h=0.2$.  However, the Fixman scheme rejects almost all trajectories generated at the low temperature $\beta=100$ even if the time step size is reduced to $h=0.0001$, and its performance is much worse at the moderate temperature $\beta=10$.


 \begin{table}[Ht!]
\centering
\begin{tabular}{|c|c|c|}
\hline
\multicolumn{1}{|c|}{\bf Parameter} & \multicolumn{1}{|c|}{\bf Description} &  \multicolumn{1}{|c|}{\bf Value}   \\
\hline
\hline
\multicolumn{3}{|c|}{ {\em Physical Parameters} } \\
\hline
\hline
$\mu$ & solvent viscosity & $1$ \\
\hline
$L$ & chain length & $1$ \\
  \hline
  $\epsilon$ & energy constant & $1$  \\
  \hline
$\beta$ & inverse temperature factor & $\{ 10^{12}, 10^3, 10^2, 10 \}$  \\
  \hline
$n$ & \# of beads & $8$  \\
\hline
$T$ & \begin{tabular}{c} time-span \\ for simulation \end{tabular} & $50$ \\
\hline
\hline
\multicolumn{3}{|c|}{{\em  Numerical Parameters }} \\
\hline
\hline
$h$ & time step size & $\{ 0.2, 0.1, 0.05, 0.025, 0.0125 \}$ \\
\hline
$N_s$ & \begin{tabular}{c} \# of sample \\ trajectories \end{tabular} & $10^5$ \\
\hline
$(q_1(0), \cdots, q_n(0))$ & initial positions & \begin{tabular}{c} $( 0.1,0.11,0.33,0.34,$ \\
$0.56,0.57, 0.79, 0.81)$ \end{tabular} \\
\hline
$(\tilde q_1, \cdots, \tilde q_n)$ & \begin{tabular}{c} equilibrium  \\ positions of springs \end{tabular} & \begin{tabular}{c} $(0.1,0.21,0.33, 0.44, $ \\
$0.55,  0.67,    0.79,   0.9)$ \end{tabular} \\
\hline
\end{tabular}
\caption{  {\bf Simulation parameters for a one-dimensional bead-spring chain in a solvent.}   
}
\label{tab:brownianlattice}
\end{table}


\clearpage

\subsection{Example (E3): DNA dynamics in a solvent}  The following example applies the Metropolis integrator to a BD simulation of a DNA molecule with hydrodynamic interactions. Consider a bead-spring model of a bacteriophage DNA molecule with $N_b = 11$ spherical beads, where each spring approximates the effect of 4850 base pairs, so that ten springs contain approximately the number of base pairs in a DNA molecule with a contour length of $21 \mu m$ \cite{ZhDoWeAlGrde2009}.  In this case the energy function $U(x)$ is simply the sum of the spring potential energies, whose formula is written below in \eqref{eq:wlc}.  Assume the beads are spherical with radius $R_b$ and move in a Stokesian solvent with viscosity $\eta_s$.   To be specific, the solvent velocity $u(q) \in \mathbb{R}^3$ and pressure $p(q) \in \mathbb{R}$ satisfy: \begin{equation} \label{stokes}
 \eta_s (\nabla^2 u)(q) - (\nabla p)(q) + f(q) = 0   \;, ~~( \nabla \cdot u)(q) = 0 \;,~~\text{for all $q \in \Omega \subset \mathbb{R}^3$}
\end{equation}
where $\Omega$ is the domain of the solvent and $f(q) \in\mathbb{R}^3 $ is the applied force density due to the beads.  We augment these equations with the following boundary conditions: the fluid is at rest at infinity and satisfies no-slip conditions on the surfaces of each bead.  Let $q_i$, $v_i$, and $F_i$ denote the position, translational velocity, and force of the $ith$ bead where $i \in \{ 1, \cdots, N_b \}$.  Introduce the $3 N_b$ dimensional vectors of bead positions $x = (q_1, \cdots, q_{N_b})^T$, bead translational velocities $V = (v_1, \cdots, v_{N_b})^T$ and bead forces $F = (F_1, \cdots, F_{N_b})^T$.  An immediate consequence of the linearity of the Stokes equation~\eqref{stokes} is that:  \[
V = M(x) F 
\]  
where $M(x)$ is the $ 3 N_b \times 3 N_b$ mobility matrix.    This linear relationship always holds for bodies moving in a Stokes fluid.  Moreover, the matrix $M(x)$ is always symmetric and positive definite for every $x \in \mathbb{R}^{3 N_b}$.  

Determining the entries of the mobility matrix requires solving the Stokes equation \eqref{stokes} for the solvent velocity field which is a very complicated boundary value problem.  
This difficulty motivates using the Rotne-Pragner-Yamakawa (RPY) approximation of the solvent velocity field, which leads to the following approximate mobility matrix: \begin{equation} \label{RPYmobility}
M(x) = \begin{bmatrix} \Omega_{1,1} & \cdots & \Omega_{1,N_b} \\
\vdots & \ddots & \vdots \\
\Omega_{N_b,1} & \cdots & \Omega_{N_b, N_b} 
\end{bmatrix} \;,~~
\Omega_{i,j} = \begin{cases} \frac{1}{\zeta} I_{3 \times 3} \;,  & \text{if} ~~ i=j  \\
\Omega_{RPY}(q_i - q_j) \;, & \text{otherwise} 
\end{cases}
\end{equation}
for all $x \in \mathbb{R}^{3 N_b}$.  Here, we have introduced the $3 \times 3$ matrix $\Omega_{RPY}(q)$ defined as:
\begin{equation} \label{rpy}
\Omega_{RPY}(q) =  \frac{1}{\zeta} \left(  C_1(q) I_{3 \times 3} + C_2(q) \frac{q}{|q|} \otimes \frac{q}{|q|}  \right) 
\end{equation} where $C_1(q)$ and $C_2(q)$ are the following scalar-valued functions: \[
C_1(q) = \begin{cases} \frac{3}{4} \left( \frac{R_b}{|q|} \right) + \frac{1}{2} \left( \frac{R_b}{|q|} \right)^3 \;, & \text{if} ~~|q|>2 R_b   \\
 1 - \frac{9}{32} \left( \frac{|q|}{R_b} \right) \;, & \text{otherwise}   \end{cases} 
 \]  and, \[
 C_2(q) = \begin{cases}  \frac{3}{4} \left( \frac{R_b}{|q|} \right) - \frac{3}{2} \left( \frac{R_b}{|q|} \right)^3 \;, & \text{if} ~~|q|>2 R_b   \\ ~
  \frac{3}{32} \left( \frac{|q|}{R_b} \right) \;,  & \text{otherwise}  \end{cases} 
\]
The quantity $1/\zeta$ is the mobility constant produced by a single bead translating in an unbounded solvent at a constant velocity: $\zeta = 6 \pi \eta_s R_b$.   
The approximation \eqref{RPYmobility} preserves the physical property that the mobility matrix is positive definite, satisfies $(\divergence M)(x) =0$ for all $x \in \mathbb{R}^{3 N_b}$,  and is exact up to $\mathcal{O}((R_b/r_{ij})^4)$ where $R_b$ is the bead radius and $r_{ij}$ is the distance between distinct beads $i$ and $j$ \cite{RoPr1969}. 

A well-defined characteristic length of DNA is its Kuhn length $b_k$ which represents the distance along the polymer chain over which orientational correlations decay.  (The bending rigidity of a polymer decreases with increasing $b_k$.)  For DNA, the Kuhn length is approximately a tenth of a micrometer.  In terms of which, consider a worm like chain (WLC) model for the spring potential energy given by: \begin{equation} \label{eq:wlc}
U_{\text{WLC}}(r) = \frac{\beta^{-1}}{2 b_k} \left( \frac{\ell^2}{\ell - r} - r + \frac{2 r^2}{\ell} \right) 
\end{equation}
where $\ell$ is the maximum length of each spring.   This empirical potential energy captures the entropic elasticity which causes the DNA molecule to be in a tightly coiled state.  The linear behavior of this potential is given by a Hookean spring with spring constant: $H_s = 3 \beta^{-1}/(b_k \ell)$.  The strength of the hydrodynamic interactions can be quantified by using this spring constant to compute the dimensionless bead radius: $a^* = R_b / \sqrt{\beta^{-1} / H_s}$.  Using the DNA parameter values provided in Table~\ref{tab:DNAparameters}, we find that $a^* = 0.291$ which for a dimensionless Rouse model (same bead-spring chain, but with Hookean springs) signifies a moderate strength of hydrodynamic interactions \cite{Ot1996}.  Since the linear Rouse model has the same features as the nonlinear DNA model minus the steep potential, a reasonable time step size for the DNA simulation can be determined by simulating a dimensionless Rouse model at this value of $a^*$ and the same non-random initial condition.  In particular, numerical experiments reveal that a time step size of $h \approx 10^{-4}$ leads to an average acceptance probability of about  $98\% - 99\%$ as the Rouse chain transitions from a stretched to a coiled state.    

\begin{figure}[Hb!]
\begin{center} 
\includegraphics[width=1.0\textwidth]{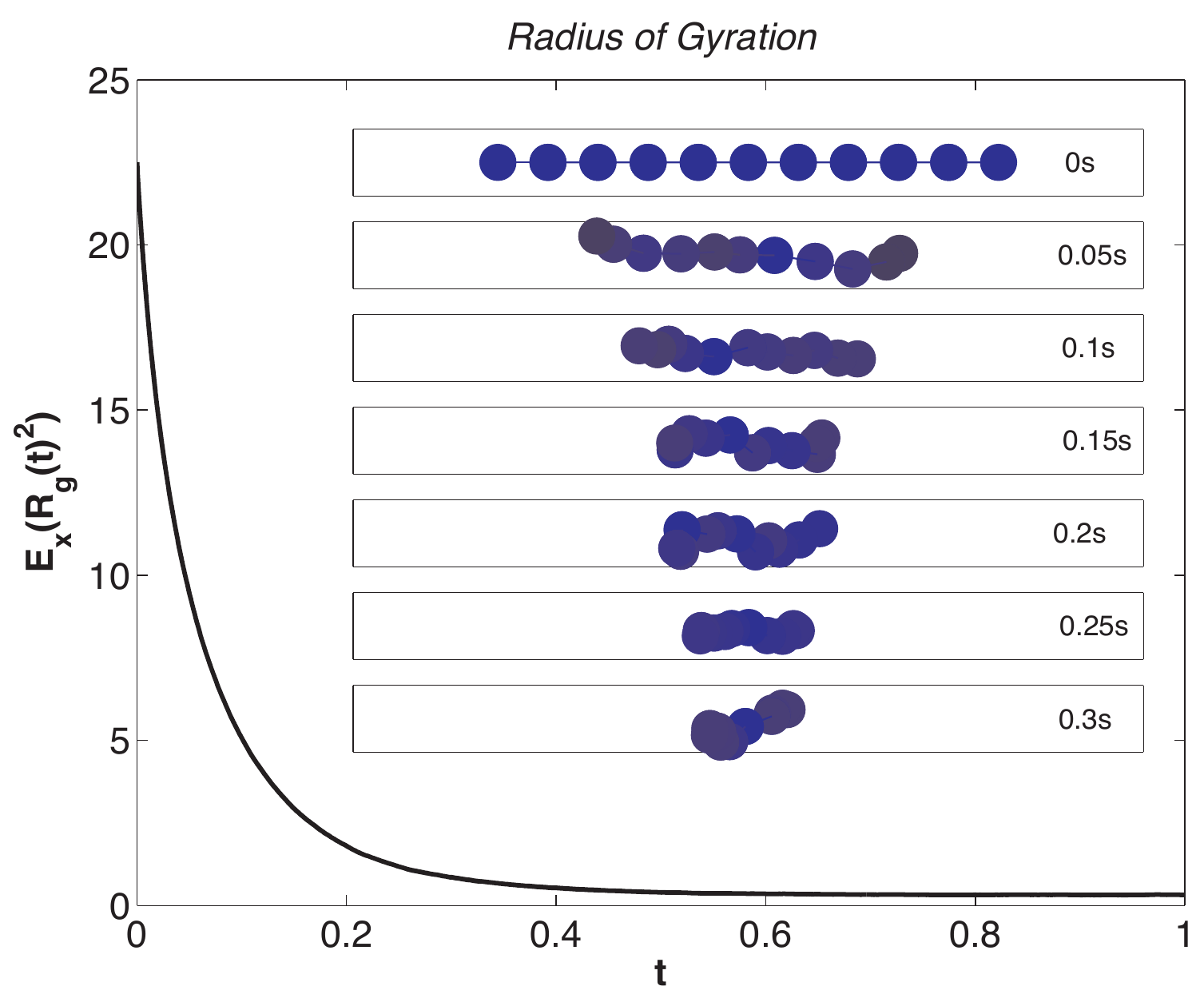} 
\caption{ \small {\bf DNA Simulation.}    The main figure plots the exponential decay in the 
mean-squared radius of gyration and the inset plots snapshots of the relaxation of a DNA molecule from an initial conformation shown 
in the top-most panel of the inset.  As shown in that panel,  the beads are initially evenly spaced on the x-axis with a spacing 
of $1.5 \mu m$.  In this example, the Metropolis integrator is able to take about one to two orders of 
magnitude larger time step size than conventional explicit schemes. 
  } \label{fig:DNAsimulation} 
\end{center}
\end{figure}

With approximately this time step size $h=( 10^{-4} )/ 2$, we use the Metropolis integrator to generate one thousand trajectories of the DNA chain from the initial conformation shown in the inset of Figure~\ref{fig:DNAsimulation} using the values provided in Table~\ref{tab:DNAparameters} over the interval $[0, 1]$.    This time-span is sufficiently long to capture the relaxation dynamics of the chain.   The average acceptance probability at this time step size is about $98 \%$, which implies that approximately the same time step size works for both the Rouse and DNA model, and emphasizes that the Metropolis integrator is to some extent insensitive to the stiffness in the WLC potential.   From the simulation data, the radius of gyration was estimated and this estimate is plotted in Figure~\ref{fig:DNAsimulation}.    Repeating this experiment at higher time resolution led to no noticeable change in this graph.   The Metropolis integrator seems able to compute dynamics for this system at a time step size that is more than $10 \times $ larger than what is possible using the Fixman scheme combined with `the method of rejecting exploding trajectories' described in the 1D bead-spring example.  


 \begin{table}[Ht!]
\centering
\begin{tabular}{|c|c|c|c|}
\hline
\multicolumn{1}{|c|}{\bf Parameter} & \multicolumn{1}{|c|}{\bf Description} &  \multicolumn{1}{|c|}{\bf Value} & \multicolumn{1}{|c|}{\bf Units}   \\
\hline
\hline
\multicolumn{4}{|c|}{{\em  Physical Parameters }} \\
\hline
\hline
$N_b$ &  \# of beads & $11$ & \\
\hline
$b_k$ & Kuhn length & $1 \times 10^{-1}$ & $\mu m$ \\
\hline
$R_b$ &  bead radius & $7.7 \times 10^{-2}$ & $\mu m$ \\
\hline
$\ell$  &  maximum spring length & $2.1$ & $\mu m$ \\
\hline
$\eta_s$  & solvent viscosity & $1 \times 10^{-9} $ &  $\frac{kg}{\mu m \cdot s}$ \\
\hline
$\beta^{-1}$ &  thermal energy & $4.11 \times 10^{-9}$  & $ \frac{kg (\mu m)^2}{s^2}$ \\
\hline
$T$ & time-span & $1$ & $s$ \\
\hline
\hline
\multicolumn{4}{|c|}{{\em  Numerical Parameters }} \\
\hline
\hline
$h$ & time step size & $0.5 \times 10^{-4} $ & $s$ \\
\hline
$N_s$ & \begin{tabular}{c} \# of sample \\ trajectories \end{tabular} & $10^3$  & \\
\hline
\end{tabular}
\caption{ \small {\bf DNA Simulation Parameters.}   }
\label{tab:DNAparameters}
\end{table}


\clearpage

\subsection{Example (E4): Brownian particle with a double well potential energy}  \label{ex:2Ddoublewell}

In the small noise limit, the proposal move \eqref{proposal} converges to a deterministic update.  If this update is a higher order accurate discretization of the zero noise limit of \eqref{sde}, then the proposal move will be deterministically accurate to that order.    However, it does not follow from this alone that the actual update in \eqref{actualupdate} is deterministically accurate too.  Indeed, Property (P3) also requires that all moves are accepted in the small noise limit.  This requirement may appear to contradict Property (P2): that the algorithm is weakly accurate so long as the noise is handled correctly regardless of the proposal move used, but it does not.  To be perfectly clear, the properties of the scheme in the small noise limit are asymptotic statements in $\beta$ while keeping the time step size $h$ fixed (and sufficiently small for certain estimates to be valid), whereas Property (P2) assumes that $\beta$ is fixed and $h$ is sufficiently small.

In \S\ref{sec:deterministicaccuracy} we show that for general $G_h(x)$ and $B_h(x)$, \begin{equation} \label{asymptotic_alpha}
\alpha_h(X_0, \xi) \sim  1 \wedge \exp\left(-\beta \mathcal{E}(X_0, h) \right) ~~ \text{as} ~~ \beta \to \infty 
\end{equation}
where we have introduced  \begin{equation} \label{fbeta}
\begin{cases}
\mathcal{E}(x, h) = U(x + h G_h(x)) - U(x) + h G_h(x)^T M_h(x + h G_h(x))^{-1} G_h(x) \;, \\
M_h(x) = B_h(x) B_h(x)^T \;.
\end{cases}
\end{equation}  
From \eqref{asymptotic_alpha} it follows that if $\mathcal{E}(x, h) < 0$ for all $x \in \mathbb{R}^n$, then all moves are accepted in the small noise limit, and the Metropolis integrator acquires the deterministic accuracy of its proposal move.   On the other hand, proposal moves emanating from points where $\mathcal{E}(x, h) > 0$ will most likely be rejected if the noise is small enough, and therefore, in the small noise limit the Metropolis integrator is not deterministically accurate at these points.  In \S\ref{sec:deterministicaccuracy} it is proved that when $G_h(x)$ and $B_h(x)$ are given respectively by \eqref{rk2_Gh} and~\eqref{rk2_Bh} with parameter values given in \eqref{rk2_optimal_parameters}, the function $\mathcal{E}(x,h)$ is negative definite for $h$ small enough, and hence, the Metropolis integrator operated with this choice of $G_h(x)$ and $B_h(x)$ is second-order deterministically accurate.  A sufficient condition for this statement to be true is introduced in \S\ref{second_order}.   The following example illustrates these ideas, and in particular, the sufficiency of this condition.

\begin{figure}[Ht!]
\centering
\includegraphics[width=0.49\textwidth]{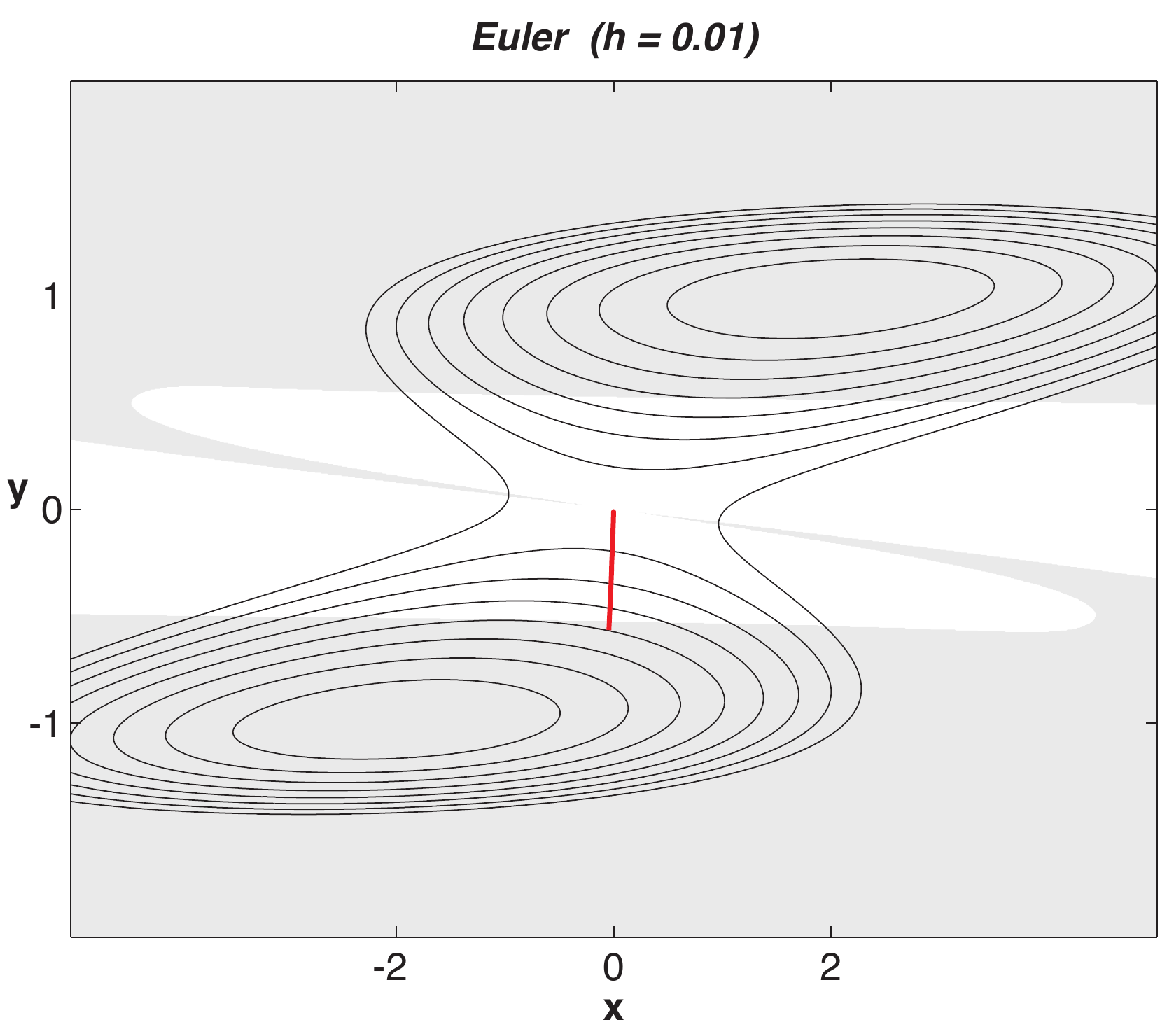}  
\includegraphics[width=0.49\textwidth]{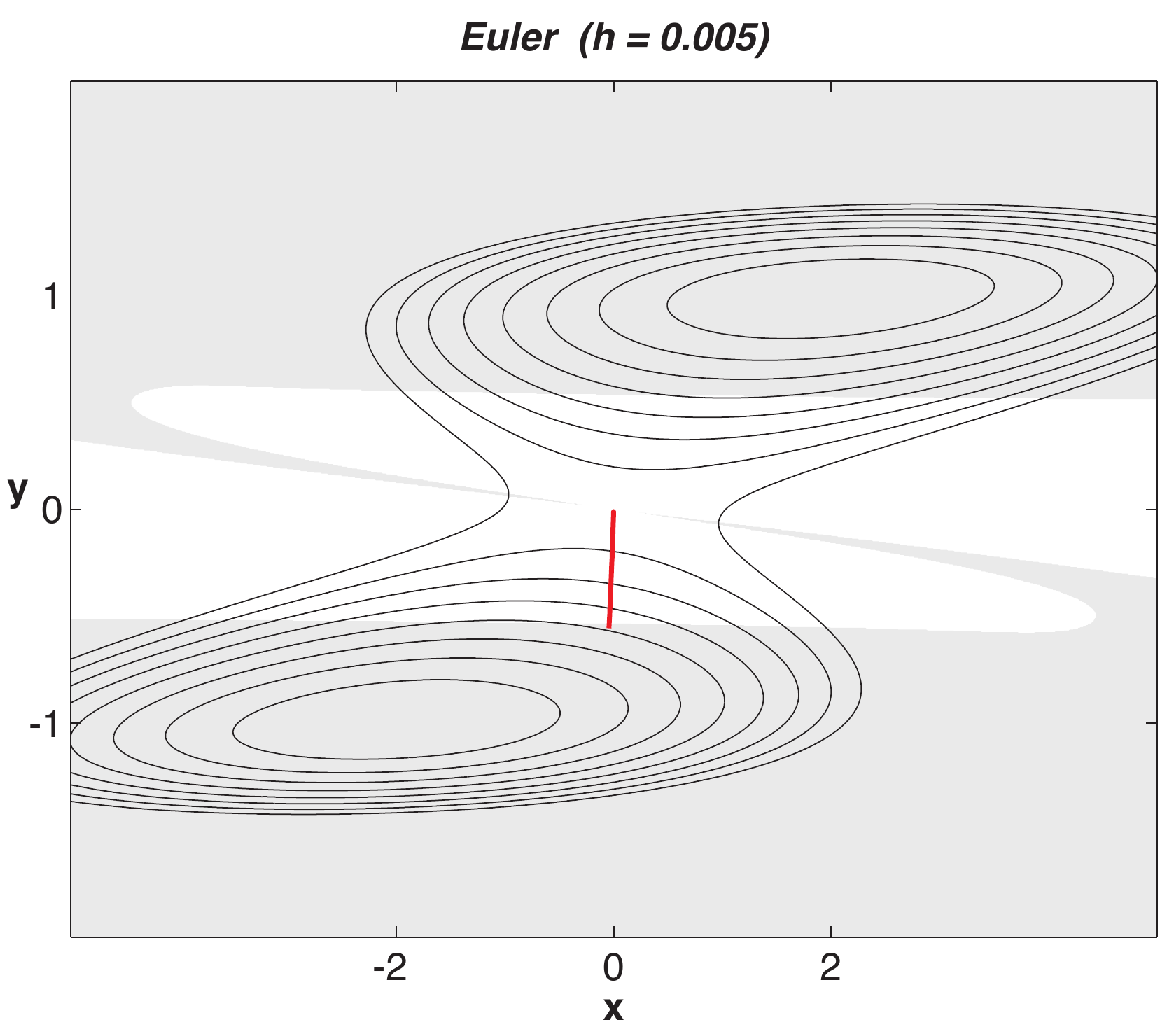}   \\
\hbox{
           \hspace{0.75in} (a) $G_h(x)=$ \eqref{eulerG}
        \hspace{1.5in} (b) $G_h(x)=$ \eqref{eulerG}
        } 
        
        \bigskip
        
\includegraphics[width=0.49\textwidth]{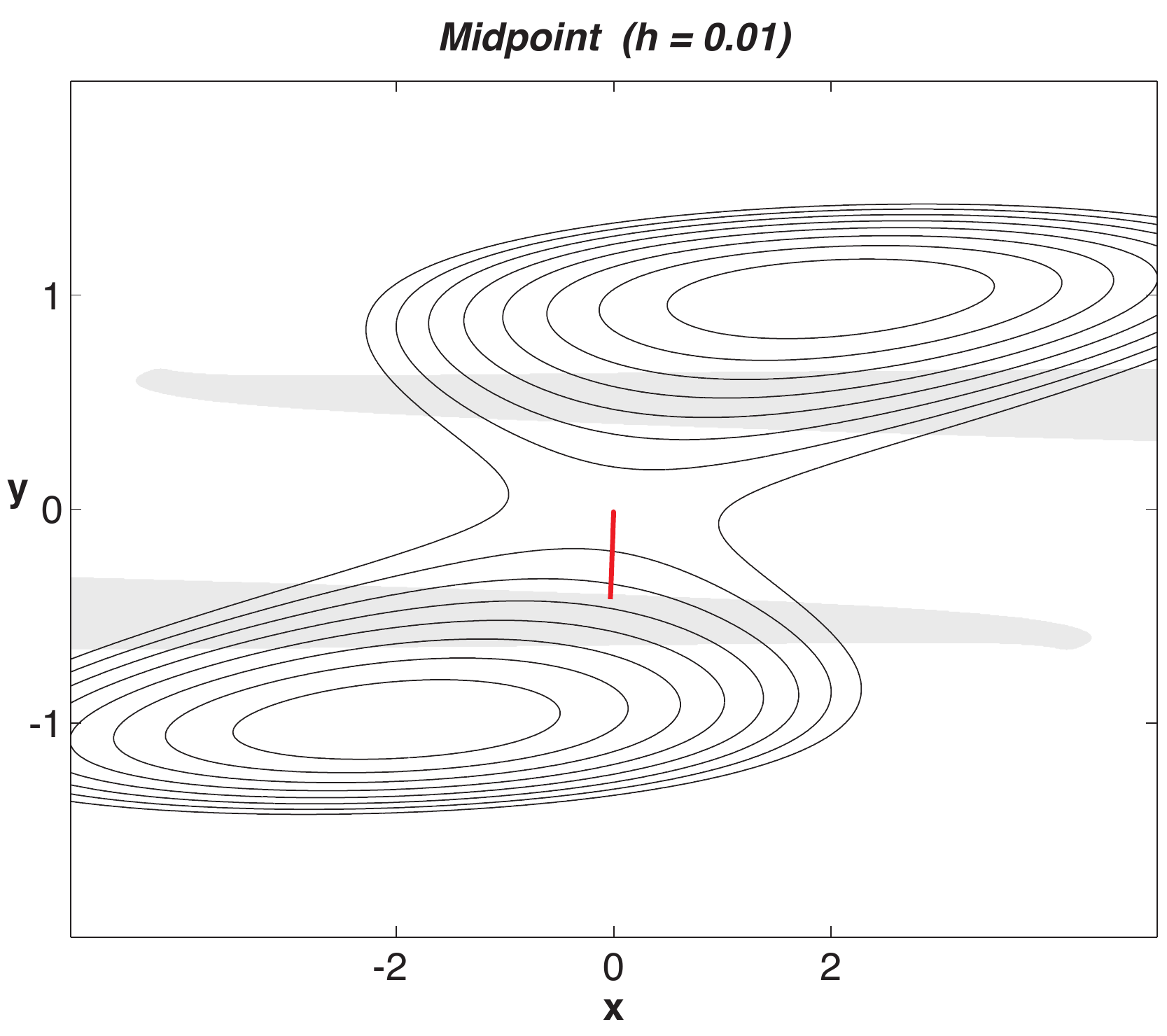}    
\includegraphics[width=0.49\textwidth]{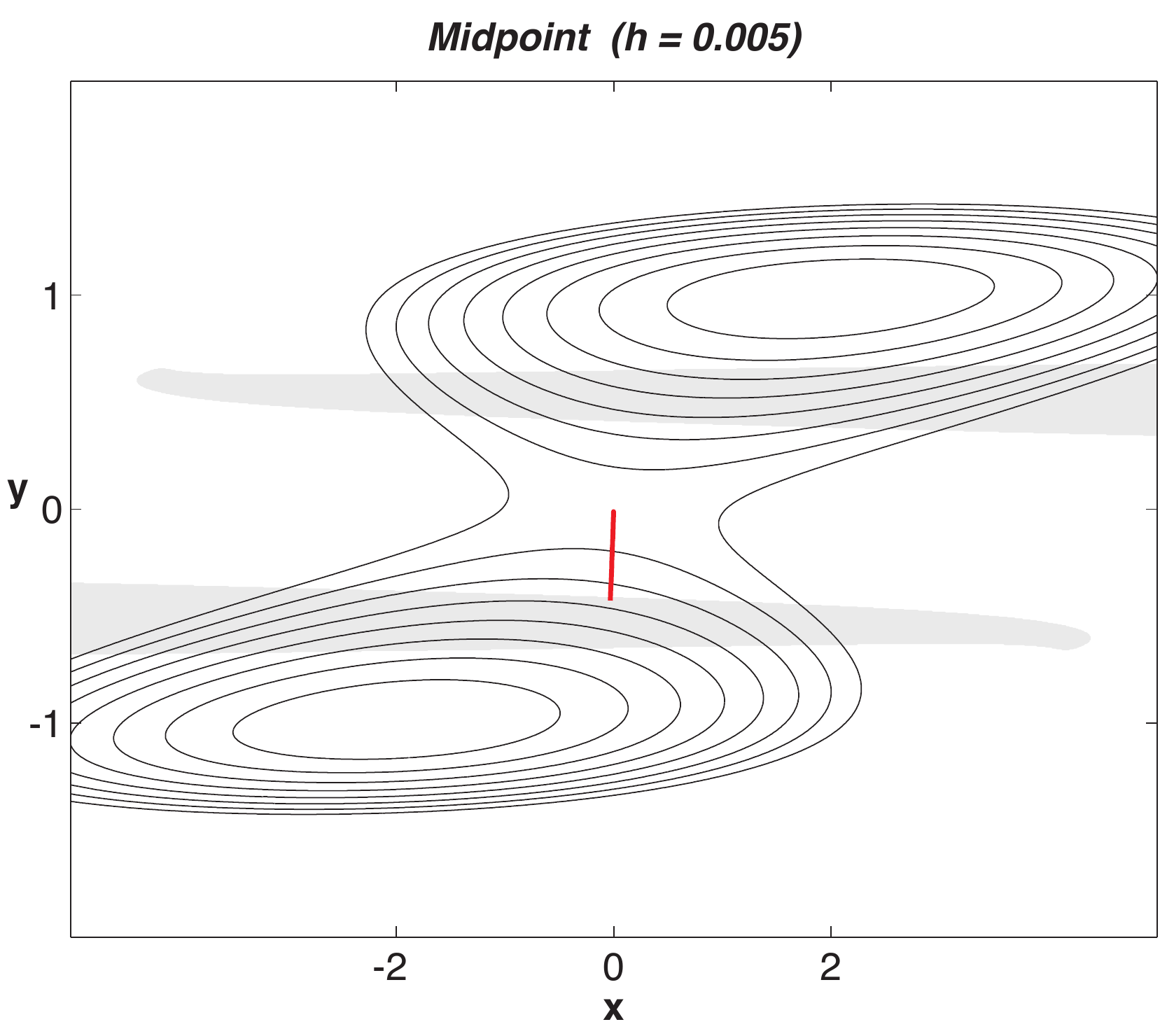}  \\
\hbox{
           \hspace{0.75in} (c) $G_h(x)=$ \eqref{midpointG}
        \hspace{1.5in} (d) $G_h(x)=$ \eqref{midpointG}
        } 
        
        \bigskip
       
\caption{ \small {\bf Deterministically Inaccurate Metropolis Integrators.}
The shaded regions in (a) - (d) are areas where the $\mathcal{E}(x,h)$ function of the Metropolis integrator is positive. In the background are contours of the potential energy in black.   
Observe that the $\mathcal{E}(x,h)$ functions of the Euler and Midpoint schemes are positive in a region which persists even if the time step size is halved.
Thus, sample trajectories produced by these schemes  -- starting at $(0, -0.01)^T$ with $\beta=10^{8}$ -- incorrectly terminate once they enter the shaded regions even when the time step size is reduced.  The  initial condition $(0, -0.01)^T$ is selected rather than $(0,0)^T$ for visualization purposes, and specifically, so that trajectories are more likely to move to the bottom well rather than the top one.
}
\label{fig:dw2da}
\end{figure}

\begin{figure}[Ht!]
\centering
 \includegraphics[width=0.49\textwidth]{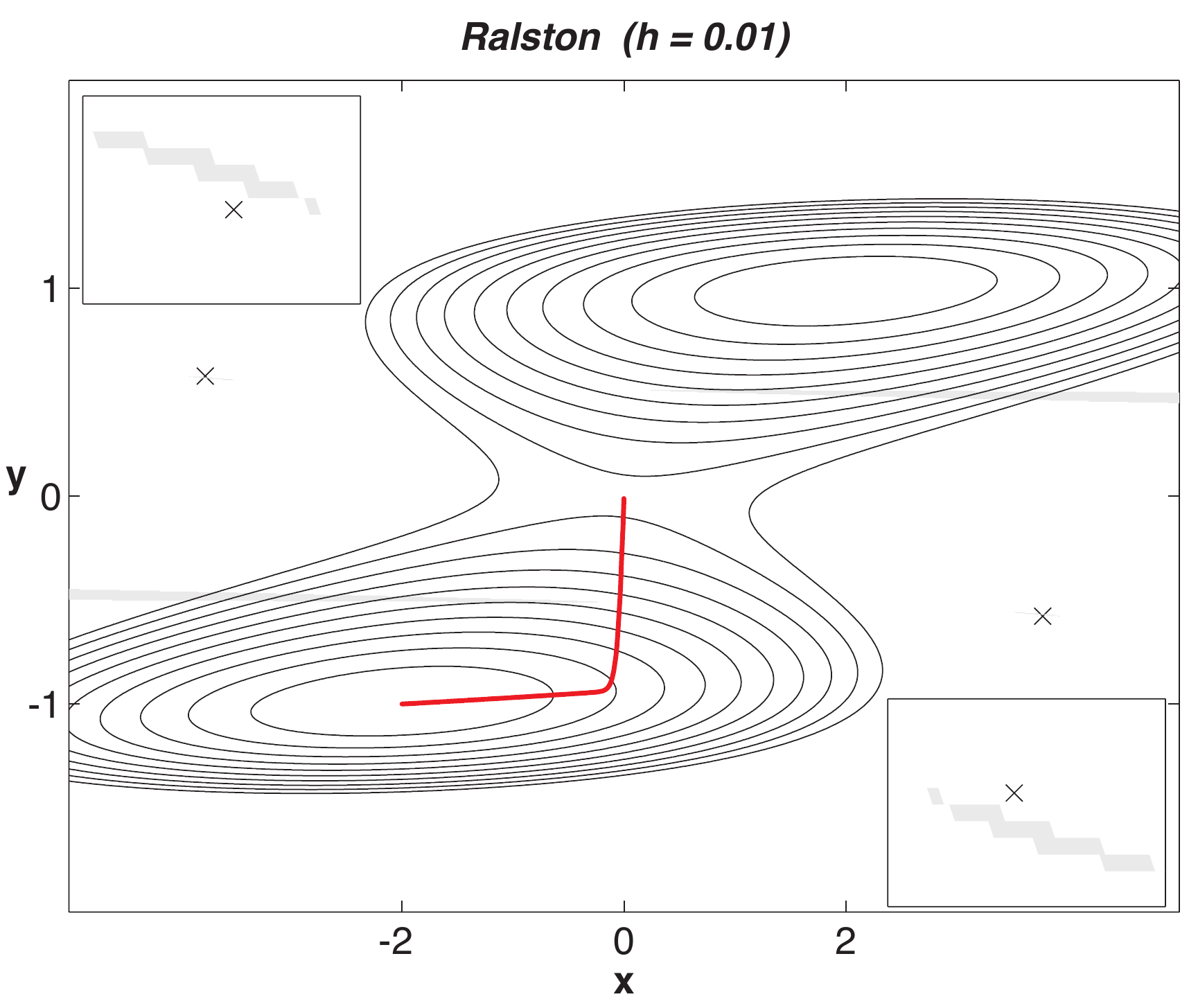}  
\includegraphics[width=0.49\textwidth]{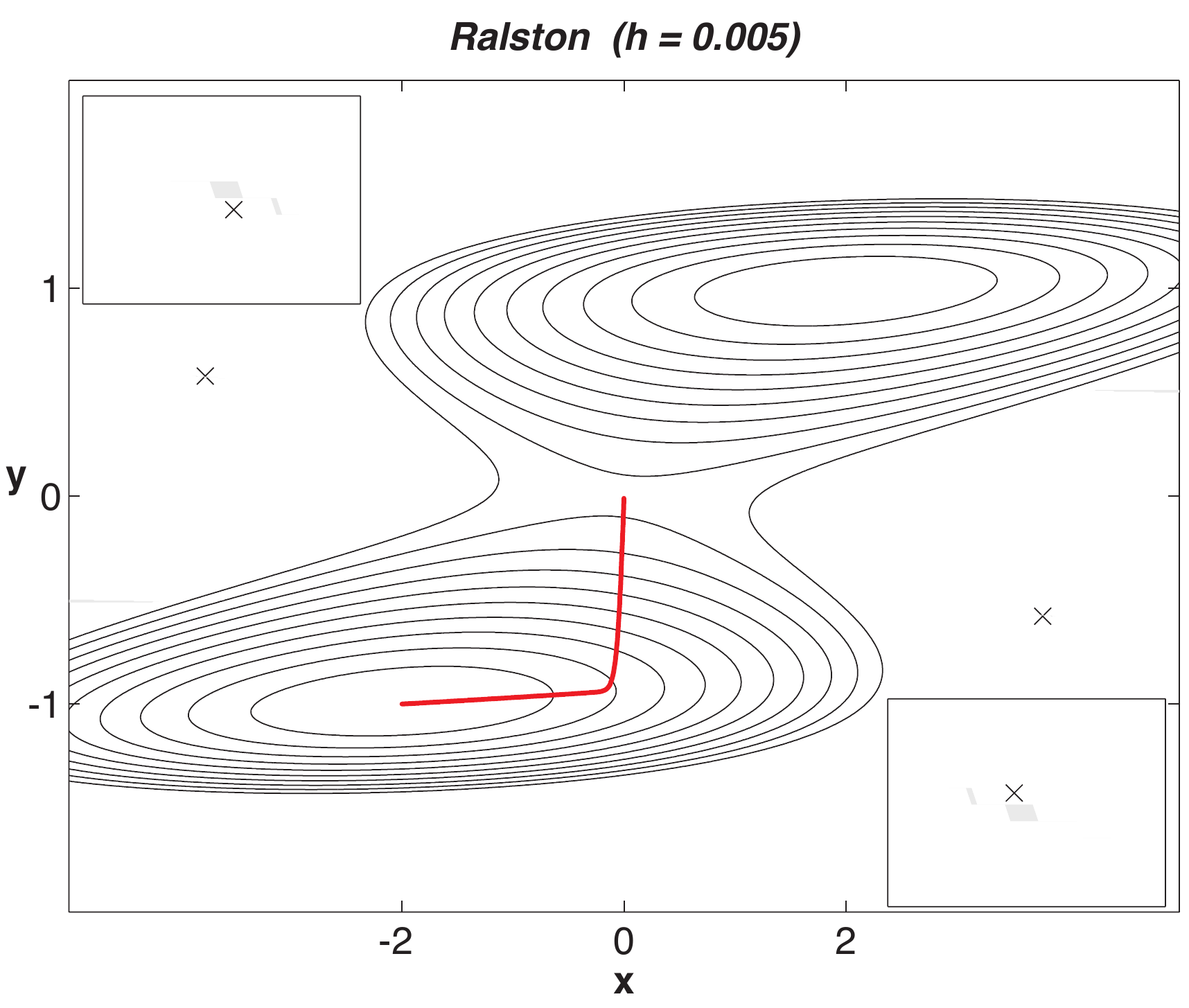}   \\
\hbox{
           \hspace{0.75in} (a) $G_h(x)=$ \eqref{e4ralstonG}
        \hspace{1.5in} (b) $G_h(x)=$ \eqref{e4ralstonG}
        } 
        
        \bigskip
        
\includegraphics[width=0.49\textwidth]{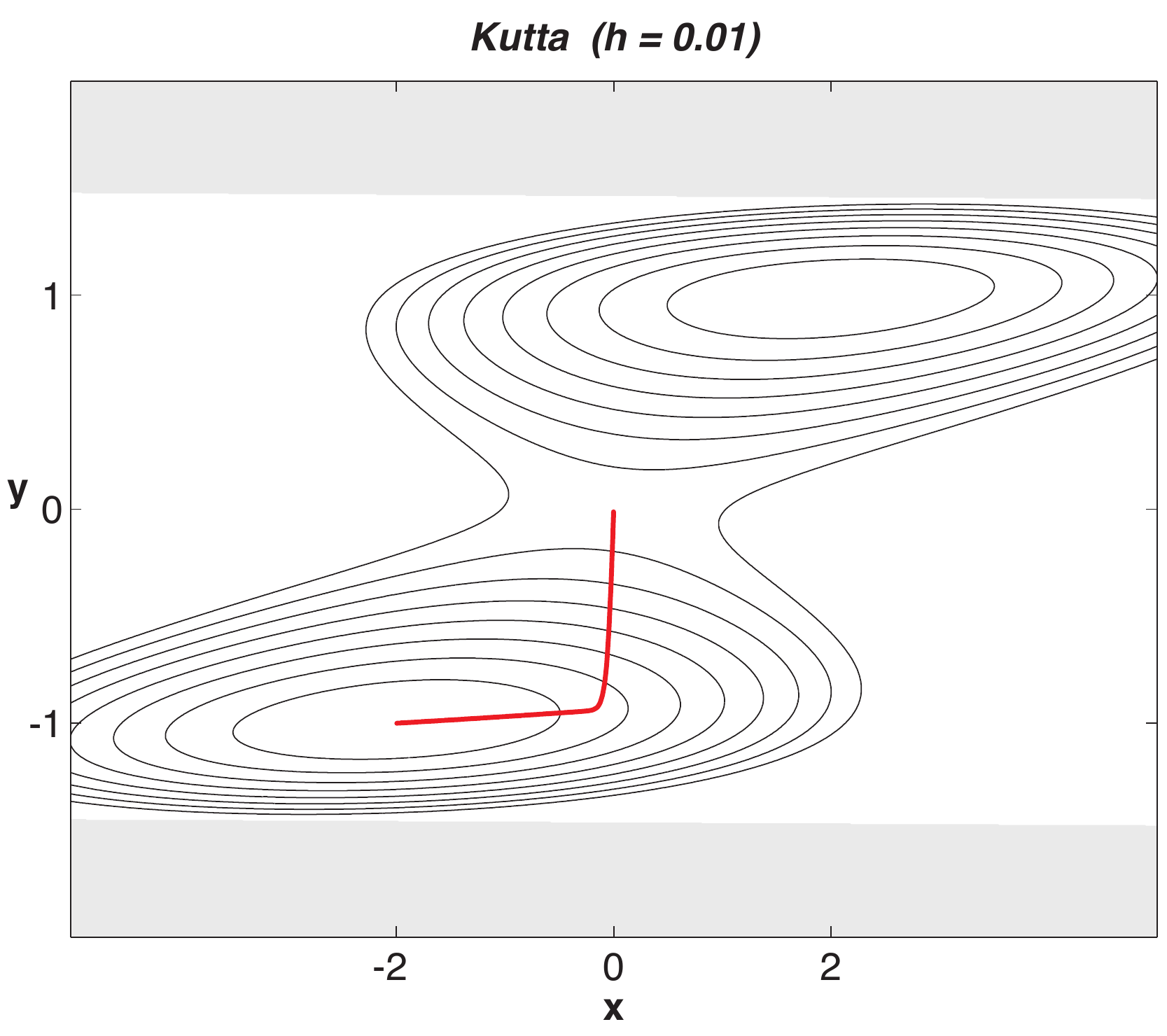}    
\includegraphics[width=0.49\textwidth]{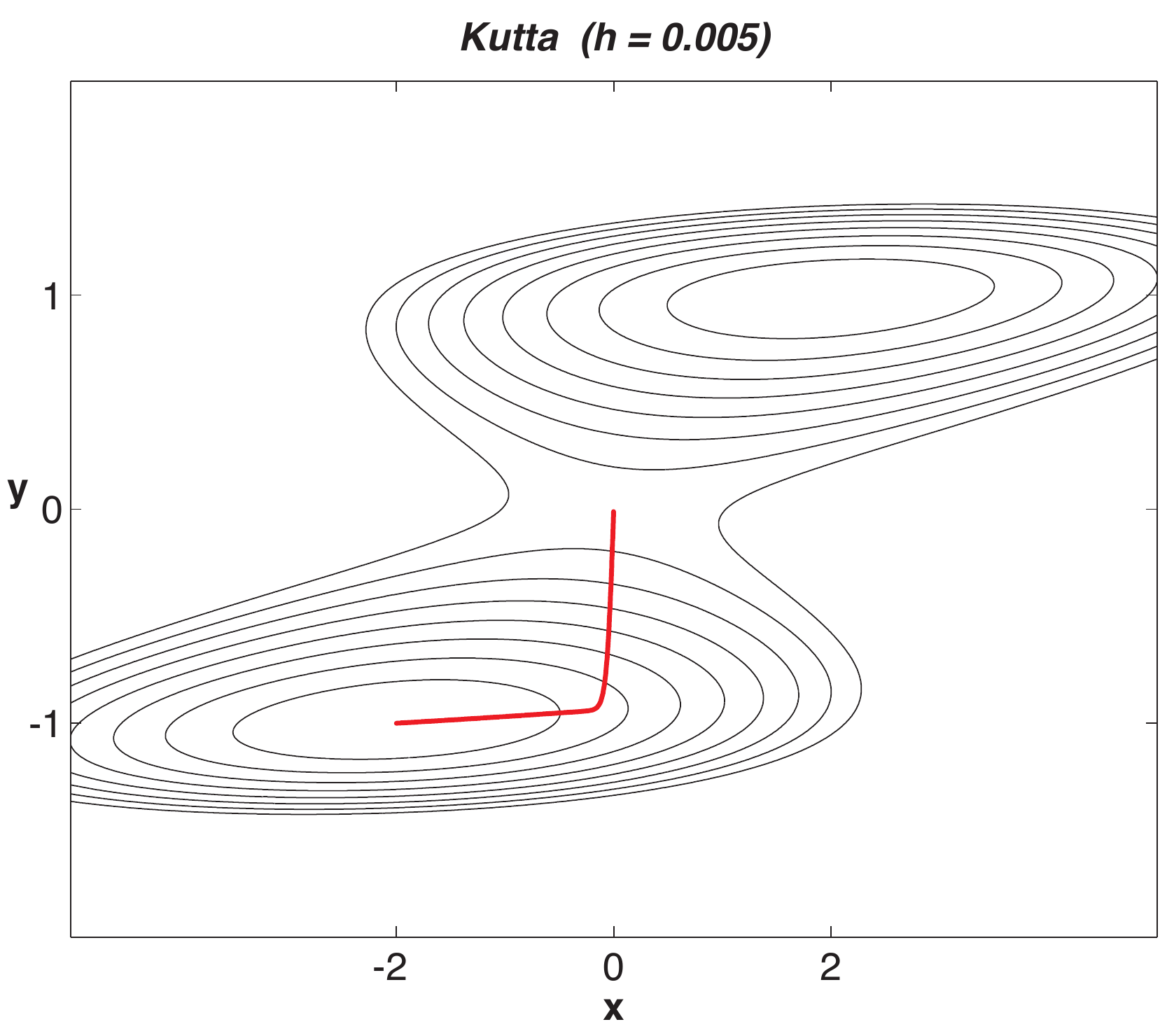}  \\
\hbox{
           \hspace{0.75in} (c) $G_h(x)=$ \eqref{kuttaG}
        \hspace{1.5in} (d) $G_h(x)=$ \eqref{kuttaG}
        } 
        
        \bigskip
       
\caption{ \small {\bf Deterministically Accurate Metropolis Integrators.}
Gray shading represents regions where the $\mathcal{E}(x,h)$ function of the Metropolis integrator is positive. In the background are contours of the potential energy in black.   The singular points of the Hessian matrix of $U(x)$ are x-marked in the main figure of (a) and (b).  The insets in (a) and (b) zoom in on these points, revealing that the $\mathcal{E}(x,h)$ function of the Ralston scheme can be positive in neighborhoods of these singular points.  Moreover, the insets show that when the time step size is reduced in the Ralston scheme, these gray regions become smaller.  The Kutta scheme -- being third-order accurate -- does not have a problem at these points, confirming the theory provided in \S\ref{sec:deterministicaccuracy}.  As shown in (a)-(d), sample trajectories -- starting at $(0, -0.01)^T$ with $\beta=10^{8}$ -- are accurately represented by both schemes.
}
\label{fig:dw2db}
\end{figure}

Consider a Brownian particle with the following two-dimensional double-well potential energy function: \[
U(x)=5 (x_2^2-1)^2 + 1.25 \left( x_2 - \frac{x_1}{2} \right)^2 \;, \qquad  x = (x_1, x_2)^T 
\] and a mobility matrix set equal to the $2 \times 2$ identity matrix for simplicity.  
The Hessian of this potential energy becomes singular at the two points marked by an `x' in Figures~\ref{fig:dw2db} (a) and (b).  For this example, $G_h(x)$ in~\eqref{rk2_Gh} with parameter values given in~\eqref{rk2_optimal_parameters} simplifies to the Ralston Runge-Kutta combination: \begin{equation} \label{e4ralstonG}
\begin{cases}
G_h(x) = - \frac{1}{4} DU(x) - \frac{3}{4} D U (\tilde x)  \;, \\
 \tilde x = x - \frac{2}{3} h DU(x) \;.
 \end{cases}
\end{equation}
We will compare this choice against the following choices of $G_h(x)$: \begin{equation} \label{eulerG}
G_h(x) = - D U(x) 
\end{equation} \begin{equation} \label{midpointG}
G_h(x) = -  D U \left(x - \frac{h}{2} DU(x) \right) 
\end{equation} 
which correspond to a first-order forward Euler and second-order midpoint discretization of \eqref{sde} in the 
zero noise limit, respectively. We also consider a third-order accurate Kutta approximation 
given by: \begin{equation} \label{kuttaG}
\begin{cases}
G_h(x) = - \frac{1}{6} DU(x) - \frac{2}{3} D U(\tilde x)  - \frac{1}{6}  D U( \bar x  )   \;,  \\
\tilde x = x - \frac{h}{2} DU(x) \;, ~~ \bar x  = x + h DU(x) - 2 h DU(\tilde x) \;.
\end{cases}
\end{equation}
Figure~\ref{fig:dw2da} plots the $\mathcal{E}(x,h)$ function and a sample trajectory for the forward Euler and midpoint schemes.  In the gray-shaded regions, the $\mathcal{E}(x,h)$ function is positive, and therefore, if the noise is small enough, the sample trajectory of the Metropolis scheme will likely get stuck in these regions.  In the simulation we pick $\beta=10^{8}$. The Metropolis integrator based on $G_h(x)$ given by \eqref{eulerG} or \eqref{midpointG} stops accepting proposal moves once they hit the shaded region.

Figure~\ref{fig:dw2db} plots the $\mathcal{E}(x,h)$ function of the Ralston and Kutta schemes, which as expected are negative definite if $h$ is small enough.  The inset in Figures~\ref{fig:dw2db} (a) and (b) zoom into neighborhoods of the two points where the matrix $D^2 U(x)$ is singular, showing that the $\mathcal{E}(x,h)$ function of the Ralston scheme can be positive in these neighborhoods. The insets also show that these regions become smaller as the time step size is reduced.    (On the other hand, the $\mathcal{E}(x,h)$ function of the  midpoint-based Metropolis schemes can be positive in regions that do not shrink with decreasing time step size -- illustrating the optimality of the Ralston scheme among two-stage Runge-Kutta methods.)  The Kutta scheme does not appear to have a problem at these singular points of the Hessian.  

We have also considered the same system, but with the non-constant mobility matrix given by: \[
M(x) = \begin{bmatrix} x_1^2+x_2^2+1 & 0 \\ 0 & x_1^2+x_2^2+1 \end{bmatrix} \;.
\]  Experiments revealed that choosing $B_h(x) = B(x)$ leads to an $\mathcal{E}(x,h)$ function that is not  negative definite, whereas using $B_h(x)$ given by \eqref{rk2_Bh} gives the desired property.  Moreover, as shown in \S\ref{sec:deterministicaccuracy}, $G_h(x)$ given by a third-order accurate Runge-Kutta scheme \eqref{rk3_Gh}, and $B_h(x)$ correspondingly selected, as shown in \eqref{rk3_Bh}, also leads to a deterministically third-order accurate scheme.


 \setlength\tabcolsep{1pt}

 \begin{table}[Ht!]
\centering
\begin{tabular}{|c|c|c|}
\hline
\multicolumn{1}{|c|}{\bf Parameter} & \multicolumn{1}{|c|}{\bf Description} &  \multicolumn{1}{|c|}{\bf Value}   \\
\hline
\hline
\multicolumn{3}{|c|}{ {\em Physical Parameters} } \\
\hline
  \hline
$\beta$ & inverse temperature factor & $10^{8}$  \\
\hline
  \hline
\multicolumn{3}{|c|}{{\em  Numerical Parameters }} \\
\hline
\hline
$h$ & time step size & $\{0.01, 0.005 \}$  \\
\hline
$Y(0)$ & initial condition & $(0, -0.01)^T$ \\
\hline
\end{tabular} 
\caption{  {\bf Simulation Parameters for a Brownian Particle with a  2D Double Well Potential.}   }
\label{tab:2Ddoublewell}
\end{table}

 
%
%
%
%
%
%
%

\clearpage

\section{Ergodicity for Normalizable $\boldsymbol \nu(  \boldsymbol x)$}\label{sec:ergodicity}

In this section it is shown that Algorithm~\ref{MetropolisIntegrator} satisfies Property (P1) when the following assumptions hold.

\medskip

\begin{assumption} \label{P1_nu_assumption}
The function $\nu(x)$ is positive and continuous for all $x \in \mathbb{R}^n$, 
and integrable, i.e., $\int_{\mathbb{R}^n} \nu(x) dx < \infty$.   
\end{assumption}

\medskip

\begin{assumption} \label{P1_Gh_assumption}
The vector-valued function $G_h(x)$ is continuously differentiable for all $x \in \mathbb{R}^n$,
and for $h$ sufficiently small, the matrix $2 I_{n \times n} + h D G_h(x)$ is invertible for all $x \in \mathbb{R}^n$,
where $I_{n \times n}$ is the $n \times n$ identity matrix.
\end{assumption}

\medskip

\begin{assumption} \label{P1_Bh_assumption}
The matrix-valued function $M_h(x) = B_h(x) B_h(x)^T$ is positive definite and continuous for all $x \in \mathbb{R}^n$.
\end{assumption}

\medskip

\noindent
We begin by showing that the method is $\nu$-symmetric by writing the algorithm first as a Metropolis-Hastings method, and second, as a hybrid Monte-Carlo (HMC) method.  By using tools from the theory of Metropolis-Hastings methods, we use the former viewpoint to prove that this stationary probability density is unique, and that the $k$-step transition probability of the Metropolis integrator converges to the stationary distribution with density $\nu(x)/Z$ as $k \to \infty$ for arbitrary initial distributions.

\subsection{Integrator as a Metropolis-Hastings algorithm}

Here it is shown that the Markov chain produced by Algorithm~\ref{MetropolisIntegrator} can be put in the frame of a Metropolis-Hastings method.   A key idea in this framework is the notion of Markov chain reversibility: a Markov chain with probability transition kernel $P(x,dy)$ is said to be {\em reversible} with respect to a measure $\mu(dx)$ if and only if the following equality of measures holds, \begin{equation}  \label{eq:gdb}
P(x,dy) \mu(dx) = P(y, dx) \mu(dy) \;.
\end{equation}
This condition is a generalization of the $\nu$-symmetry condition~\eqref{eq:db} to measures.  Indeed, setting $P(x,dy) = p_t(x,y) dy$ and $\mu(dy) = \nu(y) dy$, \eqref{eq:gdb} implies \eqref{eq:db}.  In other words, \eqref{eq:db} is a special case of \eqref{eq:gdb} when both the probability transition kernel and the measure $\mu(dx)$ have a common dominating measure like Lebesgue measure.  A reversible Markov chain automatically admits $\mu(dx)$ as an invariant measure since: \[
\int_{\mathbb{R}^n} \int_{\mathbb{R}^n} f(y) P(x,dy) \mu(dx) = \int_{\mathbb{R}^n} \int_{\mathbb{R}^n} f(y)  P(y, dx) \mu(dy) = \int_{\mathbb{R}^n} f(y) \mu(dy) 
\]
where $f(x)$ is an arbitrary test function.  

The Metropolis-Hastings algorithm constructs a Markov chain with a specified stationary distribution $\mu(dx)$ by enforcing condition \eqref{eq:gdb} at every step of the chain \cite{MeRoRoTeTe1953,Ha1970}.  The method is made up of a proposal move with probability transition kernel $Q(x,dy)$ and an acceptance probability $\alpha(x,y)$.  If the current state is $x$, the algorithm updates this state in two sub-steps: first, a proposal move is generated from $Q(x,dy)$; and second, this proposal move is accepted with probability $\alpha(x,y)$, and otherwise, the proposal move is rejected and the chain remains at $x$.   Thus, the  probability transition kernel of a Metropolis-Hastings chain can be written as: \begin{equation} \label{eq:mhkernel}
P(x,dy) = Q(x,dy) \alpha(x,y) + \delta_x(dy) \int_{\mathbb{R}^n} (1 - \alpha(x,z)) Q(x,dz) 
\end{equation}
where $\delta_x( dy )$ denotes the Dirac-delta distribution on $\mathbb{R}^n$ concentrated at the current state $x \in \mathbb{R}^n$.  
The Metropolis-Hastings method requires that $Q(x,dy)$ and $\alpha(x,y)$ are selected so that $P(x,dy)$ is reversible with respect to $\mu(dx)$.   The necessary and sufficient condition on the acceptance probability and the proposal kernel for reversibility to hold is: \begin{equation}
Q(x,dy) \alpha(x,y) \mu(dx) = Q(y,dx) \alpha(y,x) \mu(dy) 
\end{equation}
which is an identity statement about measures \cite{tierney1998note}.  In the case where $Q(x,dy) = q(x,y) dy$ and $\mu(dy) = \nu(y) dy$, it can be shown that this requirement is fulfilled by: \begin{equation} \label{eq:mhratio}
\alpha(x,y) = \min\left( 1, \frac{q(y,x) \nu(y) }{ q(x,y) \nu(x)} \right) \;.
\end{equation}
A quick glance at Algorithm~\ref{MetropolisIntegrator} suggests that the scheme is a Metropolis-Hastings method, but in order to establish this, we derive the probability transition density of the proposal move in \eqref{proposal}, and use this density to obtain the acceptance probability in \eqref{alphah} from \eqref{eq:mhratio}.   Once we have written the Metropolis integrator as a Metropolis-Hastings method, it immediately follows that it is reversible with respect to $\mu(dy) = \nu(y) dy$, and with a few additional steps, it can be shown that this stationary density is unique, and that the scheme is ergodic.   The proof of the following theorem fleshes out the details in this sketch.

\medskip

\begin{theorem} \label{thm:ergodicityofscheme}
Suppose that Assumptions~\ref{P1_nu_assumption}, \ref{P1_Gh_assumption}, and \ref{P1_Bh_assumption} hold.  Let $P_h(x,dy)$ denote the probability transition kernel of the integrator Algorithm~\ref{MetropolisIntegrator}.  Then the Markov chain induced by Algorithm~\ref{MetropolisIntegrator} preserves the equilibrium distribution $\mu(dx) = (\nu(x)/Z) dx$, and the k-step transition probability distribution of the scheme converges to $\mu$ in the total variation norm: \begin{equation} \label{ergodicity}
\lim_{k \to \infty} \| P_h^k(x, \cdot) - \mu \|_{\TV} = 0\;, \qquad \text{for all $x \in \mathbb{R}^n$} 
\end{equation}
for every $x \in \mathbb{R}^n$ and $h$ sufficiently small.
\end{theorem}

\medskip

\noindent
Here, we have used the total variation (TV) distance between two probability measures, which is defined as
\begin{equation}
\|\mu_1 - \mu_2 \|_\TV = 2 \; \sup_A |\mu_1(A) - \mu_2(A)|
\end{equation}
where the supremum ranges over all measurable sets. In particular, the TV distance between two probability measures is two if and only if they are mutually singular.

\medskip

\begin{proof}  Let $q_h(x,y)$ denote the transition probability density of the proposal move \eqref{proposal}.  Consider the transformation $y = \psi_x(z)$: \begin{equation} \label{psi}
 \psi_x(z) =  2 z - x + h G_h(z) 
\end{equation}
where $x$ (the current state) is a fixed parameter.   Let $\chi_x$ be the inverse transformation defined implicitly via \begin{equation} \label{psi_inverse}
y + x = 2 \chi_x(y) + h G_h( \chi_x(y)) \;.
\end{equation}   By the inverse function theorem, since $G_h(x)$ satisfies Assumption~\ref{P1_Gh_assumption}, this inverse map exists and is continuously differentiable.  Assumption~\ref{P1_Bh_assumption} implies that the distribution of $\tilde X_1$ is a Gaussian measure with full support, and since $X_1^{\star} = \psi_{X_0}(\tilde X_1)$, the transition probability distribution of the Metropolis integrator is the pullback by $\chi_x$ of a multivariate Gaussian measure with mean $X_0$ and covariance matrix $\Sigma = (\beta^{-1} h)/2  M_h(X_0)$, with transition density $q_h(x,y)$ given by:
 \begin{align} \label{proposaldensity}
 & q_h(x, y) = \frac{( \pi h \beta^{-1})^{-n/2} }{ \det(B_h(x))}    \exp\left(  - \frac{\beta}{h}  (  \chi_x(y) - x)^T M_h(x)^{-1} (  \chi_x(y) - x ) \right)  | \det(D \chi_x( y )  ) |   \;.
\end{align}  
Given $q_h(x,y)$ in \eqref{proposaldensity} and $\nu(x)$ in \eqref{eq:nu},  the acceptance probability~\eqref{alphah} can be derived from the Metropolis-Hastings ratio~\eqref{eq:mhratio} to obtain: \begin{equation} \begin{aligned} \quad & \alpha_h  (x,y) = 1 \wedge \frac{q_h(y,x) \nu(y) } 
		       {q_h(x,y) \nu(x) }  = 1 \wedge \frac{\det(B_h(x))}{\det(B_h(y))}   \\
	& \times \exp\left( - \beta \left( \frac{(z-y)^T M_h(y)^{-1} (z-y) - (z-x)^T M_h(x)^{-1} (z-x)}{h} + U(y) - U(x) \right) \right)  \;, \\
\end{aligned} \end{equation}
where $z = \chi_x(y) = \chi_y(x)$.  Note that the identity $ \chi_x(y) = \chi_y(x)$ follows from \eqref{psi_inverse} and also implies that the Jacobian determinant of the inverse transformation $\chi_x$ appearing in \eqref{proposaldensity} drops out of this ratio.

The transition probability kernel of the algorithm can now be written in the Metropolis-Hastings form~\eqref{eq:mhkernel}: \begin{equation}
P_h(x,d y) = p_h(x,y) d y + r_h(x) \delta_x(dy) 
\label{mhkernel}
\end{equation}
where $p_h(x,y)$ is the off-diagonal transition density, \[
p_h(x,y) = q_h(x,y)  \alpha_h(x,y) 
\]  and $r_h(x)$ is the probability of remaining at the same point, \[
r_h(x) = 
1 - \int_{\mathbb{R}^n} p_h(x,z) dz  \;.
\]
Thus, we conclude that the Metropolis integrator, Algorithm~\ref{MetropolisIntegrator}, is reversible, and therefore, preserves the stationary distribution $\mu(dx) = (\nu(x)/Z) dx$.

We now turn to proving that the $k$-step transition probability of the integrator converges to this stationary distribution in the limit $k \to \infty$.
The hypotheses of the theorem imply that the stationary probability density $\nu(x)/Z$ and the transition density $q_h(x,y)$ in \eqref{proposaldensity} are strictly positive and continuous everywhere.  Hence, the algorithm is irreducible with respect to Lebesgue measure and aperiodic; see,  e.g., Lemma 1.2 of~\cite{MenTw1996} and references therein.   
According to Corollary 2 of \cite{Ti1994},  a Metropolis-Hastings algorithm that is irreducible with respect to the same measure it is designed 
to preserve is positive Harris recurrent.  Consequently, the algorithm is irreducible, aperiodic, and positive Harris recurrent.
According to Proposition 6.3 of \cite{Nu1984}, these properties are equivalent to ergodicity of the chain.
\end{proof}

\medskip

While Theorem~\ref{thm:ergodicityofscheme} assumes that $\nu(x)$ is normalizable, neither the $\nu$-symmetry condition~\eqref{eq:db}  -- which is a special case of \eqref{eq:gdb} -- nor the Metropolis-Hastings algorithm require that $\nu(x)$ is normalizable.   In fact, the identity \eqref{eq:gdb} does not require that $\mu(dx)$ is a probability measure, and is automatically satisfied by a Metropolis-Hastings chain with  probability transition kernel  \eqref{eq:mhkernel}. 

\subsection{Integrator as an HMC algorithm}

The HMC method is a general technique to sample from a given probability distribution with a density function \cite{RoDoFr1978, DuKePeRo1987, Li2008}; in our case this density is the stationary density of the diffusion $\nu(x)/Z$.    The algorithm is set in an extended position-velocity space, $\{  (x,v) \in \mathbb{R}^n \times \mathbb{R}^n \}$, where positions $x \in \mathbb{R}^n$ belong to the domain of $\nu(x)$, and velocities $v \in \mathbb{R}^n$ are an auxiliary variable.   On this extended space, an extended density, $\nu_{\text{extended}}(x,v)$, is introduced which is a product of $\nu(x)/Z$ and a Gaussian density in velocities (which may depend on the position $x$), and such that, the marginal density in position space is $\nu(x)/Z$.  The algorithm produces a Markov chain on $\mathbb{R}^n$ with the required stationary distribution by using a volume-preserving and time-reversible map on the extended space $\theta: \mathbb{R}^n \times \mathbb{R}^n  \to \mathbb{R}^n \times \mathbb{R}^n$ and iterating the following update rule.  We stress that this map does not have to be a volume-preserving discretization of the Hamiltonian dynamics associated with the free energy of $\nu_{\text{extended}}(x,v)$.

\medskip

\paragraph{HMC Algorithm} The one-step $X_0 \mapsto X_1$ update in HMC is as follows.
\begin{itemize}
\item Draw $V_0$ from the measure $\nu_{\text{extended}}(X_0, v) dv$.
\item Set $(X_1^{\star}, V_1^{\star}) = \theta(X_0, V_0)$. 
\item Take as actual update: \[
X_1 = \gamma X_1^{\star} + (1-\gamma) X_0 
\] where $\gamma$ is a Bernoulli random variable with parameter: \[
\alpha(X_0, V_0) = 1 \wedge \frac{\nu_{\text{extended}}(X_1^{\star}, V_1^{\star}) }{\nu_{\text{extended}}(X_0, V_0) } \;.
\]
\end{itemize}

\medskip

\noindent
The above description of the HMC algorithm is a slight adaptation of the usual description; see, e.g., the description given in a recent tutorial introduction to HMC \cite[Section 9]{Sa2014}.   Note that the updated velocity is discarded by the algorithm.  The Markov chain produced by HMC preserves $\nu(x)/Z$ and is $\nu$-symmetric \cite{Sc1999,CaLeSt2007}.

The Metropolis integrator, Algorithm~\ref{MetropolisIntegrator}, can be written as an HMC algorithm by defining: \begin{equation} \label{eq:extendednu}
\nu_{\text{extended}}(x,v) = \frac{(2 \pi \beta^{-1})^{-n/2}}{ \det(B_h(x))} \exp\left( - \beta \frac{v^T M_h(x)^{-1} v}{2} \right) \frac{\nu(x)}{Z} 
\end{equation} 
and by setting $\theta(x,v)$ in the HMC algorithm equal to a position Verlet discretization of the following extended dynamics:
\begin{equation} \label{eq:extended}
\begin{cases}
\dot{Q} = V  \\
\dot{V} = G_h(Q) 
\end{cases}
\end{equation}
where $G_h: \mathbb{R}^n \to \mathbb{R}^n$ is the function appearing in Algorithm~\ref{MetropolisIntegrator}. As expected the marginal probability density function of $\nu_{\text{extended}}(x,v)$ in the original position space is the correct one: $\nu(x)/Z$.  A position Verlet discretization of \eqref{eq:extended} is given explicitly by \cite{HaLuWa2010}: \begin{equation} \label{Position_Verlet}
\begin{cases}
Q_{1/2} = Q_0 + \frac{\delta t}{2} \; V_0  \\
V_1 = V_0 + \delta t \; G_h(Q_{1/2} )  \\
Q_1 = Q_{1/2} + \frac{\delta t}{2} \; V_1 
\end{cases}
\end{equation}
where $\delta t$ is an artificial time-step size that is related to the physical time-step size via $\delta t = \sqrt{2 h}$.  Notice that $\theta(x,v)$ is not a volume-preserving discretization of the Hamiltonian dynamics associated with the free energy of $\nu_{\text{extended}}(x,v)$.  In fact, this Hamiltonian function is non-separable, and volume-preserving discretizations of non-separable Hamiltonian systems are implicit in general \cite{HaLuWa2010}.  In addition to this issue of implicitness, the resulting vector fields will involve derivatives of the mobility matrix whose evaluation we seek to avoid.  In contrast, the artificial dynamics \eqref{eq:extended} can be readily discretized using an explicit method that satisfies the conditions that HMC assumes on $\theta(x,v)$, namely time-reversibility and incompressibility.   We emphasize that nowhere in this derivation did we assume a specific form of $G_h(x)$ or $B_h(x)$, and that by casting the integrator as an HMC algorithm, it immediately follows that the Metropolis integrator preserves the probability density $\nu(x)/Z$ and is $\nu$-symmetric.  

\medskip

In the next section, we show how to control the rejection rate in this algorithm, and discuss the benefits of using  an integration scheme that exactly preserves the $\nu$-symmetry condition~\eqref{eq:db}.

\section{Weak Accuracy at Constant Temperature}\label{sec:accuracy}

By using the theory of Monte-Carlo methods, the previous section showed that the Metropolis integrator satisfies the $\nu$-symmetry condition \eqref{eq:db} by design.  For the diffusion process $Y(t)$ that solves \eqref{sde}, this $\nu$-symmetry property is equivalent to the self-adjoint property \eqref{eq:idb}, which imposes a fairly strong constraint on the dynamics.  To see this, expand the generator of the process $Y(t)$ acting on a test function $f(x)$ in~\eqref{eq:generator} as follows: \begin{align} \label{generator}
 (L f )( x ) =  ( - M(x)DU(x) + \beta^{-1} \divergence M(x) )^T D f(x)  +\beta^{-1} \; \tr( M(x) D^2 f(x) ) 
 \end{align}   and consider another diffusion process $\tilde Y(t)$ driven by the same noise, yet with a different drift: \begin{equation} \label{eq:tildeYsde}
d \tilde Y = a(\tilde Y) dt + \sqrt{2 \beta^{-1}} B(\tilde Y) d W \;, \qquad \tilde Y(0)  \in \mathbb{R}^n \;.
\end{equation}
Let $\tilde L$ represent the generator of $\tilde Y(t)$ whose action on a test function is given by: \[
(\tilde L f)(x) =  a(x)^T D f(x) + \beta^{-1} \; \tr( M(x) D^2 f(x) ) \;. 
\]
 If the generator of this diffusion is self-adjoint with respect to $\nu(x)$, then its drift is {\em uniquely determined}: \begin{equation} \label{eq:uniquedrift}
a(x) = - M(x) DU(x) + \beta^{-1} \; \divergence{M}(x) 
\end{equation}
and is identical to the drift appearing in \eqref{sde} \cite{HaPa1986}.  
This constraint on the dynamics seems to motivate using approximations which preserve the $\nu$-symmetry property of the continuous process, like the Metropolis integrator.   We must consider, however, what precisely is gained -- if anything -- from a discretization that satisfies $\nu$-symmetry especially when $\nu(x)$ is not normalizable.   

For this purpose, consider a process $\tilde X(t)$ that satisfies the SDE: \begin{equation} \label{eq:tildeXsde}
d \tilde X = \tilde a( \tilde X ) dt + \sqrt{2 \beta^{-1}} B_h( \tilde X ) dW \;, \qquad \tilde X(0) \in \mathbb{R}^n 
\end{equation}
for all $t \in [0, h]$.  The infinitesimal generator $L_h$ of this process is given by: \[
(L_h f)(x) =  \tilde a(x)^T D f(x) + \beta^{-1} \tr( M_h(x) D^2 f(x) ) 
\]
where $M_h(x) = B_h(x) B_h(x)^T$.  Assuming that the noise in \eqref{eq:tildeXsde} is a single step approximation to the exact noise appearing in \eqref{sde} and that $L_h$ is self-adjoint with respect to $\nu(x)$, then $\tilde a(x)$ is a single step approximation to the true drift in \eqref{sde}: \begin{equation} \label{eq:tildea}
\tilde a(x) = - M(x) DU(x) + \beta^{-1} \divergence{M}(x) + \mathcal{O}(h) \;.
\end{equation}  We interpret \eqref{eq:tildea} as saying that: 

\medskip

\begin{quote}
{\em noise accuracy is sufficient for accuracy of a $\nu$-symmetric integrator.}
\end{quote}

\medskip

\noindent
This result enables the design of integrators that satisfy Property (P4) -- avoid computing the divergence of the mobility matrix -- and is one of the 
main advantages of using a $\nu$-symmetric integrator to simulate a self-adjoint diffusion.   
 
To state this result precisely for Algorithm~\ref{MetropolisIntegrator},  we will make the following assumptions.

\medskip

\begin{assumption} \label{P2_U_Assumption}
The first two derivatives of $U(x)$ are bounded for all $x \in \mathbb{R}^n$.
\end{assumption}

\medskip

\begin{assumption} \label{P2_Gh_assumption}
The vector-valued function $G_h(x)$ and its first derivative are bounded for all $x \in \mathbb{R}^n$ and for $h$ sufficiently small.
\end{assumption}

\medskip

\begin{assumption} \label{P2_Bh_assumption}
The matrix-valued function $B_h(x)$, its inverse, first and second derivatives are bounded
for all $x \in \mathbb{R}^n$ and for $h$ sufficiently small.  Additionally, there exists 
a bounded matrix-valued function $C(x)$ with bounded first derivative such that, \[
B_h(x) B_h(x)^T = B(x) B(x)^T + C(x) h^{1/2}  
\]
for all $x \in \mathbb{R}^n$ and for $h$ sufficiently small.
\end{assumption}

\medskip

\noindent
The latter assumption requires that $B_h(x)$ is chosen so that the Brownian force appearing in \eqref{sde} is approximated to up to $\mathcal{O}(h^{1/2})$.  
Note that a bound on $B_h(x)$ automatically implies a bound on $M_h(x) = B_h(x) B_h(x)^T$, since in the $2$-norm $\| B_h(x) \| = \| M_h(x)  \|^{1/2}$.  
We emphasize that nowhere in these assumptions is there a similar requirement on $G_h(x)$, and just to be clear, 
$G_h(x)$ does not need to approximate any part of the drift appearing in \eqref{sde} to satisfy Assumption~\ref{P2_Gh_assumption}.

\medskip

\begin{theorem}[Weak Accuracy] \label{thm:infinitesimalfdt}
Let Assumptions~\ref{P2_U_Assumption}, \ref{P2_Gh_assumption}, and \ref{P2_Bh_assumption} hold.
 For every time interval~$T>0$,  there exists a $C(T)>0$ such that the Metropolis integrator, Algorithm~\ref{MetropolisIntegrator}, satisfies:  \[
| \E_x ( f(Y( \lfloor t/h \rfloor h))) - \E_x (f (X_{\lfloor t/h \rfloor }) ) | \le C(T) h^{1/2} 
\] 
for all $t \in [0, T]$, for every $x \in \mathbb{R}^n$,  for $h$ sufficiently small, and for sufficiently regular test function $f(x)$.
\end{theorem}

\medskip

\begin{proof}
In order to prove Theorem~\ref{thm:infinitesimalfdt}, we analyze the transition probability distribution of the Metropolis integrator.  To this end let us fix some notation.  Let $\mathcal{P}_t$ and $P_h$ respectively denote  the Markov operators of the true solution and Metropolis integrator, whose actions on a test function $f(x)$ are given by: \begin{equation} \label{eq:Ptf}
(\mathcal{P}_t f)( x ) = \E_{x} f(Y(t)) \;, \quad \text{with $Y(0)=x$} 
\end{equation}  \begin{equation} \label{eq:Phf}
(P_h f)(x) =  \E_{x} ( f(X_1^{\star}) \alpha_h(X_0, \xi) ) + f(X_0) \E_{x} (1 - \alpha_h (X_0, \xi) ) \;, \quad \text{with $X_0=x$} \;.
\end{equation}
To express the derivative of the mobility tensor, we adopt the following shorthand: \[
DM(x)(u,v,w) = \frac{\partial M_{ij}}{\partial x_k} u_k v_j w_i \;.
\] 
In this notation there are three vectorial inputs to the tensor of order three $DM(x)$.  Notice that the first input gives the direction of the derivative of $M(x)$.  Also notice that since $M(x)$ is symmetric, $DM(x)(u, v, w) = DM(x)(u, w, v)$.  When evaluated at two vectors, this third-order tensor returns a vector with $ith$ component given by: \[
(DM(x)(u,v))_i = \frac{\partial M_{ij}}{\partial x_k} u_k v_j \;.
\] With this notation consider the following approximation (in the $L^2$ sense) to the acceptance probability: \begin{equation} \label{tilde_alphah}
 \tilde \alpha_h(x, \xi) =   1 \wedge \exp\left( -\beta \Gamma_h(x, \xi) \right) 
\end{equation}
where we have introduced:
\begin{equation}   \label{Gammah}
\begin{aligned}
& \Gamma_h(x,\xi) =    \sqrt{2 h} ( DU(x) + M_h(x)^{-1} G_h(x) )^T B_h(x) \xi  \\
& \qquad +\sqrt{2 h}  \frac{\beta^{-1}}{2} \tr (M_h(x)^{-1} D M_h(x)( B_h(x) \xi ) )   \\
& \qquad -  \sqrt{2 h} \frac{1}{2} DM_h(x) (B_h(x) \xi, B_h(x)^{-T} \xi, B_h(x)^{-T} \xi )   \;.  
\end{aligned}
\end{equation}
Recall, $\xi\in\mathbb{R}^n$ denotes a Gaussian random vector with mean zero and covariance $\E( \xi_i \xi_j ) = \beta^{-1} \delta_{ij}$. 
Lemma~\ref{lem:tilde_alphah} proves that this approximation satisfies: \[
\tilde \alpha_h(x, \xi) = \alpha_h(x, \xi ) + \mathcal{O}(h) \;.
\]
Note that when the mobility matrix is constant the choice $G_h(x) = - M_h(x) DU(x) + \mathcal{O}(h)$ implies that  $\Gamma_h(x, \xi) = 0$.  Corollary~\ref{cor:infinitesimalfdt} builds on this result by showing that the Metropolis integrator operated with this choice of $G_h(x)$ is first order weakly accurate.  We emphasize that this choice of $G_h(x)$ is not required by Theorem~\ref{thm:infinitesimalfdt}. 

The proof focuses on deriving the following upper bound on the difference between $\mathcal{P}_h$ in~\eqref{eq:Ptf} and $P_h$ in \eqref{eq:Phf}:  \begin{equation} \label{eq:onestepestimate}
| (\mathcal{P}_h f)( x ) - (P_h f)( x ) | \le C h^{3/2} 
\end{equation} for both sufficiently small $h$ and regular test functions $f(x)$.  Standard results in numerical analysis for SDEs then imply the algorithm converges weakly 
on finite time intervals with global order $1/2$; see for instance \cite[Chapter 2.2]{MiTr2004}.  For a detailed treatment of the technical issues involved in 
extending this one-step error estimate to a global error estimate when the potential force appearing in \eqref{sde} is not globally Lipschitz, 
the reader is referred to \cite{BoHa2013}.

Using Taylor's theorem and \eqref{eq:Phf}, one can write:  \begin{equation} \label{Phfexpansion}
\E_{x}  f(X_1)  = f(x) + I_1 + I_2 + I_3 + I_4  + I_5  + I_6 
 \end{equation}
 where we have introduced \begin{align*}
I_1 &=  \E_{x} \left( D f(x) (X_1^{\star} - x) \right) \;, \\
I_2 &=  \frac{1}{2} \E_{x} \left( D^2 f(x) (X_1^{\star} - x, X_1^{\star} - x)  \right) \;, \\
I_3 &=  \E_{x} \left( D f(x) (X_1^{\star} - x) (\tilde \alpha_h(x, \xi) - 1) \right) \;, \\
I_4 &= \E_{x} \left(  D f(x) (X_1^{\star} - x)  ( \alpha_h(x, \xi)   - \tilde \alpha_h(x, \xi) ) \right)   \;,  \\
I_5 &=\frac{1}{2} \E_{x} \left( D^2 f(x) (X_1^{\star} - x, X_1^{\star} - x)  ( \alpha_h(x, \xi) - 1)   \right) \;,  \\
I_6 &= \frac{1}{2}  \E_{x}  \left( \alpha_h(x, \xi) \int_0^1 (1-s)^2  D^3 f(x + s (X_1^{\star} - x) ) (X_1^{\star} - x)^3 ds    \right) \;.
\end{align*}
Recall that $X_1^{\star}$ denotes the proposal move \eqref{proposal} with $X_0 = x$.
(Here we interpret $D^3 f(x) y^3$ as being the trilinear form $D^3 f(x)$ applied to the triple $(y,y,y)$.)
The terms $f(x)$, $I_1$,  $I_2$ and $I_3$ contribute to the weak accuracy of the method.   In what follows we describe how to treat these terms.   It is straightforward to use the hypotheses to show that the remaining terms $I_4$, $I_5$ and $I_6$ in \eqref{Phfexpansion} are $\mathcal{O}(h^{3/2})$.

Notice that $I_1$ simplifies to, \begin{align}
I_1 =  h \; D f(x)^T G_h(x) + \mathcal{O}(h^{3/2}) \;. \label{I1} 
\end{align}
The leading term in $I_2$ involves a quadratic form in $\xi$, and hence, \begin{align}
I_2 = h \; \beta^{-1} \;  \tr( M_h(x) D^2 f(x)  ) + \mathcal{O}(h^{2}) \;. \label{I2}
\end{align}
By referring to \eqref{generator}, let's take stock of the analysis so far.  The term $I_2$ contributes the Ito-correction term in \eqref{generator}.   The choice $G_h(x) = - M(x) DU(x) + \mathcal{O}(h)$ implies that the term $I_1$ accurately represents the contribution of the deterministic drift in \eqref{generator}.  However, this choice of $G_h(x)$ is insufficient for accuracy since the term involving $\divergence M(x)$ in \eqref{generator} has not been accounted for.   While this missing term can be added to $G_h(x)$, recall that an explicit formula for $\divergence M(x)$ is not typically available in practical BD problems.  This discussion foretells the importance of the term $I_3$.  In fact, it will be shown that -- independent of the precise form of $G_h(x)$ -- the sum of $I_1$ and $I_3$ represent the first order effect of the drift in \eqref{generator}.

Now we show how to estimate $I_3$.  
To simplify these calculations, introduce \begin{equation} \label{abc}
\begin{cases}
c(z_1, z_2, z_3) = \frac{1}{2}  \sqrt{2 h \beta^{-1}} DM_h(B_h z_1, B_h^{-T} z_2, B_h^{-T} z_3 )   \\
b(z) = - \sqrt{2 h \beta} ( DU + M_h^{-1} G_h )^T B_h z - \frac{\sqrt{2 h \beta^{-1}}}{2} \tr (M_h^{-1} D M_h( B_h z ) )  \\
a = \sqrt{2 h \beta^{-1}} B_h^T Df  
\end{cases}
\end{equation}
where $z \in \mathbb{R}^n$.  For the sake of readability, the dependence of $Df(x)$, $G_h(x)$, $M_h(x)$ and $B_h(x)$ on $x$ is suppressed.  
In terms of these variables, $I_3$ can be written as \begin{align*}
I_3 &=  (2 \pi)^{-n/2}  \int_{\mathbb{R}^n}   \exp\left( - \frac{|z|^2}{2} \right) \left( a^T z \right)  \left(1 \wedge  e^{ b(z) + c(z,z,z) } - 1 \right)dz + \mathcal{O}(h^{3/2}) \;, \\
&=  - (2 \pi)^{-n/2}  \int_{\mathbb{R}^n} \divergence_z \left( \exp\left(- \frac{|z|^2}{2} \right) a \right) \left(1 \wedge  e^{ b(z) + c(z,z,z) }  - 1\right)dz + \mathcal{O}(h^{3/2}) \;.
\end{align*}  Let $\Omega = \{ z \in \mathbb{R}^n \mid b(z) + c(z,z,z) \le 0  \}$, 
$\partial \Omega = \{ z \in \mathbb{R}^n \mid b(z) + c(z,z,z)=0\}$ be the boundary of $\Omega$, and $\vec{\nu}$ be the outward pointing normal to the region $\Omega$. 
Set $g_1(z) =  e^{ b(z) + c(z,z,z) } - 1$ and $\vec{g}_2(z) = \exp\left(- |z|^2/2 \right) a$, and in terms of which, consider the following multidimensional integration by parts formula: \begin{equation} \label{integration_by_parts}
\int_{\partial \Omega} g_1(z) \vec{g}_2(z)^T \vec{\nu} dz = \int_{\Omega} g_1(z) ( \divergence_z \vec{g}_2)(z) dz + \int_{\Omega} \vec{g}_2(z)^T (\nabla_z g_1)(z) dz \;.
\end{equation}
The boundary term in this equation vanishes since $g_1(z)$ vanishes on $\partial \Omega$.  Using this identity, $I_3$ can be written as \begin{align*}
I_3 &= - (2 \pi)^{-n/2} \int_{\Omega}  g_1(z)  ( \divergence_z  \vec{g}_2) (z)  dz + \mathcal{O}(h^{3/2}) \;,  \\
&=  (2 \pi)^{-n/2} \int_{\Omega}  e^{- \frac{|z|^2}{2} } ( b(a) + c(a,z,z)+2 c(z,a,z) )  \left( e^{ b(z) + c(z,z,z) }  \right)  dz + \mathcal{O}(h^{3/2}) \;.
\end{align*}
Since the map $z \mapsto b(a) + c(a, z,z) + 2 c(z,a,z)$ is symmetric (or even: $f(-z) = f(z)$)
and $z \mapsto b(z) + c(z,z,z)$ is anti-symmetric (or odd: $f(-z)=-f(z)$), it follows that \begin{align*}
I_3 &=   (2 \pi)^{-n/2} \int_{\Omega}  e^{- \frac{|z|^2}{2} } ( b(a) + c(a,z,z)+2 c(z,a,z) )  dz  + \mathcal{O}(h^{3/2}) \\
& \qquad +   (2 \pi)^{-n/2} \int_{\Omega}  e^{- \frac{|z|^2}{2} } ( b(a) + c(a,z,z)+2 c(z,a,z) )  (e^{ b(z) + c(z,z,z) }   -1) dz \;, \\
&= \frac{1}{2} (2 \pi)^{-n/2} \int_{\mathbb{R}^n}  e^{- \frac{|z|^2}{2} } ( b(a) + c(a,z,z)+2 c(z,a,z) )  dz   + \mathcal{O}(h^{3/2}) \\
& \qquad +  (2 \pi)^{-n/2} \int_{\Omega}  e^{- \frac{|z|^2}{2} } ( b(a) + c(a,z,z)+2 c(z,a,z) )  (e^{ b(z) + c(z,z,z) }   -1) dz \;, \\
&= \frac{1}{2} ( b(a) + c_{ijj} a_i +2 c_{iji} a_j )  + \mathcal{O}(h^{3/2}) \;.
\end{align*}
Substituting \eqref{abc} back into this  expression and simplifying yields, \begin{align}
I_3 = h (\beta^{-1} \divergence M_h - M_h DU - G_h)^T Df + \mathcal{O}(h^{3/2}) \;.  \label{I3}
\end{align}

Combining estimates \eqref{I1}, \eqref{I2} and \eqref{I3}, and invoking Assumption~\ref{P2_Bh_assumption}, yields the following estimate:  \[
(P_h f)(x) = f(x) + h (L f)(x) + \mathcal{O}(h^{3/2}) 
\]
which agrees with an Ito-Taylor expansion of $( \mathcal{P}_h f)(x)$ to up to terms of $ \mathcal{O}(h^{3/2})$, and hence, the scheme has the desired single step accuracy.
\end{proof}

\medskip

As a corollary of Theorem~\ref{thm:infinitesimalfdt}, the Metropolis integrator operated with  $G_h(x)$ given by 
the deterministic drift up to first order is weakly accurate when the mobility matrix is constant.

\medskip

\begin{cor} \label{cor:infinitesimalfdt}
Let Assumptions~\ref{P2_U_Assumption}, \ref{P2_Gh_assumption}, and \ref{P2_Bh_assumption} hold.  
Also assume that the mobility is constant and that $B_h = B$.
Consider the Metropolis integrator
with $G_h(x)$ given by: \begin{equation} \label{rk1_Gh}
G_h(x) = - M DU(x) + \mathcal{O}(h) 
\end{equation}
 For every time interval~$T>0$,  there exists a $C(T)>0$ such that:  \[
| \E_x ( f(Y( \lfloor t/h \rfloor h))) - \E_x (f (X_{\lfloor t/h \rfloor }) ) | \le C(T) h 
\] 
for all $t \in [0, T]$, for every $x \in \mathbb{R}^n$,  for $h$ sufficiently small, and for sufficiently regular test function $f(x)$.
\end{cor}

\medskip

\begin{proof}
The proof is identical to the proof of Theorem~\ref{thm:infinitesimalfdt}, but with \eqref{eq:onestepestimate} replaced with:  \begin{equation} \label{eq:onestepestimate_additive_noise}
| (\mathcal{P}_h f)( x ) - (P_h f)( x ) | \le C h^{2} 
\end{equation} for both sufficiently small $h$ and regular test functions $f(x)$.  To derive this estimate, note that 
$\Gamma_h(x,\xi)$ in \eqref{tilde_alphah} and $R_h(x,\xi)$ in \eqref{alphah_expansion} satisfy (in the $L^2$ sense): \begin{align*}
\Gamma_h(x, \xi) = \mathcal{O}(h^{3/2})  \;, \quad R_h(x, \xi)  = \mathcal{O}(h^{3/2})  
\end{align*}
when the mobility matrix is constant and $G_h(x)$ satisfies \eqref{rk1_Gh}.  These relations imply that: \[
1- \alpha_h(x,\xi)  =\mathcal{O}(h^{3/2}) \;, \quad   1- \tilde \alpha_h(x,\xi) = \mathcal{O}(h^{3/2}) \;, \quad \tilde \alpha_h(x,\xi)  - \alpha_h(x,\xi) = \mathcal{O}(h^{3/2}) \;.
\] 
Referring back to \eqref{Phfexpansion}, it then follows that: \begin{align*}
& I_1 = - h Df(x) ( M DU(x) ) + \mathcal{O}(h^2) \\
& I_2 = h \beta^{-1} \tr( D^2 f(x) M ) + \mathcal{O}(h^2) \\
& I_3 + I_4 + I_5 + I_6 = \mathcal{O}(h^2) \;.
\end{align*}
These estimates imply the desired bound \eqref{eq:onestepestimate_additive_noise}.  
\end{proof}

\medskip

Assuming constant mobility the next theorem states higher-order weak accuracy of the Metropolis integrator operated 
with $G_h(x)$ given by any second-order Runge-Kutta combination of the deterministic drift.  
To prove this statement, it is sufficient to assume the following regularity on $U(x)$.  

\medskip

\begin{assumption} \label{P2_U_Assumption_additive_noise}
The first five derivatives of $U(x)$ are bounded for all $x \in \mathbb{R}^n$.
\end{assumption}

\medskip

\begin{theorem} \label{thm:finite_time_accuracy_additive_noise}
Assume~\ref{P2_U_Assumption_additive_noise} holds and that the mobility matrix is constant.  Consider the Metropolis integrator with $G_h(x)$
given by: \begin{equation} 
G_h(x) = - b_1 M DU(x) - b_2 M DU(x_1)  \quad
x_1 = x - h a_{12} M DU(x) 
\label{rk2_Gh_additive}
\end{equation}
and parameter values satisfying the standard second-order conditions for a two-stage Runge-Kutta method: \begin{equation} 
\label{standard_rk2_conditions}
b_1 + b_2 = 1 \;, \quad b_2 a_{12} = 1/2 \;.
\end{equation}  Then for every $T>0$ there exists $C(T)>0$ such that: \[
| \E_x ( f(Y( \lfloor t/h \rfloor h))) - \E_x (f (X_{\lfloor t/h \rfloor }) ) | \le C(T) h^{3/2} 
\]
for all $t \in [0, T]$, for every $x \in \mathbb{R}^n$,  for $h$ sufficiently small, and for sufficiently regular test function $f(x)$.
\end{theorem}

\medskip

Figure~\ref{fig:E0_heavy_tailed} confirms this result.  In particular, the figure shows that the Metropolis integrator operated with $G_h(x)$ given by 
the  Ralston Runge-Kutta combination is $3/2$-weakly accurate.

\medskip

\begin{proof}
Following the proof of Theorem~\ref{thm:infinitesimalfdt}, we will prove that the scheme is one-step weakly accurate by deriving the following bound on the difference between $\mathcal{P}_h$ in~\eqref{eq:Ptf} and $P_h$ in \eqref{eq:Phf}:  \begin{equation} \label{eq:onestepestimate_additive}
| (\mathcal{P}_h f)( x ) - (P_h f)( x ) | \le C h^{5/2} 
\end{equation} for both sufficiently small $h$ and regular test functions $f(x)$.  
Using Taylor's theorem and \eqref{eq:Phf}, one can write:  \begin{equation} \label{Phfexpansion_additive}
\E_{x}  f(X_1)  = f(x) + I_1 + I_2 + I_3 + I_4  + I_5   
 \end{equation}
 where we have introduced \begin{align*}
I_1 &=  \E_{x} \left( D f(x) (X_1^{\star} - x)  \alpha_h(x, \xi) \right) \;, \\
I_2 &=  \frac{1}{2} \E_{x} \left( D^2 f(x) (X_1^{\star} - x)^2   \alpha_h(x, \xi) \right) \;, \\
I_3 &=  \frac{1}{6} \E_{x} \left( D^3 f(x) (X_1^{\star} - x)^3  \alpha_h(x, \xi)\right) \;, \\
I_4 &= \frac{1}{24} \E_{x} \left(  D^4 f(x) (X_1^{\star} - x)^4   \alpha_h(x, \xi)   \right)   \;,  \\
I_5 &= \frac{1}{24}  \E_{x}  \left( \alpha_h(x, \xi) \int_0^1 (1-s)^4  D^5 f(x + s (X_1^{\star} - x) ) (X_1^{\star} - x)^5 ds    \right) \;.
\end{align*}
Recall that $X_1^{\star}$ denotes the proposal move \eqref{proposal} with $X_0 = x$.  
(Here $D^k f(x) (v)^k$ is a $k$-linear derivative map applied to the $k$-tuple $(v, \cdots , v)$.)
Going backwards, the term $I_5$ is $\mathcal{O}(h^{3})$.  The term $I_4$ involves a quartic form in $\xi$, which contributes the following $\mathcal{O}(h^2)$ term: \begin{equation} \label{additive_I4}
I_4 = \frac{h^2}{2} \beta^{-2} \frac{\partial^4 f}{\partial x_i \partial x_j \partial x_k \partial x_{\ell}} M_{ij} M_{k \ell} + \mathcal{O}(h^3) \;.
\end{equation}
Similarly, the term $I_3$ involves a quadratic form in $\xi$, which contributes the following $\mathcal{O}(h^2)$ term:
 \begin{equation} \label{additive_I3}
I_3 = - h^2 \beta^{-1} \frac{\partial^3 f}{\partial x_i \partial x_j \partial x_k} M_{ij} M_{k \ell} \frac{\partial U}{\partial x_{\ell}} + \mathcal{O}(h^3) \;.
\end{equation}
A direct expansion shows that the term $I_2$ contributes three terms up to $\mathcal{O}(h^{5/2})$:
 \begin{equation} \label{additive_I2}
 \begin{aligned}
& I_2 = h \E_{x} \left( D^2 f(x) (B \xi)^2  \right) + \frac{h^2}{2}  D^2 f(x) (M DU, M DU) \\
& \qquad - \frac{h^2}{2}  \E_{x} \left(  D^2 f(x)( B \xi, M D^2 U(x) B \xi) \right)  + \mathcal{O}(h^{5/2}) \;, \\
& \quad =  h  \beta^{-1} \frac{\partial^2 f}{\partial x_i \partial x_j} M_{ij}  + \frac{h^2}{2}  \frac{\partial^2 f}{\partial x_i \partial x_j}  M_{ik} M_{j \ell}  \frac{\partial U}{\partial x_k} \frac{\partial U}{\partial x_{\ell}} \\
& \qquad - h^2   \beta^{-1}  \frac{\partial^2 f}{\partial x_i \partial x_j}  M_{i \ell}  M_{j k} \frac{\partial^2 U}{\partial x_k \partial x_{\ell}}   + \mathcal{O}(h^{5/2})  \;.
\end{aligned}
\end{equation}
Note that the remainder term at $\mathcal{O}(h^{5/2})$ does not vanish because the $\operatorname{min}$ function appearing in the acceptance probability $\alpha_h(x, \xi)$ breaks the odd symmetry of the integrands at that order.  Using the second-order conditions \eqref{standard_rk2_conditions} and the integration by parts argument in Theorem~\ref{thm:infinitesimalfdt}, it can be shown that the term $I_1$ contributes three terms up to $\mathcal{O}(h^{5/2})$.  Indeed, expand $I_1$ to obtain:
\begin{equation} 
I_1 =  - h D f(x) (M DU) + \frac{h^2}{2} D f(x) ( M D^2 U(x) M DU(x) ) + I_{11} + I_{12}  
\end{equation}
where we have introduced: \begin{align*}
& I_{11} = \sqrt{2 h} \E_{x} \left( D f(x) (B \xi) ( \alpha_h(x, \xi) - 1 )\right) \;, \\
& I_{12} =  - \frac{h^2}{4}  \E_{x} \left( D f(x) ( M D^3 U(x) ( B \xi, B \xi) ) \right) + \mathcal{O}(h^{5/2}) \;.
\end{align*}
The term $I_{12}$ can be directly evaluated since it is a quadratic form in $\xi$. The term $I_{11}$ 
can be estimated by introducing the following approximation (in the $L^2$ sense) to the acceptance probability function: \[
\tilde \alpha_h(x,\xi) = 1 \wedge \exp\left( - \beta   \frac{\sqrt{2} h^{3/2}}{12} D^3 U(x) ( B \xi, B \xi, B \xi) \right) \;.
\]  This approximation satisfies: \[
 \E | \tilde \alpha_h(x,\xi) - \alpha(x, \xi) | \le C h^2 \;.
\]  The proof of this statement is nearly identical to the proof of Lemma~\ref{lem:tilde_alphah} and therefore omitted.  In terms of which, write $I_{11}$ as: \begin{equation}
I_{11} = \sqrt{2 h} \E_{x} \left( D f(x) (B \xi) ( \tilde \alpha_h(x, \xi) - 1 )\right) + \mathcal{O}(h^{5/2}) \;.
\end{equation}
The integration by parts argument given in Theorem~\ref{thm:infinitesimalfdt} with $g_1(z)$ and $\vec{g}_2(z)$ in \eqref{integration_by_parts} given by \begin{align*}
& g_1(z) = 1 - \exp\left( - \beta^{-1/2}   \frac{\sqrt{2} h^{3/2}}{12} D^3 U(x) ( B z, B z, B z) \right) \\
& \vec{g}_2(z) = (2 h \beta^{-1})^{1/2} \exp\left(- \frac{|z|^2}{2}  \right) B^T D f 
\end{align*} implies that \[
I_{11} = -  \frac{h^2}{4} \beta^{-1}  (2 \pi)^{-n/2} \int_{\mathbb{R}^n}    D^3 U(x) ( M D f, B z, B z) \exp\left(- \frac{|z|^2}{2} \right)  dz +  \mathcal{O}(h^{5/2}) \;.
\]
Remarkably, substituting this estimate back into $I_1$ gives:
\begin{equation} \label{additive_I1}
\begin{aligned}
& I_1 = - h  \frac{\partial f}{\partial x_i} M_{ij} \frac{\partial U}{\partial x_j}   + \frac{h^2}{2}  \frac{\partial f}{\partial x_i}  M_{i j} M_{k \ell} \frac{\partial^2 U}{\partial x_j \partial x_k} \frac{\partial U}{\partial x_{\ell}}  \\
& \qquad - \beta^{-1} \frac{h^2}{2} \frac{\partial f}{\partial x_i}  M_{i j} M_{k \ell} \frac{\partial^3 U}{\partial x_j \partial x_k \partial x_{\ell} } + \mathcal{O}(h^{5/2}) \;.
\end{aligned}
\end{equation}

Combining estimates \eqref{additive_I1}, \eqref{additive_I2}, \eqref{additive_I3}, and \eqref{additive_I4}, and invoking Assumption~\ref{P2_Bh_assumption}, shows that:  \[
(P_h f)(x) = f(x) + h (L f)(x) +  \frac{h^2}{2} (L^2 f)(x)  + \mathcal{O}(h^{5/2}) 
\]
which agrees with an Ito-Taylor expansion of $( \mathcal{P}_h f)(x)$ to up to terms of $ \mathcal{O}(h^{5/2})$, and hence, the scheme has the desired single step accuracy.  
\end{proof}

\section{Small Noise Limit}\label{sec:deterministicaccuracy}

\subsection{General Principle}

For every $t>0$, the solution to \eqref{sde} with non-random initial condition $Y(0) = x \in \mathbb{R}^n$
converges in the small noise limit  to the solution of the following ordinary differential equation: \begin{equation} \label{ode}
\dot{\mathcal{Y}}= - M(\mathcal{Y}) DU( \mathcal{Y}) \;,~~ \mathcal{Y}(0) = x \in \mathbb{R}^n 
\end{equation}
in the $L^2$-norm, meaning that: \[
 \E_x \left\{ | Y(t) - \mathcal{Y}(t) |^2 \right\} \to 0 \;, \quad \text{as $\beta \to \infty$} 
\] 
which we denote as $Y(t) \overset{L^2}{\to} \mathcal{Y}(t) $.  At the same time, the proposal move in \eqref{proposal} satisfies,  \begin{equation} \label{deterministicupdate}
X_1^{\star} \overset{L^2}{\to} X_0 + h G_h( X_0 ) \;.
\end{equation}
Higher order accuracy of Algorithm~\ref{MetropolisIntegrator} necessitates that the deterministic update \eqref{deterministicupdate} approximates $\mathcal{Y}(h)$ to higher-order.  However, this condition is not sufficient because the actual update in Algorithm~\ref{MetropolisIntegrator} involves a Bernoulli random variable \eqref{actualupdate}, which imposes the following additional requirement for deterministic accuracy on the acceptance probability $\alpha_h(x,y)$ in \eqref{alphah}: \[
\alpha_h(X_0, \xi) \to 1\;, \quad \text{as $\beta \to \infty$} \;.
\]  Otherwise, the proposal move will not be accepted, and consequently, the integrator will be deterministically inaccurate.    From the asymptotic relation introduced in \eqref{asymptotic_alpha}, it follows that: \[
 \mathcal{E}(X_0,h) < 0 \implies \alpha_h(X_0, \xi) \to 1 \;, \quad \text{as $\beta \to \infty$} 
 \] and the Metropolis integrator acquires the deterministic accuracy of its proposal move.
 
\medskip
 
Thus, our design philosophy is to pick $G_h(x)$ and $B_h(x)$ in Algorithm~\ref{MetropolisIntegrator} so that:
 
\medskip
 
\begin{description}
\item[(1)] the update $x \mapsto x+h G_h(x)$ generates a higher order accurate approximation to the solution of \eqref{ode}; and simultaneously,  
\item[(2)] $\mathcal{E}(x,h)$ is negative definite.
\end{description}

\medskip

\noindent
We emphasize that the time step size $h$ is held fixed in the above statements.
These ideas are formulated precisely in Theorem~\ref{thm:smallnoise} below, which requires the following regularity 
assumptions on $M(x)$, $U(x)$,  $G_h(x)$, and $B_h(x)$ to hold.

\medskip

\begin{assumption} \label{P3_M_U_assumption}
The matrix-valued function $M(x)$ and its first derivative are bounded for all $x \in \mathbb{R}^n$.
The first two derivatives of the function $U(x)$ are bounded for all $x \in \mathbb{R}^n$.
\end{assumption}

\medskip

\begin{assumption} \label{P3_Gh_assumption}
The vector-valued function $G_h(x)$ and its first derivative are bounded for all $x \in \mathbb{R}^n$ and for $h$ sufficiently small.
\end{assumption}

\medskip

\begin{assumption} \label{P3_Bh_assumption}
The matrix-valued function $B_h(x)$, its inverse, and first derivative are bounded
for all $x \in \mathbb{R}^n$ and for $h$ sufficiently small. 
\end{assumption}

\medskip

\noindent
These assumptions are sufficient to prove the following theorem.

\medskip

\begin{theorem} \label{thm:smallnoise}
Consider the solutions $Y$ and $X$ produced by \eqref{sde} and Algorithm~\ref{MetropolisIntegrator}, respectively.
Let $\mathcal{Y}$ denote the exact solution to \eqref{ode} with non-random initial condition 
$\mathcal{Y}(0) = x$.  Let $\mathcal{X}$ denote the discrete path defined by: \begin{equation} \label{rkmethod}
\mathcal{X}_{k+1} = \mathcal{X}_k + h G_h(\mathcal{X}_k) \;, \qquad \mathcal{X}_0 = x \in \mathbb{R}^n \;.
\end{equation}
Assume that Assumptions~\ref{P3_M_U_assumption},~\ref{P3_Gh_assumption}, and~\ref{P3_Bh_assumption} hold, and that 
$G_h(x)$ and $B_h(x)$ in Algorithm~\ref{MetropolisIntegrator} are selected so that the following hold.

\medskip

\begin{description}
\item[(A1)] For every $T>0$ and $x \in \mathbb{R}^n$, there exists $\tilde C(T)>0$ such that: \[
|  \mathcal{Y}( \lfloor t / h \rfloor h)  - \mathcal{X}_{\lfloor t / h \rfloor} | \le \tilde C(T) h^p \;, \qquad \text{for all $t \in [0, T]$} 
\]
and for $h$ sufficiently small.

\medskip

\item[(A2)] The function $\mathcal{E}(x,h)$~in \eqref{fbeta} is negative for all $x \in \mathbb{R}^n$ such that $M(x) DU(x) \ne 0$, 
and for $h$ sufficiently small. 
\end{description}

\medskip

\noindent
 Then there exists a constant $C(T)>0$ such that \[
\lim_{\beta \to \infty} ( \E_x | Y( \lfloor t/h \rfloor h) - X_{\lfloor t/h \rfloor} |^2 )^{1/2} \le C(T) h^p  
\]
for every $t \in [0, T]$, initial condition $x \in \mathbb{R}^n$ and sufficiently small $h$.
\end{theorem}

\medskip

In practice, to meet Assumption (A1) we select $G_h(x)$ in Algorithm~\ref{MetropolisIntegrator} to be an $n$-stage Runge-Kutta combination.  
We then pick $B_h(x)$ and the parameters in $G_h(x)$ to ensure that Assumption (A2) is also satisfied.  
 Note that the points where $M(x) DU(x) = 0$ are equilibrium points of the deterministic dynamics, where the physical solution stays forever and 
 $\mathcal{E}(x,h) = 0$.   Thus, (A2) only requires $\mathcal{E}(x,h)<0$ at all $x \in \mathbb{R}^n$ that are not critical points of $U(x)$.  

\bigskip

\begin{proof}
By the triangle inequality, the desired error can be decomposed into \begin{align} \label{errordecomposition}
 & ( \E_x | Y( \lfloor t/h \rfloor h) - X_{\lfloor t/h \rfloor} |^2 )^{1/2}  \le
\overset{\text{$\le \tilde C(T) h^p$ by Assumption (A1)}}{\overbrace{| \mathcal{Y}( \lfloor t/h \rfloor h) - \mathcal{X}_{\lfloor t/h \rfloor} |}}  \\
& \qquad    \qquad
+  \underset{\text{$\overset{L^2}{\to} 0$ by Lemma~\ref{L2Y}}}{\underbrace{ ( \E_x | Y( \lfloor t/h \rfloor h) - \mathcal{Y}({\lfloor t/h \rfloor} h) |^2 )^{1/2} }} 
+  \underset{\text{\begin{tabular}{c} $\overset{L^2}{\to} 0$ by Lemma~\ref{L2X} \\ using Assumption (A2) \end{tabular}}}{\underbrace{ ( \E_x | X_{\lfloor t/h \rfloor} - \mathcal{X}_{\lfloor t/h \rfloor} |^2 )^{1/2} }} \;. \nonumber
\end{align}
By Assumption (A1), $\mathcal{X}$ is a $pth$-order approximation of $\mathcal{Y}(t)$, and therefore, the first term in the upper bound in \eqref{errordecomposition} is bounded by $C(T) h^p$.  The remaining terms vanish in the small noise limit.  Indeed, by applying the Ito-Taylor formula to $| Y(t) - \mathcal{Y}(t) |^2$, and using Gronwall's lemma, it can be shown that \begin{equation}
 \E_x | Y(t) - \mathcal{Y}(t) |^2 \le \beta^{-1}   C_2 \exp(C_1 T) \;,~~ \forall ~~ t \in [0, T] \;.
\end{equation}
From this somewhat crude estimate, it follows that the second term appearing in the upper bound in \eqref{errordecomposition} vanishes in the small noise limit.  The third term in this upper bound vanishes because the mean acceptance probability equals one in the small noise limit as a consequence of Assumption (A2).    To prove this statement we show that the mean acceptance probability is sandwiched between one and the mean of a lower bound on the acceptance probability.  We use the dominated convergence theorem to show that this lower bound converges to one in the small noise limit, and hence, the mean acceptance probability converges to one too.   For more details see Lemmas~\ref{L2Y} and \ref{L2X} in the Appendix.
\end{proof}

\bigskip

\subsection{Necessity of Second-Order Accuracy}

Here we analytically show what can go wrong with the simple choice: \[
G_h(x)=- M(x) DU(x) \;.
\]   This $G_h(x)$ corresponds to a forward Euler approximation to \eqref{ode}.  For the sake of clarity,  assume that the mobility is constant and set $B_h=B$.  Using Taylor's theorem \eqref{fbeta} in this case can be written as: \[
\mathcal{E}(x, h) = h^2  \left( \int_0^1 (1-s) D^2 U(x_s) (M DU(x), M DU(x)) ds \right) \;, ~~ x_s = x - s h M DU(x) \;.
\]
From this expression it is clear that the function $\mathcal{E}(x,h)$ can become positive in regions where $U(x)$ is convex, no matter how small the time step size $h$ is made.  Therefore, the rejection rate of the Metropolis integrator, Algorithm~\ref{MetropolisIntegrator}, will tend to one in the small noise limit.  In other words, even though the proposal move is first-order accurate, it is never accepted, and hence, the Metropolis integrator fails to be first-order accurate.  This  prediction was verified in the numerical experiments provided in \S\ref{ex:2Ddoublewell}, and motivates using proposal moves that are deterministically second or higher-order accurate.

\subsection{Second-Order Accuracy: Property (P3) of Integrator} \label{second_order}

Here we show how Assumptions (A1) and (A2) of Theorem~\ref{thm:smallnoise} are met by the Metropolis integrator, Algorithm~\ref{MetropolisIntegrator}. To satisfy Assumption (A1) of  Theorem~\ref{thm:smallnoise} with $p=2$, we define $G_h(x)$ to be a two-stage Runge-Kutta combination \eqref{rk2_Gh} with parameters $b_1$, $b_2$, $b_3$, $b_4$, and $a_{12}$.   This staggered discretization of the deterministic drift results in second-order conditions which are not the standard conditions for a two-stage Runge-Kutta method.  To derive these order conditions, Taylor expand $G_h(x)$ and the solution to \eqref{ode} to obtain: \begin{equation} \label{second_order_conditions}
b_1 + b_2 + b_3 + b_4 = 1 \;,  ~~(b_2 + b_4) a_{12} =1/2 \;,  ~~b_3 = b_2 
\end{equation}
and so, we can choose $a_{12}$ and $b_4$ as free parameters.  


These free parameters and the matrix-valued function $B_h(x)$ are chosen to satisfy Assumption (A2).  
To illustrate the key ideas in this process, assume for simplicity that the mobility matrix is constant and that $B_h = B$.  (The proof of Theorem~\ref{secondorder} below extends these ideas to the case of non-constant mobility.)  In this case there is only one free parameter $a_{12}$.   Use Taylor's theorem to write $\mathcal{E}(x,h)$ as:
\begin{equation}  \label{fbetark2} 
\begin{aligned}
 & \mathcal{E}(x, h)  =   - \frac{h^3}{4} \left(M - \frac{h}{2} M A(x) M \right) (A(x) M DU(x), A(x) M DU(x) )  \\
&  \qquad +  \frac{h^2}{2} \left( B(x) -  A(x) \right)(G_h(x), G_h(x))  
\end{aligned}
\end{equation}
where $A(x)$ and $B(x)$ are defined as the following line integrals of $D^2 U(x)$: \begin{align}
& A(x) = \int_0^1 D^2 U(x - s h a_{12} M DU(x) )  ds  \;, \label{rk2_A} \\
& B(x) = 2 \int_0^1 (1-s) D^2 U(x+s h G_h(x) ) ds \;.
\end{align}
It is important to stress that $\mathcal{E}(x,h)$, $G_h(x)$, $B(x)$ and $A(x)$ depend on $a_{12}$,
since Assumption (A2) will be met by adjusting this parameter.
If $U(x)$ is quadratic, then these line integrals become constant matrices: \[
A = B = D^2 U(x) 
\] and $\mathcal{E}(x,h)$ reduces to:
\begin{equation}  \label{fbetark2_quadratic} 
\mathcal{E}(x, h)  =   - \frac{h^3}{4} \left(M - \frac{h}{2} M A M \right) (A M DU(x), A M DU(x) )   \;.
\end{equation}
The condition, $h < 2 / \| M  A \|$, ensures that $\mathcal{E}(x,h)<0$ for all $x\in \mathbb{R}^n$ -- as required by Assumption (A2) of Theorem~\ref{thm:smallnoise}.  In the non-quadratic case, the sign of $\mathcal{E}(x,h)$ can be controlled by a dominant term argument.  This argument is based on expanding the second term in \eqref{fbetark2} to get:
\begin{equation}  \label{fbeta_ralston} 
\begin{aligned}
 & \mathcal{E}(x, h)  =   - \frac{h^3}{4} \left(M - \frac{h}{2} M A(x) M \right) (A(x) M DU(x), A(x) M DU(x) )  \\
&  \qquad + \frac{h^3}{4} \left( a_{12} - \frac{2}{3}  \right)  D^3 U(x) ( M DU(x) , M DU(x), M DU(x) )  + \mathcal{O}(h^4) \;.
\end{aligned}
\end{equation}
The quadratic form appearing in the first term of \eqref{fbeta_ralston} is negative-definite if $h< 2 / C$, where $C$ is a uniform bound on $M D^2 U(x)$.
The second term vanishes if $a_{12}=2/3$.   Since the $\mathcal{O}(h^3)$ term in \eqref{fbeta_ralston} dominates over the $\mathcal{O}(h^4)$ term, $\mathcal{E}(x,h)<0$ for every $x \in \mathbb{R}^n$ and $h$ small enough, except in the following region: 
\begin{equation} \label{bad_region}
\{ x \in \mathbb{R}^n ~~  \mid ~~ A(x)  M DU(x) = 0 \;,~~ M DU(x) \ne 0  \;,~~\text{and}~~ a_{12}=2/3 \}  \;.
\end{equation} 
Theorem~\ref{secondorder} assumes this set is empty.

While this assumption is sufficient for the dominant term argument to be valid, it does not hold in general.  For a concrete example,  take a look back at the gray-colored regions in the insets of Figures~\ref{fig:dw2db} (a) and (b).  This figure illustrates that there are regions where the line integral $A(x)$ in \eqref{rk2_A} crosses a singularity of $D^2 U(x)$, which are the x-marked points in the insets.  Assuming that these singularities are isolated, these regions are small for small $h$ but they can (and will eventually) be reached by the Metropolis integrator if the noise is nonzero.

\begin{figure}[Ht!]
\centering
\includegraphics[width=0.9\textwidth]{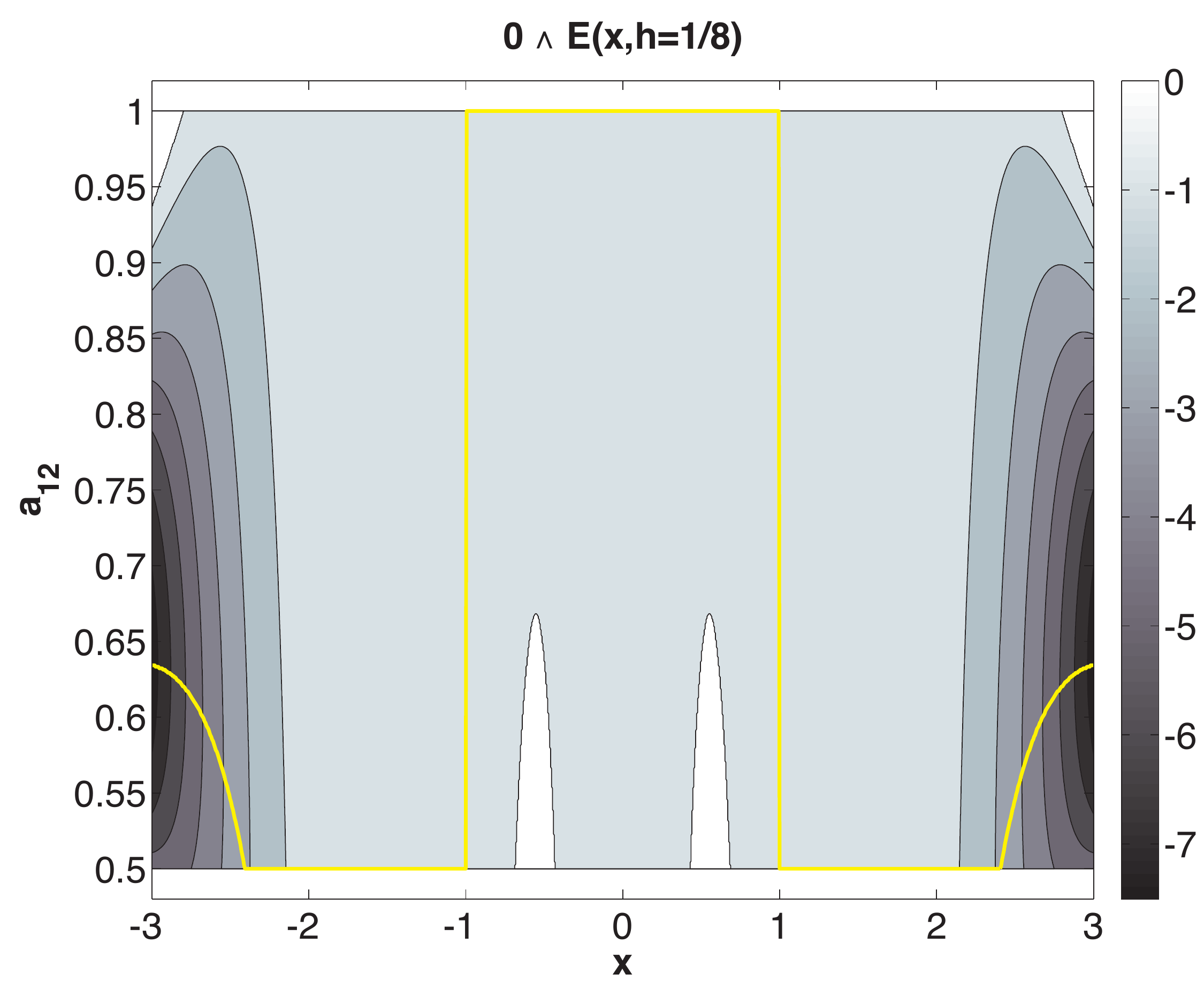} 
%
\caption{ \small {\bf Contour plot of $0 \wedge \mathcal{E}(x,h)$ for a one-dimensional double well example.}
The contour lines of $0 \wedge \mathcal{E}(x,h)$ are plotted as a function of $x$ and the parameter
$a_{12}$ for a double well potential with $U(x) = (1-x^2)^2 / 4$ and $h=1/8$.  The white-filled contours indicate regions
where $\mathcal{E}(x,h)$ is non-negative, and where the integrator operated with the corresponding $a_{12}$ can get stuck if the noise is small enough.  
The yellow curve plots the value of $a_{12}$ for which $\mathcal{E}(x,h)$ attains a minimum as a function of $x$.  
This plot confirms that both the patched and optimized proposal moves (\eqref{patched_proposal} and \eqref{optimized_proposal}, respectively) 
get around the problem associated with zeros of $U''(x)$ (located at $x=\pm 1/\sqrt{3} \approx \pm 0.58$)
by tuning the value of $a_{12}$ to make $\mathcal{E}(x,h)$ strictly negative.  
}
\label{fig:dw1d}
\end{figure}

However, if this does happen, the parameter $a_{12}$ can be tuned to preserve the dominant term argument.  Indeed, referring back to \eqref{fbeta_ralston}, Lemma~\ref{omegas_are_disjoint} implies that the first term in  \eqref{fbeta_ralston} can be made strictly negative by switching to any value of $a_{12}$ not equal to $2/3$, and the second term in \eqref{fbeta_ralston} can be made non-positive 
by choosing either $a_{12} < 2/3$ or $>2/3$.  Thus,  one way to prevent the integrator from getting stuck in these regions, is to replace the proposal move in \eqref{proposal} with:
\begin{equation}
\label{patched_proposal} 
\begin{cases}
\tilde X_1 =  X_0 + \sqrt{\frac{h}{2}} B_h(X_0) \xi \\
a_{12} = \left. \argmin_{\alpha \in \{ 2/3, 1, 1/2 \} } \mathcal{E}(\tilde X_1, h) \right|_{a_{12}=\alpha} \\
X_1^{\star} = 2 \tilde X_{1}  - X_0 + h\, G_h ( \tilde X_1)   
\end{cases}
\end{equation}
To assess the validity of this patch, Figure~\ref{fig:dw1d} plots the contours of $\mathcal{E}(x,h) $ as a function of $x$ and $a_{12}$ for the overdamped dynamics of a particle in a double-well potential: $U(x) = (1-x^2)^2 / 4$.  In this figure the white-filled contours correspond to points where $ \mathcal{E}(x,h) $ is non-negative and where the integrator operated with the corresponding $a_{12}$ can get stuck if the noise is small enough.   In this example trajectories starting near $x=0$ must pass through a zero of $U''(x)$ to end at the bottom of the well -- these zeros are located at $x=\pm 1/\sqrt{3} \approx \pm 0.58$.  In a neighborhood of these zeros, the set in \eqref{bad_region} is not empty which implies that $\mathcal{E}(x,h)$ is non-negative, as indicated by the peaks of the white-filled contours.  For this example the second term in \eqref{fbeta_ralston} takes the form: \[
 \frac{h^3}{6} \left(  \frac{2}{3} - a_{12}\right) x^4 (1-x^2)^3  \;.
\] If  $a_{12}<2/3$ then this term  is non-negative for $x \in (-1, 1)$, which is consistent with the white-filled contour for $a_{12}<2/3$.  
If $a_{12}>2/3$ then this term is non-positive for $x \in (-1,1)$, and together with the first term in \eqref{fbeta_ralston}, able to dominate the $\mathcal{O}(h^4)$ terms in \eqref{fbeta_ralston}, as illustrated by the gray-filled contours for $a_{12}>2/3$. 
(The white-filled contours in the corners of the graph are due to the finiteness of $h$ and the unboundedness of $ U''(x)$ in this example.  These regions shrink if the time step size is reduced, while the white-filled regions near the singularities persist.)



\medskip

When the mobility matrix is not constant, the choice $M_h(x) = M(x)$ does not always lead to a negative definite $\mathcal{E}(x,h)$ function.   This loss of definiteness motivates the approximation \eqref{rk2_Bh} with parameters $d_1$, $d_2$, and $c_{12}$.  First-order accuracy at constant temperature requires that  $d_1 + d_2 = 1$ which leaves $d_2$ and $c_{12}$ as free parameters.   Following the preceding dominant term argument, we will choose $a_{12}=2/3$, and the remaining free parameters ($b_4$, $d_2$, and $c_{12}$) to remove the indefinite terms appearing at $\mathcal{O}(h^2)$ and $\mathcal{O}(h^3)$ that involve the derivatives of the mobility.  By eliminating these terms, $\mathcal{E}(x,h)$ will be dominated by the leading term appearing in \eqref{fbetark2}, and hence, the argument presented when the mobility matrix is constant carries over to the non-constant case.  
The following regularity conditions on $U(x)$ and $M(x)$ are sufficient for this argument.

\medskip

\begin{assumption} \label{P3_U_Assumption}
The first four derivatives of the potential energy function $U(x)$ are bounded for all $x \in \mathbb{R}^n$.
\end{assumption}

\medskip

\begin{assumption} \label{P3_M_Assumption}
The matrix-valued function $M(x)$ and its first three derivatives are bounded for all $x \in \mathbb{R}^n$.
\end{assumption}

\medskip

With these assumptions we are now in position to state and prove Property (P3).  

\medskip

\begin{theorem} \label{secondorder}
Consider the solution $\mathcal{Y}$ to \eqref{ode} with the non-random initial condition $x \in \mathbb{R}^n$ and suppose that $U(x)$ and $M(x)$
satisfy Assumptions~\ref{P3_U_Assumption} and~\ref{P3_M_Assumption}, respectively.  In addition, suppose that 
there exists an $\epsilon >0$ such that $\| A(x) M(x) DU(x) \| \ge \epsilon \| M(x) DU(x) \|$ for all $x \in \mathbb{R}^n$ and $h$ sufficiently small. 
Let $X$ denote the approximation produced by Algorithm~\ref{MetropolisIntegrator} with $G_h(x)$ in \eqref{rk2_Gh} and $B_h(x)$ in \eqref{rk2_Bh} with 
parameter values given in \eqref{rk2_optimal_parameters}.  For every $T>0$ there exists a constant $C(T)>0$ such that \[
\lim_{\beta \to \infty} ( \E_x | Y( \lfloor t/h \rfloor h) - X_{\lfloor t/h \rfloor} |^2 )^{1/2} \le C(T) h^2  
\]
for every $t \in [0, T]$, initial condition $x \in \mathbb{R}^n$ and sufficiently small $h$.
\end{theorem}

\bigskip

In addition to the notation introduced in the proof of Theorem~\ref{thm:infinitesimalfdt}, we use the following shorthand for the second derivative of the mobility: \[
D^2 M(x) (s, u, v, w) = \frac{\partial M_{ij}}{\partial x_k \partial x_l} s_l u_k v_j w_i \;.
\]  Symmetry of second derivatives implies that $D^2 M(x) (s, u, v, w) = D^2 M(x) (u, s, v, w)$.  For brevity, if a function is evaluated at $x$, the function's input is omitted. 

\bigskip

\begin{proof}
It suffices to show that Assumption (A2) is satisfied by Algorithm~\ref{MetropolisIntegrator} operated with $G_h(x)$ in~\eqref{rk2_Gh}, $B_h(x)$ in~\eqref{rk2_Bh}, 
and the parameter values given in \eqref{rk2_optimal_parameters}, since we already established Assumption (A1).
Based on the preceding discussion, assume the second order conditions \eqref{second_order_conditions} hold and $a_{12}=2/3$.  
These conditions leave $b_4$, $d_2$, and $c_{12}$ as free parameters.  Rewrite $\mathcal{E}(x,h)$ as, \begin{equation}  \label{fbetaMh}
\begin{aligned}
 & \mathcal{E}(x,h)  =  h G_h^T M^{-1} G_h + U(x+ h G_h) - U  \\
 & \qquad  + h G_h^T M^{-1} (M - M_h(x+h G_h)) M^{-1} G_h     \\
 & \qquad  + h G_h^T M^{-1} (M - M_h(x+h G_h)) M^{-1} (M - M_h(x+h G_h)) M^{-1} G_h \\
 & \qquad + \mathcal{O}(h^4)  \;, 
\end{aligned}
\end{equation}
where we have used the identity: $M_h^{-1} = M^{-1} - M^{-1} (M_h - M) M_h^{-1}$.    
A Taylor expansion of the mobility matrix implies $G_h(x)$ can be written as, \begin{equation}   \label{ralstonG} \begin{aligned}
& G_h = - M DU + \frac{h}{2} (M A M DU + DM(M DU, DU) )  \\
& \qquad - \frac{h^2}{6} D^2 M(MDU, MDU, DU) - h^2 b_4 a_{12}^2 DM(M DU, A M DU) + \mathcal{O}(h^3) \;.  
\end{aligned}
\end{equation}
Substituting \eqref{ralstonG} into \eqref{fbetaMh} yields,
\begin{equation} \label{fbetaGMh}
 \begin{aligned}
& \mathcal{E}(x,h)  =  - \frac{h^3}{4} \left(M - \frac{h}{2} M A M \right) (A M DU, A M DU )   +  \frac{h^2}{2} \left( B -  A \right)(G_h, G_h)   \\
& \qquad - \frac{h^2}{2} DM(MDU, DU, DU) \\
& \qquad + \frac{h^3}{4} DM(MDU, DU)^T M^{-1} DM(MDU, DU)  \\
& \qquad + h^3 b_4 a_{12}^2 DM(MDU, DU, AMDU)  \\
& \qquad + \frac{h^3}{6} D^2 M(MDU, MDU, DU, DU)  \\
& \qquad +  h G_h^T M^{-1} (M - M_h(x+h G_h)) M^{-1} G_h   \\
& \qquad  + h G_h^T M^{-1} (M - M_h(x+h G_h)) M^{-1} (M - M_h(x+h G_h)) M^{-1} G_h \\
& \qquad + \mathcal{O}(h^4)     \;.
\end{aligned}
\end{equation}
Observe that the leading term in \eqref{fbetaGMh} includes \eqref{fbetark2} and additional terms involving the derivatives of the mobility.  Using \eqref{rk2_Bh} expand 
$M_h(x+h G_h)$ as: \begin{equation}  \label{Mh2}
M_h(x+ h G_h) = M + h (d_2 c_{12} - 1) DM ( MDU) + \mathcal{O}(h^2) \;.
\end{equation}
Substitute \eqref{Mh2} into \eqref{fbetaGMh} to obtain: \begin{equation} \label{fbetaGMh2} 
\begin{aligned}
& \mathcal{E}(x,h)  =  - \frac{h^3}{4} \left(M - \frac{h}{2} M A M \right) (A M DU, A M DU )   +  \frac{h^2}{2} \left( B -  A \right)(G_h, G_h)  \\
& \qquad + h^2 ( \frac{1}{2} - d_2 c_{12}  ) DM(MDU, DU, DU) + \mathcal{O}(h^3)  \;.  
\end{aligned}
\end{equation}
If $d_{2} c_{12} = 1/2$ the derivatives of the mobility up to $\mathcal{O}(h^3)$ are eliminated in this expression.   Under this condition \begin{align}
& M_h(x+ h G_h) = M - \frac{h}{2} DM ( MDU)  + c_{12}  \frac{h^2}{4} D^2M(MDU,MDU)  + \mathcal{O}(h^3) \;. \label{Mh3}
\end{align}
Substituting \eqref{Mh3} into \eqref{fbetaGMh} yields: \begin{equation} \label{fbetaGMh3} 
\begin{aligned}
& \mathcal{E}(x,h)  =  - \frac{h^3}{4} \left(M - \frac{h}{2} M A M \right) (A M DU, A M DU )   +  \frac{h^2}{2} \left( B -  A \right)(G_h, G_h)  \\
& \qquad + h^3 (b_4 a_{12}^2  - \frac{1}{2} ) DM(MDU, DU, A M DU)   \\
& \qquad + h^3 (   \frac{1}{6} - \frac{c_{12}}{4} ) D^2M(MDU, MDU, DU, DU) + \mathcal{O}(h^4)  \;.  
\end{aligned}
\end{equation}
Eliminating the terms up to $\mathcal{O}(h^3)$ requires the conditions: $b_4 a_{12}^2 = 1/2$ and $c_{12} = 2/3$.     The above conditions uniquely specify the parameters appearing in \eqref{rk2_Gh} and \eqref{rk2_Bh}.  Moreover, they imply that the leading term in \eqref{fbetaGMh3} is the same as \eqref{fbetark2}, and by the same dominant asymptotic argument, $\mathcal{E}(x,h)$ is negative definite when $h$ is small enough.   
\end{proof}

\subsection{Optimizing the Proposal Move}

Another approach to choosing the free parameter $a_{12}$ is to set it be the value at which $\mathcal{E}(x,h)$ is minimized.  The resulting algorithm is identical to the Metropolis integrator, but with the proposal move in \eqref{proposal} replaced with
\begin{equation}
\label{optimized_proposal} 
\begin{cases}
\tilde X_1 =  X_0 + \sqrt{\frac{h}{2}} B_h(X_0) \xi \\
a_{12} =  \argmin_{\alpha \in [1/2, 1]} \left. \mathcal{E}(\tilde X_1, h) \right|_{a_{12}=\alpha} \\
X_1^{\star} = 2 \tilde X_{1}  - X_0 + h\, G_h ( \tilde X_1)   
\end{cases}
\end{equation} 
The minimization problem in \eqref{optimized_proposal} adds implicitness to the algorithm, but is straightforward to solve since it requires finding the lowest point on the graph of $ \left. \mathcal{E}(x,h) \right|_{a_{12}=\alpha}$ for $\alpha \in [1/2, 1]$.   The yellow curve in Figure~\ref{fig:dw1d} shows the solution to this minimization problem as a function of $x$
for a one-dimensional double well example.  Note that this curve avoids the white-filled regions, that is, along this curve we have $\mathcal{E}(x, h) < 0$.

\subsection{Third-Order Accuracy}

Define $G_h(x)$ to be the following three-stage Runge-Kutta combination: \begin{equation}
\begin{cases}
G_h(x) = - b_1 M(x) DU(x) - b_2 M(x_1) DU(x_1) - b_3 M(x_2) DU(x_2)  \\
x_1 = x - h a_{12} M(x) DU(x)  \\
x_2 = x - h a_{31} M(x) DU(x) - h a_{32} M(x_1) DU(x_1) 
\end{cases}
\label{rk3_Gh}
\end{equation}
 This staggered discretization of the mobility matrix results in conditions for third-order accuracy which are not the standard ones for a three-stage Runge-Kutta method.  
Third-order accuracy requires that: \begin{equation} 
\label{third_order_conditions}
\begin{cases}
b_1 + b_2 + b_3 = 1   \\
b_2 a_{12} + b_3 (a_{31} + a_{32}) =1/2 \\
b_2 a_{12}^2 + b_3 (a_{31} + a_{32})^2 = 1/3    \\
 a_{12} a_{32} b_3 = 1/6 
\end{cases}
\end{equation}
and so, we can choose $a_{31}$ and $a_{32}$ as free parameters.  In a corresponding manner, the covariance matrix of the approximation to the noise is defined as: \begin{equation} \label{rk3_Bh}
\begin{cases}
B_h(x) B_h(x)^T = d_1 M(x) + d_2 M( \bar x_1) + d_3 M( \bar x_2)    \\
\bar x_1 = x + h c_{12} M(x) DU(x)  \\
\bar x_2 = x + h c_{31} M(x) DU(x) + h c_{32} M( \bar x_1) DU( \bar x_1) 
\end{cases}
\end{equation}
First-order accuracy at constant temperature requires that  $d_1 + d_2 + d_3 = 1$ which leaves five free parameters.  

The following three assumptions are sufficient to select the free parameters in  \eqref{rk3_Gh} and \eqref{rk3_Bh} so that (A2) of Theorem~\ref{thm:smallnoise} is also satisfied.  

\medskip

\begin{assumption} \label{P3_U_assumption_2}
The first five derivatives of the potential energy function $U(x)$ are bounded for all $x \in \mathbb{R}^n$.
\end{assumption}

\medskip

\begin{assumption} \label{P3_M_assumption_2}
The matrix-valued function $M(x)$ and its first four derivatives are bounded for all $x \in \mathbb{R}^n$.
\end{assumption}

\medskip

Interestingly, every third-order Runge-Kutta method yields a deterministically third-order Metropolis algorithm provided the parameters in \eqref{rk3_Bh} are coupled to the parameters in \eqref{rk3_Gh} in the manner prescribed in the following theorem.
To state this theorem, we require the following matrix: \begin{equation} \label{rk3_A}
A_{\mathtt{RK3}}(x)  = 2 (b_2 a_{12} A_1(x) + b_3 (a_{31} + a_{32}) A_2(x))  
\end{equation}
where we have introduced the line integrals: \begin{align*}
& A_1(x) = \int_0^1 D^2 U( x - s h a_{12} M(x) DU(x)) ds   \\
& A_2(x) = \int_0^1 D^2 U(x - s h (a_{31} M(x) DU(x) + a_{32} M(x_1) DU(x_1) ) ) ds  \\
& x_1 = x - h a_{12} M(x) DU(x) 
\end{align*}

\begin{theorem} \label{thirdorder}
Assume that $U(x)$ and $M(x)$ satisfy Assumptions~\ref{P3_U_assumption_2} and~\ref{P3_M_assumption_2}, respectively.
 In addition, suppose that 
there exists an $\epsilon >0$ such that $\| A_{\mathtt{RK3}}(x) M(x) DU(x) \| \ge \epsilon \| M(x) DU(x) \|$ for all $x \in \mathbb{R}^n$ and $h$ sufficiently small. 
Let $X$ denote the approximation Algorithm~\ref{MetropolisIntegrator} with $G_h(x)$ given by \eqref{rk3_Gh}, $B_h(x)$ given by \eqref{rk3_Bh}, 
and parameter values satisfying: \begin{equation} \label{optimalMhparameters}
d_1 = b_1 \;, ~~ d_2 = b_2 \;, ~~ d_3 = b_3 \;, ~~ c_{12} = a_{12} \;,~~ c_{31} = a_{31} \;, ~~ c_{32} = a_{32} \;.  
\end{equation}
For every $(b_1, b_2, b_3, a_{12}, a_{31}, a_{32})$ satisfying \eqref{third_order_conditions}, and for every $T>0$, there exists a constant $C(T)>0$ such that \[
\lim_{\beta \to \infty} ( \E_x | Y( \lfloor t/h \rfloor h) - X_{\lfloor t/h \rfloor} |^2)^{1/2} \le C(T) h^3 
\]
for every $t \in [0, T]$, initial condition $x \in \mathbb{R}^n$ and sufficiently small $h$.
\end{theorem}

\medskip

\begin{proof}
According to Theorem~\ref{thm:smallnoise}, to prove this statement it suffices to show that Assumption (A2) is met by 
Algorithm~\ref{MetropolisIntegrator} operated using $G_h(x)$ in~\eqref{rk3_Gh}, $B_h(x)$ in~\eqref{rk3_Bh}, and 
parameter values  satisfying \eqref{third_order_conditions} and \eqref{optimalMhparameters}.  A Taylor expansion of $G_h(x)$ 
yields: \begin{align}   \label{Grk3expansion}  
 G_h =  -M DU + \frac{h}{2} ( M A_{\mathtt{RK3}} MDU + DM(MDU, DU) ) + R_1 h^2 + \mathcal{O}(h^3) \;. 
\end{align}
where $A_{\mathtt{RK3}}$ is given in \eqref{rk3_A} and $R_1$ is defined as: \begin{equation} \label{R1}
\begin{aligned}
& R_1 = - \frac{h^2}{6}  M A_2 M A_1 M DU - \frac{h^2}{6} D^2 M (MDU, M DU, DU)   \\
& \qquad - \frac{h^2}{2} DM(MDU, D^2 U M DU) \\
& \qquad - \frac{h^2}{6} DM(M D^2 U M DU + DM(MDU, DU), DU) \;. 
\end{aligned}
\end{equation}
Observe that \eqref{Grk3expansion} matches \eqref{ralstonG} up to $\mathcal{O}(h^2)$.  It follows that \begin{equation}
\label{fbetark3}
\begin{aligned} 
& \mathcal{E}(x,h)  =  - \frac{h^3}{4} \left(M - \frac{h}{2} M A_{\mathtt{RK3}} M \right) (A_{\mathtt{RK3}} M DU, A_{\mathtt{RK3}} M DU )   \\
& \qquad +  \frac{h^2}{2} \left( B -  A_{\mathtt{RK3}} \right)(G_h, G_h)  - \frac{h^2}{2} DM(MDU, DU, DU) \\
& \qquad + \frac{h^3}{4} DM(MDU, DU)^T M^{-1} DM(MDU, DU)  \\
& \qquad  - h^3 R_1^T DU +  h G_h^T M^{-1} (M - M_h(x+h G_h)) M^{-1} G_h   \\
& \qquad  + h G_h^T M^{-1} (M - M_h(x+h G_h)) M^{-1} (M - M_h(x+h G_h)) M^{-1} G_h \\
& \qquad + \mathcal{O}(h^4)     \;,
\end{aligned}
\end{equation} where as before $B = 2 \int_0^1 (1-s) D^2 U(x + s h G_h) ds$.  In the case when $U$ is quadratic: \[
A_{\mathtt{RK3}} = B = A_1 = A_2 
\]
and \eqref{fbetark3} simplifies to: \[
\mathcal{E} = - h^3 \left(\frac{1}{12} M - \frac{h}{8} M A_{\mathtt{RK3}} M \right) (A_{\mathtt{RK3}} M DU, A_{\mathtt{RK3}} M DU ) \;,~~ \text{(if $U$ is quadratic)} \;.
\]
This expression is negative if $h < 2/ ( 3 \| M A_{\mathtt{RK3}} \| ) $.  In the non-quadratic case, and as before, the norm of $M A_{\mathtt{RK3}}$ is replaced with a local Lipschitz constant on $DU$.

Assuming \eqref{optimalMhparameters} an expansion of $M_h$ yields: \begin{equation} \label{Mhrk3}
\begin{aligned} 
& M_h(x+h G_h) = M - \frac{h}{2} DM(MDU) \\
& \qquad  + \frac{h^2}{6} \left( D^2M( MDU, MDU) + DM( M D^2 U M DU + DM(MDU, DU) ) \right) \;. 
\end{aligned} \end{equation} Substituting \eqref{R1} and \eqref{Mhrk3} into \eqref{fbetark3} and simplifying yields: \begin{equation}  \label{fbetark3final}
\begin{aligned}
& \mathcal{E}(x,h)  =  - \frac{h^3}{4} \left(M - \frac{h}{2} M A_{\mathtt{RK3}} M \right) (A_{\mathtt{RK3}} M DU, A_{\mathtt{RK3}} M DU )  \\
& \qquad  +  \frac{h^2}{2} \left( B -  A_{\mathtt{RK3}} \right)(G_h, G_h)  \\
& \qquad + \frac{h^3}{6} DU^T M A_2 M A_1 M DU  \;.  
\end{aligned}
\end{equation} This function is negative definite provided $h$ is small enough.
\end{proof}

\bigskip

Although all third-order Runge Kutta methods lead to a negative definite $\mathcal{E}(x,h)$ function in the deterministic limit for sufficiently small
time step (under the assumptions stated in the theorem),  the following theorem specifies an optimal choice of free parameters.

\bigskip

\begin{prop} \label{optimalthirdorder}
Assume that $U(x)$ and $M(x)$ satisfy Assumptions~\ref{P3_U_assumption_2} and~\ref{P3_M_assumption_2}, respectively.
Let $X$ denote the approximation induced by Algorithm~\ref{MetropolisIntegrator} with $G_h(x)$ in \eqref{rk3_Gh} 
with parameter values given by \begin{equation} \label{kuttaparameters}
b_1 = 1/6 \;, ~~ b_2 = 2/3 \;, ~~ b_3 = 1/6 \;, ~~ a_{12} = 1/2 \;, ~~ a_{31} = -1 \;, ~~ a_{32} = 2 \;.
\end{equation}
 In addition, suppose that there exists an $\epsilon >0$ such that $\| A_{\mathtt{RK3}}(x) M(x) DU(x) \| \ge \epsilon \| M(x) DU(x) \|$ for all $x \in \mathbb{R}^n$ and $h$ sufficiently small.  When $M$ is constant, and $h$ sufficiently small, this choice of parameters minimizes the indefinite remainder terms in \eqref{fbetark3}.
\end{prop}

\bigskip

\begin{proof}  
When $M$ is constant, \eqref{fbetark3} simplifies to:
\begin{equation} \label{fbetark3B} 
\begin{aligned}  
&  \mathcal{E}(x,h)  =  - \frac{h^3}{4} \left(M - \frac{h}{2} M A_{\mathtt{RK3}} M \right) (A_{\mathtt{RK3}} M DU, A_{\mathtt{RK3}} M DU )   \\
 & \qquad +  \frac{h^2}{2} \left( B -  A_{\mathtt{RK3}} \right)(G_h, G_h)   \;.
\end{aligned}
\end{equation}
A Taylor expansion of the indefinite remainder in \eqref{fbetark3B} yields \begin{align*}
& B-A_{\mathtt{RK3}} =  h \left(  b_2 a_{12}^2 + b_3 (a_{31} + a_{32})^2 - \frac{1}{3}  \right) D^3 U (M DU)  \\
& \qquad   + h^2 \left( \frac{1}{6} - \frac{a_{31} + a_{32}}{6} \right) D^3 U(x) ( D^2 U M DU ) \\
& \qquad   + h^2 \left( \frac{1}{12} - \frac{b_2 a_{12}^3}{3} - \frac{b_3 (a_{31} + a_{32})^3}{3} \right)  D^4 U (M DU, M DU) + \mathcal{O}(h^3) \;.
\end{align*}
The third-order conditions~\eqref{third_order_conditions} imply the $\mathcal{O}(h)$ term in $B-A_{\mathtt{RK3}}$  vanishes.   Within this two-parameter family of schemes it appears that one is optimal.  In particular, to eliminate the $\mathcal{O}(h^2)$ error in $B-A_{\mathtt{RK3}}$, 
the following additional conditions must be satisfied: \begin{equation}  \label{optimalrk3}
\begin{cases}
a_{31} + a_{32} = 1 \;, \\
b_2 a_{12}^3 + b_3 (a_{31} + a_{32})^3 = \frac{1}{4}\;.
\end{cases}
\end{equation}
The six equations in \eqref{third_order_conditions} and \eqref{optimalrk3} uniquely specify a three stage Runge-Kutta scheme.  It is given by the following choice of coefficients: \begin{equation}
b_1 = 1/6 \;, ~~ b_2 = 2/3 \;, ~~ b_3 = 1/6 \;, ~~ a_{12} = 1/2 \;, ~~ a_{31} = -1 \;, ~~ a_{32} = 2 \;.
\end{equation}
This scheme is known as Kutta's third order method \cite{lambert1973computational}.  
\end{proof}

 \section{Conclusion} \label{sec:conclusion}

This paper presented a $\nu$-symmetric integrator for self-adjoint diffusions.  This integrator is a Metropolis-Hastings algorithm with an optimized Runge-Kutta based proposal move and target density set equal to $\nu(x)$.  Since the Metropolis-Hastings ratio does not involve the normalization constant of $\nu(x)$, the algorithm is well-defined even in situations where this density is not normalizable (i.e., it is not a probability density).  In the context of non-normalizable $\nu(x)$, the paper proved that the algorithm is weakly accurate for finite noise as a direct consequence of its $\nu$-symmetry and its consistent approximation of the noise.  Through an asymptotic analysis of the integrator's rejection rate in the small noise limit, second-order deterministic accuracy was also established.   For BD simulations of self-adjoint diffusions, the scheme shares the nice properties of the Fixman scheme (explicitness, finite-time accuracy, second-order deterministic accuracy and avoids divergence of the mobility).  In addition, it is ergodic if $\nu(x)$ is normalizable.  These properties imply that the scheme is able to stably calculate dynamic quantities at reasonable time step sizes and  generate long trajectories.  The paper verified these claims on a collection of low-dimensional toy problems and a more realistic simulation of DNA in an unbounded solvent.  These features make the scheme appealing for the simulations of realistic BD applications, which will be the topic of future investigations.  They also motivate generalizing the algorithm to situations where the $\nu$-symmetry condition does not hold, e.g., to BD of systems that are driven out of equilibrium in ways that are more general than those considered in the present paper.

\medskip

\section{Appendix}

Here we provide Lemmas required in the proofs of Theorems~\ref{thm:infinitesimalfdt} and~\ref{thm:smallnoise}.

\medskip

\begin{lemma} \label{lem:tilde_alphah}
Consider $\alpha_h( X_0,\xi)$ and $\tilde \alpha_h( X_0, \xi)$ given in \eqref{alphah} and \eqref{tilde_alphah}, respectively.
Let Assumptions~\ref{P2_U_Assumption}, \ref{P2_Gh_assumption}, and \ref{P2_Bh_assumption} hold.
Then, \[
\E | \alpha_h( X_0,\xi)  - \tilde \alpha_h( X_0, \xi) | \le C h \;,
\]
for all $X_0 \in \mathbb{R}^n$.
\end{lemma}

\medskip

\begin{proof}
A first-order Taylor expansion of the exponent in \eqref{alphah} shows that: \begin{equation} \label{alphah_expansion}
\alpha_h( X_0,\xi) = 1 \wedge \exp\left( - \beta \left[  \Gamma_h( X_0, \xi) + R_h( X_0, \xi) \right] \right) \;,
\end{equation}
where $\Gamma_h( X_0, \xi)$ is the function introduced in \eqref{Gammah}, and $R_h( X_0, \xi)$ is a remainder term given by: \begin{align*}
& R_h( X_0,\xi) = h G_h( \tilde X_1)^T M_h( X_1^{\star})^{-1} G_h( \tilde X_1) + h DU( X_0)^T G_h( \tilde X_1) \\
& \quad + \int_0^1 (1-s) D^2U(X(s)) (X_1^{\star} - X_0, X_1^{\star} - X_0) ds \\
& \quad + \sqrt{2 h} (G_h(\tilde X_1) - G_h(X_0))^T M_h(X_1^{\star})^{-1} B_h( X_0) \xi \\
& \quad + \sqrt{2 h} G_h(X_0)^T M_h(X_1^{\star})^{-1} ( M_h( X_0) - M_h(X_1^{\star}) ) B_h( X_0)^{-T} \xi \\
& \quad + \frac{h \beta^{-1}}{2} \tr( M_h(X_0)^{-1} DM_h(X_0) ( G_h(\tilde X_1) ) ) \\
& \quad + \frac{\beta^{-1} }{2} \int_0^1 \tr( ( M_h(X(s))^{-1} D M_h(X(s)) - M_h(X_0)^{-1} D M_h(X_0) ) (X_1^{\star} - X_0) ) ds \\
& \quad - \frac{h}{2} DM_h( X_0) (G_h( \tilde X_1), B_h(X_0)^{-T} \xi, B_h(X_0)^{-T} \xi) \\
& \quad - \frac{1}{2} \int_0^1 (1-s) D^2M_h(X(s)) ) (X_1^{\star} - X_0, X_1^{\star} - X_0, B_h( X_0)^{-T} \xi, B_h( X_0)^{-T} \xi) ds \\
& \quad + \frac{1}{2} \xi^T B_h(X_0)^T M_h(X_1^{\star})^{-1} (M_h(X_0) - M_h(X_1^{\star}) )  \\
& \qquad \times M_h( X_0)^{-1} ( M_h( X_0) - M_h(X_1^{\star}) ) B_h( X_0)^{-T} \xi \;,
\end{align*}
where $X(s) = X_0 + s (X_1^{\star} - X_0)$.  This remainder term involves $B_h$, $B_h^{-T}$, $D B_h$, $D^2 B_h$, $DU$, $D^2U$, $G_h$, and $DG_h$, which are bounded by hypothesis.    The Lipschitz property of the map  $x \mapsto 1 \wedge \exp(x)$ implies that: \begin{equation}
| \alpha_h( X_0,\xi)  - \tilde \alpha_h( X_0,\xi) |^2 \le \beta^2 | R_h( X_0, \xi) |^2 \;.
\end{equation}
Taking the expectation of both sides of this inequality and invoking Assumptions~\ref{P2_U_Assumption}, \ref{P2_Gh_assumption}, and \ref{P2_Bh_assumption} 
gives the desired result.
\end{proof}

\bigskip

\begin{lemma} \label{L2Y}
Consider the solutions $Y$ to \eqref{sde} and $ \mathcal{Y}$ to \eqref{ode} with initial condition: $Y(0) =  \mathcal{Y}(0) = x \in \mathbb{R}^n$.  
For every $T>0$, then \[
 Y(t) \overset{L^2}{\to} \mathcal{Y}(t)   \;, \quad  \text{for all $t \in [0, T]$ and $x \in \mathbb{R}^n$} \;.
\]
\end{lemma}

\medskip

\begin{proof}
By the Ito-Taylor formula, \begin{align*}
& | Y(t) -  \mathcal{Y}(t) |^2 \le 2 \int_0^t \langle Y(s) -  \mathcal{Y}(s) , -M(Y(s)) DU(Y(s)) + M( \mathcal{Y}(s)) DU( \mathcal{Y}(s)) \rangle ds \\
& \qquad + 2 \beta^{-1} \int_0^t \langle Y(s) - \mathcal{Y}(s), \divergence M(Y(s)) \rangle ds + 2 \beta^{-1} \int_0^t \tr M(Y(s)) ds \\
& \qquad + 2 \sqrt{2 \beta^{-1}} \int_0^t \langle Y(s) - \mathcal{Y}(s) , B(Y(s)) dW(s) \rangle 
\end{align*} Using bounds on $M(x)$ and $U(x)$ and their derivatives, and 
taking the expectation of both sides of this inequality, it follows from Gronwall's Lemma that
\begin{align*}
 \E_x | Y(t) - \mathcal{Y}(t) |^2 &\le C_1 \int_0^t  \E_x  |Y(s) - \mathcal{Y}(s)|^2 ds + \beta^{-1} \; C_2 \; T \;, \\
& \le \beta^{-1} \; \exp(C_1 T) \; C_2 \; T \;.
\end{align*} Passing to the small noise limit  produces the desired $L^2$ convergence.   \end{proof}

\medskip

The following lemma is a discrete analog of Lemma~\ref{L2Y}.

\medskip

\begin{lemma} \label{L2X}
Let $X$ and $\mathcal{X}$ denote the numerical solutions generated by Algorithm~\ref{MetropolisIntegrator} and 
 \eqref{rkmethod} respectively, with: $X_0 = \mathcal{X}_0 = x \in \mathbb{R}^n$.  If the $\mathcal{E}(x,h)$
function of Algorithm~\ref{MetropolisIntegrator} satisfies Assumption~(A2) in Theorem~\ref{thm:smallnoise},
then for every $T>0$, \[
X_{\lfloor t/h \rfloor} \overset{L^2}{\to} \mathcal{X}_{\lfloor t/h \rfloor}   \;, \quad \text{for all $t \in [0, T]$ and $x \in \mathbb{R}^n$} \;.
\]
\end{lemma}

\medskip

\begin{proof}
The proof goes by induction over the number of steps, so it suffices to consider the difference between a single step of both schemes conditioned
on the previous time step:  \begin{align*}
 \E |X_1 - \mathcal{X}_1|^2 &= \E |X_1^{\star} - \mathcal{X}_1|^2  \;  \alpha_h(X_0, \xi) + | X_0- \mathcal{X}_1|^2 \; \E  \; (1 - \alpha_h(X_0, \xi) ) \;, \\
&  \le  (1+C_1 \; h) \;  |X_0 - \mathcal{X}_0 |^2 + C_2 \; h^{3/2}  \; \beta^{-1} + C_3 \; h \;  \E \;  (1 - \alpha_h(X_0, \xi) )  \;.
\end{align*}
The desired $L^2$ convergence follows from applying Discrete Gronwall's Lemma to this recurrence inequality, passing to the small noise limit, 
and invoking Lemma~\ref{L2alphah}.  
\end{proof}

\bigskip

\begin{lemma} \label{L2alphah}
Let $\alpha_h(X_0, \xi)$ denote the acceptance probability in \eqref{alphah}.  If the $\mathcal{E}(x,h)$
function of Algorithm~\ref{MetropolisIntegrator} satisfies Assumption~(A2) of Theorem~\ref{thm:smallnoise}, then \[
\E \; \alpha_h(X_0, \xi) \to 1   \;, \qquad \text{as $\beta \to \infty$} \;,
\]
 for all $X_0 \in \mathbb{R}^n$.
\end{lemma}

\medskip

\begin{proof}
Let $\mathcal{X}_1 = X_0 + h G_h(X_0)$ and $M_h(x) = B_h(x) B_h(x)^T$.
Write  $\alpha_h(X_0, \xi)$ as: \begin{align}
 \alpha_h(X_0, \xi) =  1 \wedge \exp\left( - \beta \left[ \mathcal{E}(X_0,h)  +  R_h(X_0,\xi) \right]  \right)
\end{align}
where $\mathcal{E}(X_0,h)$ is the function introduced in \eqref{fbeta}, and $R_h(X_0,\xi)$ is a remainder term given by: \begin{align*}
& R_h(X_0,\xi) = U(X_1^{\star}) - U(\mathcal{X}_1)  + \frac{\beta^{-1}}{2} \log \left( \frac{\det M_h(X_1^{\star})}{\det M_h(X_0) } \right)  \\
& \quad + h \left( G_h( \tilde X_1)^T M_h( X_1^{\star} )^{-1}  G_h( \tilde X_1) - G_h(X_0)^T M_h( \mathcal{X}_1 )^{-1} G_h( X_0) \right) \\
&  \quad + \sqrt{2 h} (B_h(X_0) \xi)^T M_h( X_1^{\star} )^{-1} G_h( \tilde X_1 )  \\
& \quad + \frac{1}{2} \xi^T B_h(X_0)^T M_h(X_1^{\star})^{-1} (M_h( X_0 ) - M_h(X_1^{\star}) ) B_h(X_0)^{-T} \xi  \;.
\end{align*}
Recall that, $\xi\in\mathbb{R}^n$ denotes a Gaussian random vector with 
mean zero and covariance $\E( \xi_i \xi_j ) = \beta^{-1} \delta_{ij}$.   Since: \begin{equation}
\begin{cases}
X_1^{\star} - \mathcal{X}_1 = 2 ( \tilde X_1 - X_0 ) + h ( G_h( \tilde X_1) - G( X_0 ) ) \;,  \\
\tilde X_1 - X_0 = \sqrt{\frac{h}{2}} B_h(X_0) \xi \;,
\end{cases}
\end{equation}
it follows from bounds on $G_h(x)$, $D G_h(x)$, $B_h(x)$, $B_h(x)^{-T}$, $D B_h(x)$, and $DU(x)$ that, \[
R_h( X_0, \xi) \le C_1 \; | \xi |^2 \; h^{1/2} \;.
\]
Thus, the acceptance probability satisfies  \begin{equation} \label{squeezelemma} 
\tilde \alpha_h(X_0, \xi) \le \alpha_h(X_0, \xi) \le 1 \;,
\end{equation}
where we have introduced: \[
\tilde \alpha_h(X_0, \xi) = 
 1 \wedge \exp\left( - \beta \left[ \mathcal{E}(X_0,h) + C_1 \; | \xi |^2 \; h^{1/2} ) \right] \right) \;.
\]
Let \[
\Omega_{\beta} = \left\{ z \in \mathbb{R}^{n} ~:~ C_1 \; h^{1/2} \; | z |^2 \le -  \beta \mathcal{E}(X_0,h)  \right\} \;,
\]
and in terms of which, the expectation of the lower bound in \eqref{squeezelemma} can be written as: \begin{align*}
& \E \; \tilde \alpha_h(X_0, \xi)  = \int_{\Omega_{\beta}} \exp\left( - \frac{|z|^2}{2} \right) (2 \pi)^{-n/2} d z \\
& \qquad + \exp\left( - \beta \mathcal{E}(X_0,h) \right)  \int_{\Omega_{\beta}^c} \exp\left( - \frac{|z|^2}{2} (1 + C_1 h^{1/2}) \right) (2 \pi)^{-n/2} d z \;, \\ 
& \quad \ge \int_{\Omega_{\beta}} \exp\left( - \frac{|z|^2}{2} \right) (2 \pi)^{-n/2} d z +  \int_{\Omega_{\beta}^c} \exp\left( - \frac{|z|^2}{2} (1 + C_1 h^{1/2}) \right) (2 \pi)^{-n/2} d z  \;,
\end{align*} 
where we used Assumption (A2) of Theorem~\ref{thm:smallnoise} in the last step.
Passing to the small noise limit $\beta \to \infty$, the dominated convergence theorem implies that $\E  \tilde \alpha_h(X_0, \xi)  \to 1$  
and the desired limit statement follows from applying the Squeeze Lemma to \eqref{squeezelemma}.
\end{proof}

\bigskip

In loose terms, the following Lemma says that a problem point for the Ralston Runge-Kutta combination (i.e., where $A(x) M(x) DU(x) = 0$ and $M(x) DU(x) \ne 0$) will not be a problem point for any other Runge-Kutta combination.

\medskip

\begin{lemma} \label{omegas_are_disjoint}
Assume that all of the singularities of $D^2U(x)$ are uniformly isolated, i.e., there exists an $R>0$ such that for every singularity $x$ the region: \[
\{ y \in \mathbb{R}^n ~~ \mid ~~  \| y - x \| < R \}
\] contains no other singularity.   Also assume that $DU(x)$ is bounded for all $x \in \mathbb{R}^n$.   Then for every $x \in \mathbb{R}^n$
there is at most one $\alpha \in [1/2, 1]$ such that \[
\left. A(x) \right|_{a_{12} = \alpha} ( M(x) DU(x) , M(x) DU(x) ) = 0
\] for $h$ sufficiently small.  
\end{lemma}

\medskip

\begin{proof}
Let $R$ be the minimum distance between singularities.  Let $h_c$ be small enough such that $h \| M(x) DU(x) \| < R/2$ for all $h < h_c$ and for all $x \in \mathbb{R}^n$.  
Suppose -- for the sake of contradiction -- that there exists an $x \in \mathbb{R}^n$ and $\alpha_1 \ne \alpha_2$ such that: \begin{align*}
& \left. A(x) \right|_{a_{12} = \alpha_1} ( M(x) DU(x) , M(x) DU(x) ) = 0 \;, \quad \text{and}  \\
& \left. A(x) \right|_{a_{12} = \alpha_2} ( M(x) DU(x) , M(x) DU(x) ) = 0 \;.
\end{align*}  Or equivalently: \[
\int_0^{\alpha_1} g(s) ds = 0 \;, \quad 
\int_0^{\alpha_2} g(s) ds = 0 \;,
\]  where we have introduced \[
g(s) = D^2 U(x - s h M(x) DU(x)) (M DU(x), M DU(x)) \;.
\]  Assume that $\alpha_1 < \alpha_2$.  These zero integrals imply that: \[
\int_0^{\alpha_1} g(s) ds = 0 \;, \quad 
\int_{\alpha_1}^{\alpha_2} g(s) ds = 0 \;.
\] By the mean value theorem for integrals, there exists an $s_1 \ne s_2$ such that: \[
g(s_1) = 0 \;, \quad g(s_2) = 0 \;.
\] In other words, $D^2U(x)$ possesses two singular points in the ball of radius $R/2$ centered at $x$, which is a contradiction.
\end{proof}

\newpage

\bibliographystyle{amsplain}
\bibliography{nawaf}

\providecommand{\bysame}{\leavevmode\hbox to3em{\hrulefill}\thinspace}
\providecommand{\MR}{\relax\ifhmode\unskip\space\fi MR }
\providecommand{\MRhref}[2]{%
  \href{http://www.ams.org/mathscinet-getitem?mr=#1}{#2}
}
\providecommand{\href}[2]{#2}
\begin{thebibliography}{10}

\bibitem{AkBoRe2009}
E.~Akhmatskaya, N.~Bou-Rabee, and S.~Reich, \emph{A comparison of generalized
  hybrid {M}onte {C}arlo methods with and without momentum flip}, J Comput Phys
  \textbf{228} (2009), 2256--2265.

\bibitem{AkRe2008}
E.~Akhmatskaya and S.~Reich, \emph{{GSHMC}: An efficient method for molecular
  simulation}, J Comput Phys \textbf{227} (2008), 4937--4954.

\bibitem{BiArHaCu1987}
R.~B. Bird, R.~C. Armstrong, O.~Hassager, and C.~F. Curtiss, \emph{Dynamics of
  polymeric liquids}, Wiley, 1987.

\bibitem{BoHa2013}
N.~Bou-Rabee and M.~Hairer, \emph{Non-asymptotic mixing of the {MALA}
  algorithm}, IMA J of Numer Anal \textbf{33} (2013), 80--110.

\bibitem{BoVa2010}
N.~B{ou-Rabee} and E.~V{anden-Eijnden}, \emph{Pathwise accuracy and ergodicity
  of {M}etropolized integrators for {SDE}s}, Comm Pure and Appl Math
  \textbf{63} (2010), 655--696.

\bibitem{CaLeSt2007}
E.~Canc\'{e}s, F.~Legoll, and G.~Stoltz, \emph{Theoretical and numerical
  comparison of some sampling methods for molecular dynamics}, Mathematical
  Modelling and Numerical Analysis \textbf{41} (2007), 351--389.

\bibitem{ChLa2002}
M.~Chopra and R.~G. Larson, \emph{{B}rownian dynamics simulations of isolated
  polymer molecules in shear flow near adsorbing and nonadsorbing surfaces}, J
  Rheol \textbf{46} (2002), 831.

\bibitem{de2011}
J.~J. de~Pablo, \emph{Coarse-grained simulations of macromolecules: From {DNA}
  to nanocomposites}, Annu Rev Phys Chem \textbf{62} (2011), 555--574.

\bibitem{DuKePeRo1987}
S.~Duane, A.~D. Kennedy, B.~J. Pendleton, and D.~Roweth, \emph{Hybrid
  {M}onte-{C}arlo}, Phys Lett B \textbf{195} (1987), 216--222.

\bibitem{EVa2004}
W.~E and E.~Vanden-Eijnden, \emph{Metastability, conformation dynamics, and
  transition pathways in complex systems}, Multiscale Modelling and Simulation
  (S.~Attinger and P.~Koumoutsakos, eds.), Lecture Notes in Computational
  Science and Engineering, Springer, 2004, pp.~35--68.

\bibitem{EVa2010}
\bysame, \emph{Transition-path theory and path-finding algorithms for the study
  of rare events}, Annual Review of Physical Chemistry \textbf{61} (2010),
  391--420.

\bibitem{ErMc1978}
D.~L. Ermak and J.~A. McCammon, \emph{{B}rownian dynamics with hydrodynamics
  interactions}, J Chem Phys \textbf{69} (1978), 1352--1360.

\bibitem{Fi1978}
M.~Fixman, \emph{Simulation of polymer dynamics}, J Chem Phys \textbf{69}
  (1978), 1527--1545.

\bibitem{Fi1986}
\bysame, \emph{Implicit algorithm for {B}rownian dynamics of polymers},
  Macromolecules \textbf{19} (1986), 1195.

\bibitem{GiCa2011}
M.~Girolami and B.~Calderhead, \emph{Riemann manifold {L}angevin and
  {H}amiltonian {M}onte {C}arlo methods}, J R Statist Soc B \textbf{73} (2011),
  123--214.

\bibitem{gobet2004exact}
E.~Gobet and S.~Menozzi, \emph{Exact approximation rate of killed hypoelliptic
  diffusions using the discrete euler scheme}, Stoch Proc Appl \textbf{112}
  (2004), no.~2, 201--223.

\bibitem{gobet2007discrete}
\bysame, \emph{Discrete sampling of functionals of it{\^o} processes},
  S{\'e}minaire de probabilit{\'e}s XL, Springer, 2007, pp.~355--374.

\bibitem{gobet2010stopped}
\bysame, \emph{Stopped diffusion processes: boundary corrections and
  overshoot}, Stoch Proc Appl \textbf{120} (2010), no.~2, 130--162.

\bibitem{Gr2011}
M.~D. Graham, \emph{Fluid dynamics of dissolved polymer molecules in confined
  geometries}, Annu Rev Fluid Mech \textbf{43} (2011), 273--298.

\bibitem{HaLuWa2010}
E.~Hairer, C.~Lubich, and G.~Wanner, \emph{Geometric numerical integration},
  Springer, 2010.

\bibitem{hairer2009hot}
M.~Hairer, \emph{How hot can a heat bath get?}, Comm Math Phys \textbf{292}
  (2009), no.~1, 131--177.

\bibitem{Ha1970}
W.~K. Hastings, \emph{{M}onte-{C}arlo methods using {M}arkov chains and their
  applications}, Biometrika \textbf{57} (1970), 97--109.

\bibitem{HaPa1986}
U.~G. Haussman and E.~Pardoux, \emph{Time reversal for diffusions}, Annals of
  Probability \textbf{14} (1986), no.~4, 1188--1205.

\bibitem{HedeGr2007}
J.~P. Hern\'{a}ndez-Ortiz, J.~J. de~Pablo, and M.~D. Graham, \emph{Fast
  computation of many-particle hydrodynamic and electrostatic interactions in a
  confined geometry}, Phys Rev Lett \textbf{98} (2007), 140602.

\bibitem{HeMadeGr2006b}
J.~P. Hern\'{a}ndez-Ortiz, H.~Ma, J.~J. de~Pablo, and M.~D. Graham,
  \emph{Cross-stream-line migration in confined flowing polymer solutions:
  theory and simulation}, Phys Fluids \textbf{18} (2006), 123101.

\bibitem{HeMadeGr2006a}
\bysame, \emph{${N} \log{N}$ method for hydrodynamic interactions of confined
  polymer systems: {B}rownian dynamics}, J Chem Phys \textbf{125} (2006),
  164906.

\bibitem{HeMadeGr2008}
\bysame, \emph{Concentration distributions during flow of confined flowing
  polymer solutions at finite concentration: slit and grooved channel},
  Korea-Aust Rhelo J \textbf{20} (2008), 143.

\bibitem{HeOt1997}
M.~Herrchen and H.~C. \"{O}ttinger, \emph{A detail comparison of various {FENE}
  dumbbell models}, J Non-Newton Fluid \textbf{68} (1997), 17.

\bibitem{Hi2011}
D.~J. Higham, \emph{Stochastic ordinary differential equations in applied and
  computational mathematics}, IMA J Appl Math \textbf{76} (2011), 449--474.

\bibitem{HiMaSt2002}
D.~J. Higham, X.~Mao, and A.~M. Stuart, \emph{Strong convergence of
  {E}uler-type methods for nonlinear stochastic differential equations}, IMA J
  Num Anal \textbf{40} (2002), 1041--1063.

\bibitem{Ho1991}
A.~M. Horowitz, \emph{A generalized guided {M}onte-{C}arlo algorithm}, Phys
  Lett B \textbf{268} (1991), 247--252.

\bibitem{HsLiLa2003}
C.-C. Hsieh, L.~Li., and R.~G. Larson, \emph{Modeling hydrodynamic interaction
  in {B}rownian dynamics: simulations of extensional flows of dilute solutions
  of {DNA} and polystyrene}, J Non-Newton Fluid \textbf{113} (2003), 147.

\bibitem{hsieh2011simulation}
C.-C. Hsieh and T.-H. Lin, \emph{Simulation of conformational preconditioning
  strategies for electrophoretic stretching of {DNA} in a microcontraction},
  Biomicrofluidics \textbf{5} (2011), no.~4, 044106.

\bibitem{hsieh2012simulation}
C.-C. Hsieh, T.-H. Lin, and C.-D. Huang, \emph{Simulation guided design of a
  microfluidic device for electrophoretic stretching of {DNA}},
  Biomicrofluidics \textbf{6} (2012), 044105.

\bibitem{HuShBaCh2002}
J.~S. Hur, E.~S.~G. Shaqfeh, H.~P. Babcock, and S.~Chu, \emph{Dynamics and
  configurational fluctuations of single {DNA} molecules in linear mixed
  flows}, Phys Rev E \textbf{66} (2002), 011915.

\bibitem{HuShLa2000}
J.~S. Hur, E.~S.~G. Shaqfeh, and R.~G. Larson, \emph{{B}rownian dynamics
  simulations of single {DNA} molecules in shear flow}, J Rheol \textbf{44}
  (2000), 713.

\bibitem{HuOt1998}
M.~H\"{u}tter and H.~C. \"{O}ttinger, \emph{Fluctuation-dissipation theorem,
  kinetic stochastic integral, and efficient simulations}, J Chem Soc, Faraday
  Trans \textbf{94} (1998), 1403--1405.

\bibitem{HuJeKl2012}
M.~Hutzenthaler, A.~Jentzen, and P.~E. Kloeden, \emph{Strong convergence of an
  explicit numerical method for sdes with non-globally {L}ipschitz continuous
  coefficients}, Ann Appl Probab \textbf{22} (2012), 1611--1641.

\bibitem{JedeGr2000}
R.~M. Jendrejack, J.~J. de~Pablo, and M.~D. Graham, \emph{Hydrodynamic
  interactions in long chain polymers: application of the chebyshev polynomial
  approximation in stochastic simulations}, J Chem Phys \textbf{113} (2000),
  7752.

\bibitem{JedeGr2002}
\bysame, \emph{Stochastic simulations of {DNA} in flow: Dynamics and the
  effects of hydrodynamic interaction}, J Chem Phys \textbf{116} (2002), 7752.

\bibitem{JeDiScGrde2003a}
R.~M. Jendrejack, E.~T. Dimalanta, D.~C. Schwartz, M.~D. Graham, and J.~J.
  de~Pablo, \emph{{DNA} dynamics in a microchannel}, Phys Rev Lett \textbf{91}
  (2003), 038102.

\bibitem{JeDiScGrde2003b}
\bysame, \emph{Effect of confinement on {DNA} dynamics in microfluidic
  devices}, J Chem Phys \textbf{119} (2003), 2894.

\bibitem{JeScdeGr2004}
R.~M. Jendrejack, D.~C. Schwartz, J.~J. de~Pablo, and M.~D. Graham,
  \emph{Shear-induced migration in flowing polymer solutions: simulation of
  long-chain {D}eoxyribose {N}ucleic {A}cid in micro channels}, J Chem Phys
  \textbf{120} (2004), 2513.

\bibitem{KePe2001}
A.~D. Kennedy and B.~Pendleton, \emph{Cost of the generalized hybrid {M}onte
  {C}arlo algorithm for free field theory}, Nucl. Phys. B \textbf{607} (2001),
  456--510.

\bibitem{Ke1978}
J.~Kent, \emph{Time-reversible diffusions}, Adv. Appl. Prob. \textbf{10}
  (1978), 819--835.

\bibitem{KiYoMaWa1991}
K.~Kikuchi, M.~Yoshida, T.~Maekawa, and H.~Watanabe, \emph{Metropolis {M}onte
  {C}arlo method as a numerical technique to solve the {F}okker-{P}lanck
  equation}, Chem Phys Lett \textbf{185} (1991), 335--338.

\bibitem{KiYoMaWa1992}
\bysame, \emph{Metropolis {M}onte {C}arlo method for {B}rownian dynamics
  simulation generalized to include hydrodynamic interactions}, Chem Phys Lett
  \textbf{196} (1992), 57--61.

\bibitem{kim2006brownian}
J.~M. Kim and P.~S. Doyle, \emph{A {B}rownian dynamics-finite element method
  for simulating {DNA} electrophoresis in nonhomogeneous electric fields}, J
  Chem Phys \textbf{125} (2006), 074906.

\bibitem{kim2007design}
\bysame, \emph{Design and numerical simulation of a {DNA} electrophoretic
  stretching device}, Lab on a Chip \textbf{7} (2007), no.~2, 213--225.

\bibitem{lambert1973computational}
J.~D. Lambert, \emph{Computational methods in ordinary differential equations},
  Wiley New York, 1973.

\bibitem{La2005}
R.~G. Larson, \emph{The rheology of dilute solutions of flexible polymers:
  Progress and problems}, J Rheol \textbf{49} (2005), 1.

\bibitem{LaHuSmCh1999}
R.~G. Larson, H.~Hu, D.~E. Smith, and S.~Chu, \emph{{B}rownian dynamics
  simulations of a {DNA} molecule in an extensional flow field}, J Rheol
  \textbf{43} (1999), 267.

\bibitem{LeLe2010}
T.~Lelievre and F.~Legoll, \emph{Effective dynamics using conditional
  expectations}, Nonlinearity \textbf{23} (2010), 2131--2163.

\bibitem{LeRoSt2010A}
T.~Leli\`{e}vre, M.~Rousset, and G.~Stoltz, \emph{Free energy computations: A
  mathematical perspective}, 1st ed., Imperial College Press, 2010.

\bibitem{Li2008}
J.~S. Liu, \emph{{M}onte {C}arlo strategies in scientific computing}, 2nd ed.,
  Springer, 2008.

\bibitem{CaGa1991}
J.~J. L\'{o}pez~Cascales and J.~Garc'a de~la Torre, \emph{Shear-rate dependence
  of the intrinsic viscosity of bead-and-spring chains: hydrodynamic
  interaction and excluded-volume effects}, Polymer \textbf{32} (1991), 3359.

\bibitem{Ma1999}
R.~Manella, \emph{Absorbing boundaries and optimal stopping in a stochastic
  differential equation}, Physics Letters A \textbf{254} (1999), 257--262.

\bibitem{MaFiVaCi2006}
L.~Maragliano, A.~Fischer, E.~Vanden-Eijnden, and G.~Ciccotti, \emph{String
  method in collective variables: Minimum free energy paths and isocommittor
  surfaces}, J Chem Phys \textbf{125} (2006), 024106.

\bibitem{MaVa2006}
L.~Maragliano and E.~Vanden-Eijnden, \emph{A temperature accelerated method for
  sampling free energy and determining reaction pathways in rare events
  simulations}, Chem Phys Lett \textbf{426} (2006), 168Ð175.

\bibitem{mattingly2010convergence}
J.~C. Mattingly, A.~M. Stuart, and M.~V. Tretyakov, \emph{Convergence of
  numerical time-averaging and stationary measures via {P}oisson equations},
  SIAM J Num Anal \textbf{48} (2010), no.~2, 552--577.

\bibitem{MenTw1996}
K.~L. Mengersen and R.~L. Tweedie, \emph{Rates of convergence of the {H}astings
  and {M}etropolis algorithms}, Ann Stat \textbf{24} (1996), 101--121.

\bibitem{MeRoRoTeTe1953}
N.~Metropolis, A.~W. Rosenbluth, M.~N. Rosenbluth, A.~H. Teller, and E.~Teller,
  \emph{Equations of state calculations by fast computing machines}, J Chem
  Phys \textbf{21} (1953), 1087--1092.

\bibitem{MiTr2004}
G.~N. Milstein and M.~V. Tretyakov, \emph{Stochastic numerics for mathematical
  physics}, Springer, Berlin, 2004.

\bibitem{MiTr2005}
\bysame, \emph{Numerical integration of stochastic differential equations with
  nonglobally {L}ipschitz coefficients}, IMA J Num Anal \textbf{43} (2005),
  1139--1154.

\bibitem{Nu1984}
E.~Nummelin, \emph{General irreducible {M}arkov chains and non-negative
  operators}, Cambridge University Press, New York, NY, 1984.

\bibitem{Ot1996}
H.~C. \"{O}ttinger, \emph{Stochastic processes in polymeric fluids},
  Springer-Verlag, 1996.

\bibitem{ralston1962runge}
A.~Ralston, \emph{Runge-kutta methods with minimum error bounds}, Mathematics
  of Computation \textbf{16} (1962), no.~80, 431--437.

\bibitem{randall2006methods}
G.~C. Randall, K.~M. Schultz, and P.~S. Doyle, \emph{Methods to
  electrophoretically stretch {DNA}: microcontractions, gels, and hybrid
  gel-microcontraction devices}, Lab on a Chip \textbf{6} (2006), no.~4,
  516--525.

\bibitem{ReVaLiHaRuPe2001}
P.~Reimann, C.~Van~den Broeck, H.~Linke, P.~H\"{a}nggi, J.M. Rubi, and
  A.~P\'{e}rez-Madrid, \emph{Giant acceleration of free diffusion by use of
  tilted periodic potentials}, Phys Rev Lett \textbf{87} (2001), 010602.

\bibitem{ReFrGa1992}
A.~Rey, J.~J. Freire, and J.~Garc'a de~la Torre, \emph{{B}rownian dynamics
  simulation of flexible polymer chains with excluded volume and hydrodynamic
  interactions. a comparison with {M}onte {C}arlo and theoretical results},
  Polymer \textbf{33} (1992), 3477--3481.

\bibitem{RoTw1996B}
G.~O. Roberts and R.~L. Tweedie, \emph{Exponential convergence of {L}angevin
  distributions and their discrete approximations}, Bernoulli \textbf{2}
  (1996), 341--363.

\bibitem{RoTw1996A}
\bysame, \emph{Geometric convergence and central limit theorems for
  multidimensional {H}astings and {M}etropolis algorithms}, Biometrika
  \textbf{1} (1996), 95--110.

\bibitem{RoDoFr1978}
P.~J. Rossky, J.~D. Doll, and H.~L. Friedman, \emph{{B}rownian dynamics as
  smart {M}onte {C}arlo simulation}, J Chem Phys \textbf{69} (1978), 4628.

\bibitem{RoPr1969}
J.~Rotne and S.~Prager, \emph{Variational treatment of hydrodynamic interaction
  in polymers}, J Chem Phys \textbf{50} (1969), 4831.

\bibitem{Sa2014}
J.~M. Sanz-Serna, \emph{Markov chain monte carlo and numerical differential
  equations}, Current Challenges in Stability Issues for Numerical Differential
  Equations (L.~Dieci and N.~Guglielmi, eds.), vol. 2082, Springer, 2014,
  pp.~39 -- 88.

\bibitem{ScBaShCh2003}
C.~M.~S. Schroeder, H.~B. Babcock, E.~S.~G. Shaqfeh, and S.~Chu,
  \emph{Observation of polymer conformation hysteresis in extensional flow},
  Science \textbf{301} (2003), 1515.

\bibitem{ScShCh2004}
C.~M.~S. Schroeder, E.~S.~G. Shaqfeh, and S.~Chu, \emph{The effect of
  hydrodynamic interactions on {DNA} dynamics in extensional flow: simulation
  and single molecule experiment}, Macromolecules \textbf{37} (2004), 9242.

\bibitem{Sc1999}
C.~Sch\"{u}tte, \emph{Conformational dynamics: Modeling, theory, algorithm, and
  application to biomolecules}, Habilitation, Free University Berlin, 1999.

\bibitem{Sh2005}
E.~S.~G. Shaqfeh, \emph{The dynamics of single-molecule {DNA} in flow}, J
  Non-Newton Fluid \textbf{130} (2005), 1.

\bibitem{SoKhWoHuSh2002}
M.~Somasi, B.~Khomami, N.~J. Woo, J.~S. Hur, and S.~G. Shaqfeh,
  \emph{{B}rownian dynamics simulations of bead-rod and bead-spring chains:
  numerical algorithms and coarse-graining issues}, J Non-Newtonian Fluid Mech
  \textbf{108} (2002), 227--255.

\bibitem{St1958}
R.~L. Stratonovich, \emph{Synchronization of a self-excited oscillator in the
  presence of noise}, Radiotekh~Electron \textbf{3} (1958), 497.

\bibitem{Ta2002}
D.~Talay, \emph{Stochastic {H}amiltonian systems: Exponential convergence to
  the invariant measure, and discretization by the implicit {E}uler scheme},
  Markov Processes and Related Fields \textbf{8} (2002), 1--36.

\bibitem{Ti1994}
L.~Tierney, \emph{{M}arkov chains for exploring posterior distributions}, Ann
  Stat \textbf{22} (1994), 1701--1728.

\bibitem{tierney1998note}
\bysame, \emph{A note on metropolis-hastings kernels for general state spaces},
  Ann Appl Probab (1998), 1--9.

\bibitem{trahan2009simulation}
D.~W. Trahan and P.~S. Doyle, \emph{Simulation of electrophoretic stretching of
  {DNA} in a microcontraction using an obstacle array for conformational
  preconditioning}, Biomicrofluidics \textbf{3} (2009), 012803.

\bibitem{weinan2011principles}
E~Weinan, \emph{Principles of multiscale modeling}, Cambridge University Press,
  2011.

\bibitem{WoShKh2004a}
N.~Woo, E.~S.~G. Shaqfeh, and B.~Khomami, \emph{The effect of confinement on
  the dynamics and rheology of dilute {DNA} solutions. i. entropic spring force
  under confinement and numerical algorithm}, J Rheol \textbf{48} (2004), 281.

\bibitem{WoShKh2004b}
\bysame, \emph{The effect of confinement on the dynamics and rheology of dilute
  {DNA} solutions. ii. effective rheology and single chain dynamics}, J Rheol
  \textbf{48} (2004), 299.

\bibitem{ZhdeGr2012}
Y.~Zhang, J.~J. de~Pablo, and M.~D. Graham, \emph{An immersed boundary method
  for {B}rownian dynamics simulation of polymers in complex geometries:
  Application to {DNA} flowing through a nano slit with embedded nanopits}, J
  Chem Phys \textbf{136} (2012), 014901.

\bibitem{ZhDoWeAlGrde2009}
Y.~Zhang, A.~Donev, T.~Weisgraber, B.~J. Alder, M.~D. Graham, and J.~J.
  de~Pablo, \emph{Tethered {DNA} dynamics in shear flow}, J Chem Phys
  \textbf{130} (2009), 234902.

\end{thebibliography}


\end{document}